\documentclass[11pt,          
english, 
]{article}

\usepackage[english,german,strings]{babel}     
\usepackage{amsmath}           
\usepackage[utf8]{inputenc}    
\usepackage[T1]{fontenc}      
\usepackage[sort]{cite}        
\usepackage{fancyhdr}          
\usepackage[a4paper]{geometry} 
\usepackage{xspace}            
\usepackage{tikz}             
\usetikzlibrary{arrows,positioning}

\usepackage{nchairx}
\usepackage[expansion=false   
]{microtype}       
\usepackage[%
final=true,        
pdfpagelabels,	
hypertexnames=false 
]{hyperref}

\makeatletter
\newcommand{\chairxauthorbibfont}{\textsc}
\newcommand{\chairxtitlebibfont}{\textit}
\newcommand{\chairxseriesbibfont}{\textit}

\newcommand{\bibnote}[2]{\nocite{#1}\@namedef{#1chairxnote}{#2}}
\makeatother

\newcommand{\Stiefel}[2]{\mathrm{St}_{{#1}, {#2}}}
\newcommand{\Gr}[2]{\mathrm{Gr}_{{#1}, {#2}}}

\newcommand{\symMat}[1]{\field{R}^{{#1} \times {#1}}_{\mathrm{sym}}}
\newcommand{\matR}[2]{\field{R}^{{#1} \times {#2}}}

\newcommand{\Stab}{\mathrm{Stab}}
\newcommand{\e}{\mathrm{e}}
\newcommand{\nablaCan}{\nabla^{\mathrm{can1}}}
\newcommand{\nablaCanSecond}{\nabla^{\mathrm{can2}}}

\newcommand{\nablaAlpha}{\nabla^{\alpha}}

\newcommand{\Qspace}{Q}

\newcommand{\QspaceLift}{\overline{Q}}

\newcommand{\prQLift}{\overline{\pi}}

\newcommand{\prQ}{\pi}

\newcommand{\liniso}[1]{{{#1}}}

\newcommand{\liegroupsub}{\mathrm{G}}
\newcommand{\liealgsub}{\mathfrak{g}}

\newcommand{\liegroup}[1]{\mathrm{{#1}}}

\newcommand{\keywords}[1]{\noindent \textbf{Keywords: {#1}}} 
\newcommand{\mathematicssubjectlassification}[1]{\noindent\textbf{2020 Mathematics Subject Classification: {#1}}}

\begin{document}
	
	\selectlanguage{english}

	\title{Rolling Reductive Homogeneous Spaces}
	
	\author{Markus Schlarb}
	\date{Department of Mathematics, \\
		Julius-Maximilians-Universit{\"a}t W{\"u}rzburg, \\
		Germany \\ \medskip
		\texttt{markus.schlarb@mathematik.uni-wuerzburg.de}\\~\\
		\today}
	
	\maketitle

\begin{abstract}
	Rollings of reductive homogeneous spaces are investigated.
	More precisely,
	for a reductive homogeneous space $G / H$ 
	with reductive decomposition
	$\liealg{g} = \liealg{h} \oplus \liealg{m}$,
	we consider rollings of $\liealg{m}$ over $G / H$
	without slip and without twist, where $G / H$ is equipped
	with an invariant covariant derivative.
	To this end, an intrinsic point of view is taken, meaning that
	a rolling is a curve in the configuration space $Q$
	which is tangent to a certain distribution.
	By considering
	a $H$-principal fiber bundle
	$\prQLift \colon \QspaceLift \to \Qspace$
	over the configuration space
	equipped with a suitable principal connection,
	rollings of $\liealg{m}$ over $G / H$
	can be expressed
	in terms of horizontally lifted curves on $\QspaceLift$.
	The total space of $\prQLift \colon \QspaceLift \to \Qspace$
	is a product of Lie groups.
	In particular, for a given control curve,
	this point of view allows for
	characterizing rollings of $\liealg{m}$ over $G / H$ 
	as solutions of an explicit, time-variant
	ordinary differential equation(ODE) on $\QspaceLift$,
	the so-called kinematic equation.
	An explicit solution for the 
	associated initial value problem 
	is obtained for rollings with respect to the canonical invariant
	covariant derivative of first and second kind if the
	development curve in $G / H$ is the projection
	of a one-parameter
	subgroup in $G$.
	Lie groups and Stiefel
	manifolds are discussed as examples.
\end{abstract}


\medskip

\mathematicssubjectlassification{%
 	53A17, 
 	53B05, 
	53C30, 
	58A30 
}


\section{Introduction}\label{sec:introduction}

Meanwhile, there is a vast literature on rolling manifolds without slip
and without twist.
First, we mention some works, where concrete expressions for
extrinsic rollings of certain
submanifolds of (pseudo-)Euclidean vector spaces
over their affine tangent spaces are derived.
Using the definition from~\cite[Ap.~B]{sharpe:1997} as starting point, 
extrinsic rollings of spheres $\mathrm{S}^n \subseteq \field{R}^{n + 1}$,
real Grassmann manifolds $\Gr{n}{k} \subseteq \symMat{n}$
and special orthogonal groups $\liegroup{SO}(n) \subseteq \matR{n}{n}$
over their affine tangent spaces are studied in~\cite{hueper.leite:2007}.
In a similar context, the Stiefel manifold
$\Stiefel{n}{k} \subseteq \matR{n}{k}$,
endowed with the Euclidean
metric, is investigated in~\cite{hueper.kleinsteuber.leite:2008}
while rollings of pseudo-orthogonal groups are considered
in~\cite{crouch.leite:2012}.
For these works, the need to solve interpolation problems on
these submanifolds in various applications seems to serve as
a motivation.
Indeed, the rolling and unwrapping technique from~\cite{hueper.leite:2007}, 
see also the more recent work~\cite{hueper.krakowski.leite:2020},
is a method to compute a $\Continuous^2$-curve
connecting a finite number of given
points on the manifolds $\mathrm{S}^n$, $\Gr{n}{k}$
and $\liegroup{SO}(n)$,
where the velocities at the initial and final
point are prescribed.
This algorithm relies on having an 
explicit expression for
the rolling of the manifold over its affine tangent
space along a curve joining the initial point with the final point. 

Beside these works, there is
the paper~\cite{molina.grong.markina.leite:2012},
where a notion of 
intrinsic rolling of an oriented Riemannian manifold $M$ over another
oriented Riemannian manifold $\widehat{M}$ is introduced
assuming $\dim(M) = \dim(\widehat{M})$.
In~\cite{markina.leite:2016}, this notion of intrinsic rolling is
generalized to
pseudo-Riemannian manifolds.
A further generalization can be found
in~\cite[Sec. 7]{grong:2012}
and~\cite[p. 35]{kokkonen:2012},
where the Levi-Civita covariant derivatives coming from the
pseudo-Riemannian metrics on
$M$ and $\widehat{M}$ are replaced by arbitrary covariant derivatives
on $M$ and $\widehat{M}$, respectively.

In this text, we investigate the following situation.
Let $G$ be a Lie group and $H \subseteq G$ a closed subgroup
such that $G / H$ is a reductive homogeneous space
with a fixed reductive decomposition $\liealg{g} = \liealg{h} \oplus \liealg{m}$.
Then $G / H$ can be equipped with
an invariant covariant derivative corresponding
to an invariant affine connection from~\cite{nomizu:1954}.
Motivated by the study of rollings of (pseudo-Riemannian) symmetric
spaces over flat spaces in~\cite{jurdjevic.markina.leite:2023},
we consider rollings of $\liealg{m}$ over $G / H$.
Here we generalize the above mentioned definition proposed
in~\cite{grong:2012} and~\cite[p. 35]{kokkonen:2012}
slightly in order take additional structures
of the involved manifolds into account.
In particular, this definition allows for considering rollings of
\emph{not} necessarily oriented
manifolds.

Moreover, if one is interested in getting rather simple formulas
describing the rollings, it might be convenient to consider rollings of
$\liealg{m}$ over $G / H$ with respect to the canonical
covariant derivative of first or second kind on $G / H$.
These covariant derivatives
can be defined independently of a pseudo-Riemannian metric
although they are in some sense similar
the Levi-Civita covariant derivatives on naturally reductive
homogeneous spaces or pseudo-Riemannian symmetric spaces,
respectively.

We now give an overview of this text.
In Section~\ref{sec:notation_and_terminology},
we start with introducing some notations and recalling some
definitions and well-known facts related to Lie groups
and principal fiber bundles.
Moreover, we recall some facts on reductive homogeneous spaces
with an emphasize on invariant covariant derivatives.

In Section~\ref{sec:rolling_general},
we briefly recall the notion of rolling intrinsically
a manifold $M$ over another manifold $\widehat{M}$
of equal dimension from the literature.
More precisely, as already announced above, a slightly
generalized definition of intrinsic rolling is introduced.

As preparation to determine the configuration space for
the intrinsic rollings considered in
Section~\ref{sec:intrinsic_rolling_reductive_space},
an explicit description of the frame bundle of a reductive homogeneous
space $G / H$ is needed.
Therefore frame bundles of reductive homogeneous spaces are investigated in Section~\ref{sec:frame_bundles}.
Here we first consider a more general situation.
The frame bundle of a vector bundle 
associated to a $H$-principal fiber bundle $P \to M$
is identified with an other fiber bundle associated to $P \to M$.
Afterwards, reductive homogeneous spaces are treated as a special case.

In Section~\ref{sec:intrinsic_rolling_reductive_space},
we turn our attention to rollings of a reductive
homogeneous space $G /H$
with reductive decomposition $\liealg{g} = \liealg{h} \oplus \liealg{m}$.
We consider the intrinsic rolling of 
$\liealg{m}$ over $G / H$ with respect to an invariant covariant
derivative $\nablaAlpha$.
To this end, the
configuration space $\Qspace \to \liealg{m} \times G / H$
is investigated
in detail.
Here we
determine a $H$-principal fiber bundle
$\prQLift \colon \QspaceLift \to \Qspace$ over $\Qspace$
which is equipped with a suitable
principal connection.
Its total space is given by
$\QspaceLift
= \liealg{m} \times G \times \liegroupsub(\liealg{m})$,
where
$\liegroupsub(\liealg{m}) \subseteq \liegroup{GL}(\liealg{m})$
is a closed subgroup,
i.e. the manifold $\QspaceLift$ is a
product of Lie groups.

For a fixed invariant covariant derivative
$\nablaAlpha$ on $G / H$
defined by an $\Ad(H)$-invariant billinear map
$\alpha \colon \liealg{m} \times \liealg{m} \to \liealg{m}$,
we determine a distribution $\overline{D^{\alpha}}$
on $\QspaceLift$ that projects
to a distribution $D^{\alpha}$ on $\Qspace$ with the following property.
A curve $q \colon I \to \Qspace$ is horizontal with respect to $D^{\alpha}$
iff it is a rolling of $\liealg{m}$ over $G / H$ with respect
to $\nablaAlpha$.
Moreover, horizontal lifts of curves on $\Qspace$
with respect to the principal
connection on $\prQLift \colon \QspaceLift \to \Qspace$ mentioned above
are horizontal with respect
to $\overline{D^{\alpha}}$ iff they are horizontal with respect to $D^{\alpha}$.
In particular, this fact allows for characterizing rollings of
$\liealg{m}$ over $G / H$ in terms
of an ODE on $\QspaceLift$.
More precisely,
for a prescribed control curve $u \colon I \to \liealg{m}$,
we obtain an
explicit, time-variant ODE on
$\QspaceLift = \liealg{m} \times G \times \liegroupsub(\liealg{m})$
whose solutions projected to $\Qspace$ are
rollings of $\liealg{m}$ over $G / H$
with respect to $\nablaAlpha$.
This ODE can be seen as a generalization of the
kinematic equation for rollings of oriented pseudo-Riemannian
symmetric spaces
over flat spaces from~\cite[Sec. 4.2]{jurdjevic.markina.leite:2023}.

In Subsection~\ref{subsec:rolling_along_special_curves},
we turn our attention to rollings
of $\liealg{m}$ over $G / H$ with respect to the canonical 
covariant derivative of first and second kind
such that the development curve is of the form
$I \ni t \mapsto \pr(\exp(t \xi)) \in G / H$ with
some $\xi \in \liealg{g}$,
i.e. a projection of a
\emph{not} necessarily horizontal one-parameter subgroup in $G$.
For this special case, an explicit solution of the kinematic equation
is obtained.

We end this text by discussing intrinsic rollings of Lie
groups and Stiefel manifolds as examples.

\section{Notations, Terminology and Background}
\label{sec:notation_and_terminology}

In this section,
we introduce the notation and terminology that is used throughout
this text.
Moreover, some facts concerning Lie groups and principal fiber bundles
are recalled. 
We end this section by discussing reductive homogeneous spaces
with an emphasize on invariant covariant derivatives.

\subsection{Notations and Terminology}

We start with introducing some notations and terminology
concerning differential geometry.
This subsection is an extended version of~\cite[Sec. 2]{schlarb:2023}
mostly copied word by word.

\begin{notation}
	Throughout this text we follow the convention
	in~\cite[Chap. 2]{oneill:1983}.
	A scalar product is defined as a non-degenerated symmetric bilinear form.
	An inner product is a positive definite symmetric bilinear form.
\end{notation}

Next we introduce some notations concerning differential geometry.
Let $M$ be a smooth (finite-dimensional) manifold.
We denote by $T M$ and $T^* M$ the tangent and cotangent bundle
of $M$, respectively.
A smooth vector subbundle $D$ of the tangent bundle $T M$ 
is called a regular distribution on $M$.
For a smooth map $f \colon M \to N$ between manifolds $M$ and $N$,
the tangent map of $f$ is denoted by
$T f \colon TM \to TN$.	
We write $\Cinfty(M)$ for the algebra of smooth real-valued functions on $M$.

Let $E \to M$ be a vector bundle over $M$ with typical fiber $V$.
The smooth sections of $E$ are denoted by $\Secinfty(E)$. 
We write $\End(E) \cong E^* \tensor E$ for the endomorphism bundle of $E$.
Moreover, we denote by $E^{\tensor k}$, $\Sym^k E$ and $\Anti^k E$
the $k$-th tensor power, the $k$-th symmetrized tensor power and the
$k$-th anti-symmetrized tensor power of $E$.
If $T \in \Secinfty\big((T^* M)^{\tensor k} \tensor (T M)^{\tensor \ell} \big)$
is a tensor field on $M$ and $X \in \Secinfty(T M)$ is a vector field,
$\Lie_X T$ denotes the Lie derivative.
If $x \colon N  \to \field{R}$ is a smooth function, 
we write $f^* x = x \circ f \colon M \to \field{R}$
for its pull-back by $f \colon M \to N$.
More generally, if $\omega \in \Secinfty\big(\Anti^k(T^* N)\big) \tensor V$
is a differential form taking values in a finite dimensional
$\field{R}$-vector space $V$, its pull-back by $f$ is denoted
by $f^* \omega$.
Next assume that $f \colon M \to N$ is a local diffeomorphism.
Then the pull-back of the tensor field
$T \in 
\Secinfty\big((T^* N)^{\tensor k} \tensor (T N)^{\tensor \ell} \big)$
by $f$ is denoted by $f^* T$, as well.

We now consider a fiber bundle $\pr \colon P \to M$ over $M$.
Its vertical bundle
is denoted by $\Ver(P) = \ker(T \pr) \subseteq T P$.
We write
$\Hor(P) \subseteq TP$ for a horizontal bundle,
i.e. a subbundle of $T P$ fulfilling
$\Ver(P) \oplus \Hor(P) = T P$.
If $\pr_P \colon P \to M$ and $\pr_Q \colon 
Q \to M$ are fiber bundles over
the same manifold $M$
with typical fiber $F_P$ and $F_Q$, respectively,
their fiber product
is denoted by $P \oplus Q \to M$.
It is the fiber bundle over $M$
given by
\begin{equation}
	P \oplus Q
	=
	\{ (p, q) \in P \times Q \mid \pr_P(p) = \pr_Q(q) \} 
\end{equation}
with typical fiber $F_P \times F_Q$.

Next let $S_1 \times \cdots \times S_k$ be a product of sets
and let $i \in \{1, \ldots, k\}$.
Then we denote by
\begin{equation}
	\pr_i \colon S_1 \times \cdots \times S_k \to S_i
\end{equation}
the projection onto the $i$-th factor.

We now recall a well-known fact on
surjective submersions.
This is the next lemma, see e.g.~\cite[Thm. 4.29]{lee:2013},
which is used
frequently without referencing it explicitly.

\begin{lemma}
	Let $\pr \colon P \to M$ be a surjective submersion
	and
	let $N$ be some manifold.
	Let $f \colon M  \to N$ be a map.
	Then $f$ is smooth iff $f \circ \pr \colon P \to N$ is smooth,
	i.e. if the diagram
	\begin{equation}
		\begin{tikzpicture}[node distance= 5.0cm, auto]
			\node (P) {$P$};
			\node (M) [below=1.5cm of P] {$M$};
			\node (N) [right=3.0cm of M] {$N$};
			\draw[->] (P) to node [left] {$\pr$} (M);	
			\draw[->] (M) to node  [below]{$f$} (N);
			\draw[->] (P) to node [pos = 0.5] {$f \circ \pr$} (N);			
		\end{tikzpicture}  
	\end{equation}
	commutes in the category of smooth manifolds.
\end{lemma}
Concerning the regularity of curves on manifolds,
we use the following convention.

\begin{notation}
	Whenever $c \colon I \to M$ denotes a curve in a manifold $M$
	defined on an interval $I \subseteq \field{R}$,
	we assume for simplicity
	that $c$ is smooth if not indicated otherwise.
	If $I$ is not open, 
	we assume that $c$ can be extended to smooth
	curve defined on
	an open interval $J \subseteq \field{R}$ containing $I$.
	Moreover, we implicitly assume that $0$ is contained in $I$
	if we write $0 \in I$.
	Nevertheless, many results
	can be generalized by requiring less regularity.
\end{notation}

\begin{notation}
	If not indicated otherwise, we use Einstein summation
	convention.
\end{notation}

\subsection{Lie groups}

Copying and extending~\cite[Sec. 3.1]{schlarb:2023},
we now introduce some notations and well-known facts
concerning Lie groups
and Lie algebras.

Let $G$ be a Lie group and denote its Lie algebra by $\liealg{g}$.
The identity of $G$ is usually denoted by $e$.
The left translation by an element $g \in G$ is denoted by
\begin{equation}
	\ell_g \colon G \to G,
	\quad 
	h \mapsto \ell_g(h) = g h 
\end{equation}
and we write 
\begin{equation}
	r_g \colon G \to G, 
	\quad
	h \mapsto r_g(h) = h g
\end{equation}
for the right translation by $g \in G$.
The conjugation by an element $g \in G$ is given by
\begin{equation}
	\Conj_{g} \colon G \to G,
	\quad
	h \mapsto \Conj_g(h) 
	=
	(\ell_g \circ r_{g^{-1}})(h)
	=
	(r_{g^{-1}} \circ \ell_g)(h)
	=
	g h g^{-1} 
\end{equation}
and the adjoint representation of $G$ is defined as
\begin{equation}
	\Ad \colon G \to \liegroup{GL}(\liealg{g}),
	\quad
	g \mapsto \Ad_g 
	= \big(\xi \mapsto \Ad_g(\xi) = T_e \Conj_g \xi \big) .
\end{equation}
Moreover, we denote the adjoint representation of $\liealg{g}$ by
\begin{equation}
	\ad \colon \liealg{g} \to \liealg{gl}(\liealg{g}),
	\quad
	\xi \mapsto \big(\eta \mapsto \ad_{\xi}(\eta) = [\xi, \eta] \big) .
\end{equation}
Next we recall~\cite[Def. 19.7]{gallier.quaintance:2020}.
A vector field $X \in \Secinfty(T G)$ is called left-invariant
or right-invariant
if for all $g, k \in G$
\begin{equation}
	T_k \ell_g X(k) = X (\ell_g(k))
	\quad
	\text{ or } \quad
	T_k r_g X(k) = X(r_g(k)),
\end{equation}
respectively,
holds.
For $\xi \in \liealg{g}$, we denote by $\xi^L \in \Secinfty(T G)$
and $\xi^R \in \Secinfty(T G)$
the corresponding left and right-invariant vector fields, respectively,
which are given by 
\begin{equation}
	\xi^L(g) = T_e \ell_g \xi 
	\quad \text{ and } \quad
	\xi^R(g) = T_e r_g \xi,
	\quad g \in G .
\end{equation} 	
The exponential map of the Lie group $G$ is denoted by 
\begin{equation}
	\exp \colon \liealg{g} \to G .
\end{equation}
One has for $\xi \in \liealg{g}$ and $t \in \field{R}$
\begin{equation}
	\label{equation:exponential_map_lie_group}
	\tfrac{\D}{\D t} \exp(t \xi)
	=
	T_e \ell_{\exp(t \xi)} \xi
	=
	T_e r_{\exp(t \xi)} \xi
\end{equation}
by the proof of~\cite[Prop. 19.5]{gallier.quaintance:2020}.

Next we recall
that the tangent map of the group multiplication 
$m \colon G \times G \ni (g, h) \mapsto g h \in G$ is given by
\begin{equation}
	\label{equation:tangentmap_group_multiplication}
	T_{(g, h)} m (v_g, w_h) = T_g r_h v_g + T_h \ell_g w_h
\end{equation}
for all $(g, h) \in G \times G$ and $(v_g, v_h) \in T_g G \times T_h G$,
see e.g.~\cite[Lem. 4.2]{michor:2008}.
The tangent map of the inversion $\inv \colon G \ni g \mapsto \inv(g) = g^{-1} \in G$ reads
\begin{equation}
	\label{equation:tangentmap_group_inversion}
	T_g \inv v_g = - (T_e \ell_{g^{-1}} \circ T_{g} r_{g^{-1}}) v_g
\end{equation}
for all $g \in G$ and $v_g \in T_g G$,
see e.g.~\cite[Cor. 4.3]{michor:2008}.

We now introduce the notation for some Lie groups that play
a crucial role in this text.
\begin{notation}
	Let $V$ be a finite dimensional $\field{R}$-vector space.
	We write $\liegroup{GL}(V)$ for the general
	linear group of $V$.
	If $V$ is a pseudo-Euclidean vector space,
	i.e. $V$ is endowed with a scalar product
	$\langle \cdot, \cdot \rangle \colon V \times V \to \field{R}$,
	we denote the
	corresponding pseudo-orthogonal group by 
	$\liegroup{O}(V, \langle \cdot, \cdot \rangle)$.
	Moreover, we often write
	$\liegroup{O}(V) = \liegroup{O}(V, \langle \cdot, \cdot \rangle)$
	for short.
	Similarly, the special (pseudo-)orthogonal group is denoted by
	$\liegroup{SO}(V) = \liegroup{SO}(V, \langle \cdot, \cdot \rangle)$.
	More generally, a closed subgroup of $\liegroup{GL}(V)$,
	which is not further specified, is
	often denoted by
	$\liegroupsub(V)$ and we write
	$\liealgsub(V) \subseteq \liealg{gl}(V)$ for the
	corresponding Lie algebra.
	Sometimes, the exponential map of $\liegroupsub(V)$ is denoted by
	\begin{equation}
		\liealgsub(V) \to \liegroupsub(V), 
		\quad
		\xi \mapsto \e^\xi = \sum_{k = 0}^{\infty} \tfrac{1}{k !} \xi^k .
	\end{equation}
	In the sequel, it is often convenient to denote the evaluation of
	$A \in \liegroup{GL}(V)$ at $v \in V$
	by $A v$ instead of writing $A(v)$.
\end{notation}

\subsection{Principal Fiber Bundles}

Next we recall
some well-known facts on principal fiber bundles and
introduce some notations.
For general facts on principal fiber bundles we refer
to~\cite[Sec. 18-19]{michor:2008}
and~\cite[Sec. 1.1-1.3]{rudolph.schmidt:2017}.
\begin{notation}
	\label{notation:principal_bundle_and_associated_bundle}
	Let $P \to M$ be a $H$-principal fiber bundle over $M$
	and let $\liealg{h}$ be the Lie algebra of $H$.
	The principal action is usually denoted by
	\begin{equation}
		\racts \colon P \times H \to P, 
		\quad (p, h) \mapsto p \racts h 
	\end{equation}
	and we denote for fixed $h \in H$ by
	$(\cdot \acts h) \colon P \ni p \mapsto p \racts h \in P$
	the induced diffeomorphism.
\end{notation}

Next, let $\eta \in \liealg{h}$. Then $\eta_P\in \Secinfty(T P)$ denotes the fundamental vector field associated to the principal action.
For $p \in P$, it is given by
\begin{equation}
	\eta_P(p) = \tfrac{\D}{\D t} (p \racts \exp(t \eta) ) \at{t = 0} .
\end{equation}
The vertical bundle $\Ver(P) = \ker( T  \pr) \subseteq T P$
of $P \to M$ is fiber-wise given by
\begin{equation}
	\label{equation:vertical_bundle_fiber_bundle_fiber_wise_description}
	\Ver(P)_p
	=
	\big\{
	\eta_P(p) \mid \eta \in \liealg{h}
	\big\}
	=
	\big\{  \tfrac{\D}{\D t} (p \racts \exp(t \eta)) \at{t = 0} \mid \eta \in \liealg{h} \big\}
	\subseteq
	T_p P,
	\quad
	p \in P
\end{equation}
and the map
\begin{equation}
	\label{equation:vertical_bundle_principal_bundle_trivialization}
	P \times \liealg{h} \to \Ver(P),
	\quad
	(p, \eta) \mapsto T_e (p \racts \cdot ) \eta
	= \tfrac{\D}{\D t} (p \racts \exp(t \eta)) \at{t = 0}
\end{equation}
is an isomorphism of vector bundles covering $\id_P$,
see e.g.~\cite[Lem. 1.3.1]{rudolph.schmidt:2017}
or~\cite[Sec. 18.18]{michor:2008}.

Recall that a complement of $\Ver(P)$, i.e. a
subbundle $\Hor(P) \subseteq T P$ fulfilling
$\Hor(P) \oplus \Ver(P) = TP$
is called horizontal bundle.
It is well-known that
such a complement defines a unique connection on $P$,
i.e. an endomorphism $\mathcal{P} \in \Secinfty\big(\End(T P)\big)$
such that
$\mathcal{P}^2 = \mathcal{P}$ and $\image(\mathcal{P})
= \Ver(P)$
as well as $\image(\mathcal{P}) = \Hor(P)$ holds.
This fact can be regarded as a consequence
of~\cite[Sec. 17.3]{michor:2008}.
Moreover, such a connection on
the principal fiber bundle $P \to M$ is in
one-to-one correspondence with a
$\liealg{h}$-valued one-form
$\omega \in \Secinfty(T^* P) \tensor \liealg{h}$
via
\begin{equation}
	\label{equation:principal_connection_connection_one_form_correspondence}
	\omega\at{p}(v_p) = \big( T_e (p \racts \cdot) \big)^{-1} \mathcal{P} \at{p}(v_p),
	\quad
	p \in P, 
	\quad 
	v_p \in T_p P,
\end{equation}
see e.g.~\cite[Sec. 19.1, Eq. (1)]{michor:2008}, due to
the identification
$\Ver(P) \cong P \times \liealg{h}$
as vector bundles via the
isomorphism~\eqref{equation:vertical_bundle_principal_bundle_trivialization}.
A connection $\mathcal{P} \in \Secinfty\big(\End(T P)\big)$
is called principal connection if
\begin{equation}
	T_p (\cdot \racts h)  \big(\mathcal{P}\at{p}(v_p) \big)
	=
	\mathcal{P}\at{p \racts h}\big( T_p (\cdot \racts h)  v_p \big),
	\quad
	p \in P, 
	\quad 
	v_p \in T_p P
\end{equation}
holds for all $h \in H$ or equivalently
$\big( \cdot \racts h \big)^* \mathcal{P} = \mathcal{P}$
is satisfied for all $h \in H$, see e.g.~\cite[Sec. 19.1]{michor:2008}.
We now recall from~\cite[Sec. 19.1]{michor:2008}
how a principal connection
is related to the corresponding
connection one-form
given by~\eqref{equation:principal_connection_connection_one_form_correspondence}.
This is the next lemma.
\begin{lemma}
	\label{lemma:principal_connection_connection_one_form}
	Let $\mathcal{P} \in \Secinfty\big(\End(T P)\big)$ be a principal connection.
	Then the corresponding
	connection one-form
	$\omega \in \Secinfty(T^* P) \tensor \liealg{h}$ satisfies:
	\begin{enumerate}
		\item
		\label{item:lemma_principal_connection_connection_one_form_1}
		For each $\eta \in \liealg{h}$, one has
		$\omega\at{p}\big(\eta_P(p) \big) = \eta$ for all $p \in P$.
		\item
		\label{item:lemma_principal_connection_connection_one_form_2}
		For each $h \in H$ one has
		$\big((\cdot \racts h)^* \omega\big)\at{p}(v_p) 
		=
		\Ad_{h^{-1}}\big(\omega\at{p}(v_p)\big) $
		for all $p \in P$ and $v_p \in T_p P$.
	\end{enumerate}
	Conversely, a $\liealg{h}$-valued one-form
	$\omega \in \Secinfty(T^* P) \tensor \liealg{h}$ fulfilling
	Claim~\ref{item:lemma_principal_connection_connection_one_form_1}
	and
	Claim~\ref{item:lemma_principal_connection_connection_one_form_2}
	defines a principal connection on $P \to M$ via
	\begin{equation}
		\mathcal{P}\at{p}(v_p) 
		=
		\big( T_e (p \racts \cdot) \big) \omega\at{p}(v_p)
		=
		\big( \omega\at{p}(v_p)\big)_P(p) 
	\end{equation}
	for $p \in P$ and $v_p \in T_p P$ with the map
	$(p \racts \cdot) \colon H \ni h  \mapsto p \racts h \in P$
	for fixed $p \in P$.
\end{lemma}

Next we recall the notion of reductions of principal fiber bundles,
see e.g.~\cite[Sec. 18.6]{michor:2008}.
Let $P\to M$ be a $H$-principal fiber
bundle.
Then a $H_2$-principal fiber bundle $P_2 \to M$
is called a reduction
of $P$
if there is a morphism of Lie groups $f \colon H_2 \to H$ and
a morphism $\Psi \colon P_2 \to P$
of principal fiber bundles along $f$
covering $\id_M \colon M \to M$.
In particular,
\begin{equation}
	\Psi(p_2 \racts h_2) = \Psi(p_2) \racts f(h_2)
\end{equation}
holds for all $h_2 \in H_2$ and $p_2 \in P_2$.
In the sequel, only reductions of $H$-principal fiber
bundles along the canonical inclusion of a closed
subgroup $H_2 \subseteq H$ into a Lie group $H$
will play a role.

Furthermore, we need the notion of an associated bundle
which we recall briefly
from~\cite[Sec. 18.7]{michor:2008}.
Let $F$ be some manifold and let $\acts \colon H \times F \to F$
be a smooth action of $H$ on $F$ from the left. 
Then the corresponding associated bundle is denoted by
\begin{equation}
	\pi \colon P \times_H F \to M ,
\end{equation}
whose elements are given by
\begin{equation}
	P \times_H F 
	=
	\big\{ [p, s] \mid (p, s) \in P \times F \big\} .
\end{equation}
Here $[p, s]$ denotes the equivalence class of
$(p, s) \in P \times F$
defined by the $H$-action on $P \times F$ given by
\begin{equation}
	\label{equation:associated_bundle_action_on_product}
	(P \times F) \times H \ni \big((p, f), h \big) 
	\mapsto
	\big(p \racts h, h^{-1} \acts f \big) \in P \times F, 
\end{equation}
i.e. 
$(p, s) \sim (p^\prime, s^\prime)$ iff there exists a $h \in H$
such that 
$(p^\prime, s^\prime) = (p \racts h, h^{-1} \acts s)$ is fulfilled.
The projection $\pi \colon P \times_H F \to M$, sometimes denoted
by
$\pi_{P \times_H F} \colon P \times_H F \to M$
to refer to $P \times_H F$ explicitly,
is given by
$\pi([p, s]) = \pr_P(p)$, where $\pr_P \colon P \to M$ denotes
the projection of the principal fiber bundle.

Moreover, we often write 
\begin{equation}
	\label{equation:associated_bundle_principal_bundle_over_it}
	\overline{\pi} \colon P \times F \to (P \times F) / H = P \times_H F
\end{equation}
for the $H$-principal fiber bundle over the
associated bundle $P \times_H F$, where
the principal action is given by~\eqref{equation:associated_bundle_action_on_product}
and the map $\overline{\pi}$ in~\eqref{equation:associated_bundle_principal_bundle_over_it}
is the canonical projection, i.e. $\overline{\pi}(p, f) = [p, f]$ for $(p, f) \in P \times F$.
We also denote this projection by
$\overline{\pi}_{P \times F} \colon P \times F \to P \times_H F$
to refer to $P \times F$ explicitly.

Moreover, we will use the following identification of
the tangent bundle of an associated bundle $P \times_H F \to M$
of a $H$-principal fiber bundle
$P \to M$
\begin{equation}
	\label{equation:identifcation_tangent_bundle_associated_bundle}
	T(P \times_H F)  
	\cong
	T P \times_{T H} T F
	=
	\big\{ [v_p, v_s] \mid (v_p, v_s) \in T_{(p, s)} (P \times F), (p, s) \in P \times F \big\},
\end{equation}
see e.g.~\cite[Sec. 18.18]{michor:2008}.
Here $T P$ is considered as $T H$-principal fiber bundle
over $T M$ with principal action
$T \racts \colon TP \times T H \ni (v_p, v_h) \mapsto (T \racts)(v_p, v_h) = v_p (T \racts) v_h \in TP$,
see e.g.~\cite[Sec. 18.18]{michor:2008},
and $T H$ acts on $T F$ via the tangent map of the $H$-action
on $F$ denoted by
$T \acts \colon T H \times T F \ni (v_h, v_s) \mapsto (T \acts)(v_h, v_s) = v_h (T  \acts) v_p \in TF $.
Moreover, $[v_p, v_s]$ denotes the equivalence class defined
by $(v_p, v_s) \in T_p P \times T_s F$ for $(p, s) \in P \times F$.
Here the equivalence relation is defined by
$(v_p, v_s) \sim (v_p^{\prime}, v_s^{\prime})$
iff there exits an $v_h \in T H$ with
$(v_p^{\prime}, v_s^{\prime}) 
=
(v_p (T \racts) v_h, v_h^{-1} (T \acts) v_s)$.

Finally, we introduce some notations concerning frame bundles
of vector bundles.
We refer to~\cite[Sec. 18.11]{michor:2008} for general information
on frame bundles.
\begin{notation}
	The frame bundle of a vector bundle $E \to M$
	with typical fiber $V$
	is denoted by $\liegroup{GL}(V, E) \to M$.
	If $E$ is equipped with a not necesarrily positive definite fiber metric, 
	we denote the corresponding (pseudo-)orthogonal frame bundle by 
	$\liegroup{O}(V, E) \to M$. 
	More generally, let $\liegroup{G}(V) \subseteq \liegroup{GL}(V)$ 
	be a closed subgroup
	of the general linear group $\liegroup{GL}(V)$.
	Then a $\liegroup{G}(V)$-reduction of $\liegroup{GL}(V,E)$
	along the canonical inclusion $\liegroupsub(V) \to \liegroup{GL}(V)$
	is often denoted by $\liegroup{G}(V, E)$ if it exists.
	We write
	$\pr_{\liegroupsub(V, E)} \colon \liegroupsub(V, E) \to M$
	for the bundle projection.
\end{notation}

\subsection{Reductive Homogeneous Spaces}

In this subsection, we
recall some well-known facts on reductive homogeneous spaces
and introduce the notation that is used throughout this text.
Afterwards, invariant covariant derivatives are considered.

\subsubsection{Some General Facts}

We now recall some general facts concerning reductive homogeneous spaces
mostly by copying~\cite[Sec. 3.2]{schlarb:2023} word by word.
We refer to~\cite[Sec. 23.4]{gallier.quaintance:2020}
or~\cite[Chap. 11]{oneill:1983}
for details.

Let $G$ be a Lie group and let $\liealg{g}$ be its Lie algebra.
Moreover, let $H \subseteq G$ a closed subgroup
whose Lie algebra is denoted by
$\liealg{h} \subseteq \liealg{g}$.
We consider the homogeneous space $G / H$.
Then
\begin{equation}
	\label{equation:action_on_homogeneous_space}
	\tau \colon G \times G / H \to G / H, 
	\quad 
	(g, g^{\prime} \cdot H) \mapsto (g g^{\prime}) \cdot H
\end{equation}
is a smooth $G$-action on $G / H$ from the left,
where $g \cdot H \in G / H$ denotes the coset defined by $g \in G$.
Borrowing the notation from~\cite[p. 676]{gallier.quaintance:2020},
for fixed $g \in G$, the associated diffeomorphism is denote by 
\begin{equation}
	\tau_g \colon G / H \to G / H, 
	\quad 
	g^{\prime} \cdot H 
	\mapsto \tau_g(g^{\prime} \cdot H) 
	=
	(g g^{\prime}) \cdot H .
\end{equation}
In addition, we write
\begin{equation}
	\pr \colon G \to G /H, \quad
	g \mapsto \pr(g) = g \cdot H
\end{equation}
for the canonical projection.

Since reductive homogeneous spaces play a central role in this text,
we recall their definition
from~\cite[Def. 23.8]{gallier.quaintance:2020},
see also~\cite[Sec. 7]{nomizu:1954}
or~\cite[Def. 21, Chap. 11]{oneill:1983}.

\begin{definition}
	Let $G$ be a Lie group and $\liealg{g}$ be its Lie algebra.
	Moreover, let
	$H \subseteq G$ be a closed subgroup and denote its Lie
	algebra by $\liealg{h} \subseteq \liealg{g}$.
	Then the homogeneous space $G / H$ is called reductive 
	if there exists a subspace
	$\liealg{m} \subseteq \liealg{g}$ 
	such that $\liealg{g} = \liealg{h} \oplus \liealg{m}$ is fulfilled and
	\begin{equation}
		\Ad_h(\liealg{m}) \subseteq \liealg{m}
	\end{equation}
	holds for all $h \in H$.
\end{definition}
Following~\cite[Prop. 23.22]{gallier.quaintance:2020}, we 
recall a well-known property of the isotropy representation of
a reductive homogeneous space.
This is the next lemma.

\begin{lemma}
	\label{lemma:isotropy_representation_equivalence}
	The isotropy representation 
	of a reductive homogeneous space $G / H$
	with reductive decomposition
	$\liealg{g} = \liealg{h} \oplus \liealg{m}$
	\begin{equation}
		H \ni h 
		\mapsto T_{\pr(e)} \tau_h 
		\in \liegroup{GL}\big(T_ {\pr(e)} G / H \big)
	\end{equation}
	is equivalent to the representation
	\begin{equation}
		H \to \liegroup{GL}(\liealg{m}),
		\quad
		h \mapsto \Ad_h\at{\liealg{m}} = \big(X \mapsto \Ad_h(X) \big),
	\end{equation}
	i.e.
	\begin{equation}
		T_{\pr(e)}  \tau_h  \circ  T_e \pr\at{\liealg{m}} 
		=
		T_e \pr \circ \Ad_h\at{\liealg{m}}
	\end{equation}
	is fulfilled for all $h \in H$.
\end{lemma}

\begin{notation}
	Let $\liealg{g} = \liealg{h} \oplus \liealg{m}$ be
	a reductive decomposition of $\liealg{g}$.
	Then the projection onto $\liealg{m}$ 
	whose kernel is given by $\liealg{h}$ is denoted by
	$\pr_{\liealg{m}} \colon \liealg{g} \to \liealg{m}$.
	We write 
	$\pr_{\liealg{h}} \colon \liealg{g} \to \liealg{h}$ for the projection whose kernel is given by $\liealg{m}$.
	Moreover, we write for $\xi \in \liealg{g}$
	\begin{equation}
		\xi_{\liealg{m}} = \pr_{\liealg{m}}(\xi)
		\quad  \text{ and } \quad
		\xi_{\liealg{h}} = \pr_{\liealg{h}}(\xi) .
	\end{equation}
\end{notation}
A scalar product $\langle \cdot, \cdot \rangle \colon \liealg{m} \times \liealg{m} \to \field{R}$
is called $\Ad(H)$-invariant if
\begin{equation}
	\big\langle \Ad_h(X), \Ad_h(Y) \big\rangle
	=
	\big\langle X, Y \big\rangle
\end{equation}
holds for all $h \in H$ and $X, Y \in \liealg{m}$,
see e.g~\cite[p. 301]{oneill:1983}
or~\cite[Sec. 23.4]{gallier.quaintance:2020}
for the positive definite case.	
Reformulating and adapting~\cite[Def. 23.5]{gallier.quaintance:2020},
we call a pseudo-Riemannian metric
$\langle \! \langle \!  \cdot, \cdot \rangle \! \rangle \in \Secinfty\big( \Sym^2 T^* (G / H)\big)$
invariant if
\begin{equation}
	\langle \! \langle v_p, w_p \rangle \! \rangle_{p}
	=
	\langle \! \langle T_p \tau_g v_p , T_p \tau_g w_p \rangle \! \rangle_{\tau_g(p)},
	\quad
	p \in G / H, \quad v_p, w_p \in T_p (G  / H)
\end{equation}
holds for all $g \in G$.
In the next lemma which is taken 
from~\cite[Chap. 11, Prop. 22]{oneill:1983}, see
also~\cite[Prop. 23.22]{gallier.quaintance:2020}
for the Riemannian case,
invariant metrics on $G / H$ are related to
$\Ad(H)$-invariant scalar products on $\liealg{m}$.

\begin{lemma}
	By requiring the linear isomorphism
	$T_e \pr\at{m} \colon \liealg{m} \to T_{\pr(e)} (G / H)$
	to be an isometry,
	there is a one-to-one correspondence between
	$\Ad(H)$-invariant scalar products on $\liealg{m}$ and
	invariant pseudo-Riemannian metrics on $G / H$.
\end{lemma}
Naturally reductive homogeneous spaces are special
reductive homogeneous spaces.
We recall their definition from~\cite[Def. 23, Chap. 11]{oneill:1983}.

\begin{definition}
	\label{definition:naturally_reductive_space}
	Let $G / H$ be a reductive homogeneous space equipped with an invariant metric and denote by  $\langle \cdot, \cdot \rangle \colon \liealg{m} \times \liealg{m} \to \field{R}$
	the corresponding $\Ad(H)$-invariant scalar product on $\liealg{m}$.
	Then $G /H$ is called naturally reductive homogeneous space if 
	\begin{equation}
		\big\langle [X, Y]_{\liealg{m}}, Z \big\rangle 
		= 
		\big\langle X, [Y, Z]_{\liealg{m}} \big\rangle
	\end{equation}
	holds for all $X, Y, Z \in \liealg{m}$.
\end{definition}
The following lemma can be considered as a generalization
of~\cite[Prop. 23.29 (1)-(2)]{gallier.quaintance:2020}
to pseudo-Riemannian metrics and
Lie groups which are not necessarily connected.

\begin{lemma}
	\label{lemma:normal_naturally_reductive_is_naturally_reductive}
	Let $G$ be a Lie group and denote by $\liealg{g}$ its Lie algebra.
	Moreover, let $G$ be equipped with a bi-invariant metric
	and let
	$\langle \cdot, \cdot \rangle \colon \liealg{g} \times \liealg{g}
	\to \field{R}$
	be the corresponding $\Ad(G)$-invariant scalar product.
	Moreover, let $H \subseteq G$ be a closed subgroup such that
	its Lie algebra 
	$\liealg{h} \subseteq \liealg{g}$
	is non-degenerated with respect to $\langle \cdot, \cdot \rangle$.
	Then $G / H$ is a reductive homogeneous space with reductive 
	decomposition $\liealg{g} = \liealg{h} \oplus \liealg{m}$, 
	where $\liealg{m} = \liealg{h}^{\perp}$
	is the orthogonal complement of $\liealg{h}$
	with respect to $\langle \cdot, \cdot \rangle$.
	Moreover, if $G / H$ is equipped with the invariant metric
	corresponding to the scalar product on $\liealg{m}$ that is obtained by restricting $\langle \cdot, \cdot \rangle$ to $\liealg{m}$,
	the reductive homogeneous space $G / H$ is naturally reductive.
	\begin{proof}
		The claim can be proven analogously to the proof 
		of~\cite[Prop. 23.29 (1)-(2)]{gallier.quaintance:2020}
		by taking the assumption
		$\liealg{h} \oplus \liealg{h}^{\perp} 
		=
		\liealg{h} \oplus \liealg{m}
		= \liealg{g}$
		into account.
	\end{proof}
\end{lemma}

\begin{remark}
	\label{remark:normal_naturally_recductive}
	Inspired by the terminology
	in~\cite[Sec. 23.6, p. 710]{gallier.quaintance:2020},
	we refer to the naturally reductive spaces from
	Lemma~\ref{lemma:normal_naturally_reductive_is_naturally_reductive}
	as normal naturally reductive spaces.
\end{remark}

We now consider another special class of
reductive homogeneous spaces.
To this end, we state the following
definition which can be found in~\cite[p. 209]{helgason:1978}.

\begin{definition}
	Let $G$ be a connected Lie group and let $H$ be a closed subgroup.
	Then $(G, H)$ is called a symmetric pair if there exists a
	smooth involutive
	automorphism $\sigma \colon G \to G$, i.e.
	an automorphism of Lie groups fulfilling $\sigma^2 = \sigma$,
	such that $(H_{\sigma})_0 \subseteq H \subseteq H_{\sigma}$ holds.
	Here $H_{\sigma}$ denotes the set of fixed points of $\sigma$ and
	$(H_{\sigma})_0$ denotes the connected component of $H_{\sigma}$
	containing the identity $e \in G$.
\end{definition}

Inspired by the terminology used
in~\cite[Def. 23.13]{gallier.quaintance:2020}, we refer to the 
triple $(G, H, \sigma)$ as symmetric pair, as well, where $(G, H)$ is a symmetric pair
with respect to the involutive automorphism
$\sigma \colon G \to G$.
These symmetric pairs lead to reductive homogeneous spaces
which are called symmetric homogeneous spaces if a certain ``canonical''
reductive decomposition is chosen, see
e.g.~\cite[Sec. 14]{nomizu:1954}.
Note that the definition 
in~\cite[Sec. 14]{nomizu:1954}
does not require an invariant pseudo-Riemannian metric on $G / H$.

The next lemma,
see e.g.~\cite[Sec. 14]{nomizu:1954}
shows that a symmetric homogeneous space is a reductive homogeneous
space with respect to the so-called canonical
reductive decomposition.
Here we also refer
to~\cite[Prop. 23.33]{gallier.quaintance:2020} for a proof.

\begin{lemma}
	\label{lemma:symmetric_pair_canonical_reductive_decomposition}
	Let $(G , H, \sigma)$ be a symmetric pair
	and define the subspaces of $\liealg{g}$ by
	\begin{equation}
		\liealg{h} = \{ X \in \liealg{g} \mid T_e \sigma X = X\}
		\subseteq \liealg{g}
		\quad  \text{ and } \quad
		\liealg{m} = \{X \in \liealg{g} \mid T_e \sigma X = - X \} 
		\subseteq \liealg{g} .
	\end{equation}
	Then $\liealg{g} = \liealg{h} \oplus \liealg{m}$ is a
	reductive decomposition of $\liealg{g}$ turning $G / H$
	into a reductive homogeneous space.
	Moreover, the inclusion
	\begin{equation}
		[\liealg{m} , \liealg{m}] \subseteq \liealg{h}
	\end{equation}
	is fulfilled.
\end{lemma}
Next we define symmetric homogeneous spaces and canonical reductive
decompositions following~\cite[Sec. 14]{nomizu:1954}.

\begin{definition}
	Let $(G, H, \sigma)$ be a symmetric pair.
	Then the reductive decomposition
	$\liealg{g} = \liealg{h} \oplus \liealg{m}$
	from Lemma~\ref{lemma:symmetric_pair_canonical_reductive_decomposition}
	is called canonical reductive decomposition.
	Moreover, the reductive homogeneous space $G / H$ with
	the reductive decomposition
	$\liealg{g} = \liealg{h} \oplus \liealg{m}$ 
	is called symmetric homogeneous space.
\end{definition}
For pseudo-Riemannian symmetric spaces we state the next remark
following~\cite[Chap. 11, p. 317]{oneill:1983}, see
also~\cite[Sec. 23.8]{gallier.quaintance:2020}
for the Riemannian case.

\begin{remark}
	\label{remark:pseudo-Riemannian-symmetric-naturally-reductive}
	Let $(G, H, \sigma)$ be symmetric pair and let $G / H$ be the
	associated symmetric homogeneous space with canonical reductive
	decomposition $\liealg{g} = \liealg{h} \oplus \liealg{m}$.
	Let $G / H$ be equipped with an invariant pseudo-Riemannian metric
	and let
	$\langle \cdot, \cdot \rangle \colon 
	\liealg{m} \times \liealg{m} \to \field{R}$
	be the associated $\Ad(H)$-invariant scalar product.
	Then $G / H$ is a naturally reductive homogeneous space
	since $[\liealg{m}, \liealg{m}] \subseteq \liealg{h}$ implies
	that the condition on the scalar product
	$\langle \cdot, \cdot \rangle$ from
	Definition~\ref{definition:naturally_reductive_space}
	is always satisfied.
	In the sequel, we refer to symmetric homogeneous spaces equipped 
	with an invariant pseudo-Riemannian metric as pseudo-Riemannian
	symmetric homogeneous space
	or pseudo-Riemannian symmetric spaces, for short.
\end{remark}

\subsubsection{Reductive Decompositions and Principal Connections}

We now consider
the principal connection defined on the principal fiber bundle
$\pr \colon G \to G / H$ by means of the reductive decomposition
This subsection can be seen as a short version 
of~\cite[Sec. 3.3]{schlarb:2023} mostly copied word by word.

Let $G$ be a Lie group and $H \subseteq G$ a closed subgroup.
It is well-known that
$\pr \colon G \to G / H$ carries the structure of a $H$-principle
fiber bundle, see e.g. \cite[Sec. 18.15]{michor:2008},
where the base is the homogeneous space $G / H$.
The $H$-principal action on $G$ is denoted by 
\begin{equation}
	\label{equation:H_principal-action_on_G_mod_H}
	\racts \colon G \times H \to G, 
	\quad 
	(g, h) \mapsto g \racts h
	= r_h(g) 
	= \ell_g(h) 
	= g h,
\end{equation}
if not indicated otherwise.
We now assume that $G / H$ is a reductive homogeneous space
and the reductive decomposition 
$\liealg{g} = \liealg{h} \oplus \liealg{m}$
is fixed. It is well-known that this reductive decomposition
can be used to obtain a
principal connection on $\pr \colon G \to G / H$,
see~\cite[Thm. 11.1]{kobayashi.nomizu:1963a}.
Before we consider this principal connection on
the principal fiber bundle $G \to G / H$,
we comment on its vertical bundle.
For fixed $g \in G$ the vertical bundle of $\pr \colon G \to G / H$ 
is given by
\begin{equation}
	\label{equation:vertical_bundle_reductive_space_description}
	\Ver(G)_g 
	=
	\big\{ \tfrac{\D}{\D t} \big( g \racts \exp(t \eta) \big) \at{t = 0} \mid \eta \in \liealg{h} \big\}
	= 
	(T_e \ell_g )\liealg{h} 
\end{equation}
according
to~\eqref{equation:vertical_bundle_fiber_bundle_fiber_wise_description}.
The next proposition describes the principal connection on $\pr \colon G \to G / H$
defined by the reductive decomposition
$\liealg{g} = \liealg{h} \oplus \liealg{m}$
mentioned above.

\begin{proposition}
	\label{proposition:principal_connection_reductive_homogeneous_space}
	Consider $\pr \colon G \to G / H$ as a $H$-principal fiber bundle,
	where $G / H$ is a reductive homogeneous space with reductive
	decomposition $\liealg{g} = \liealg{h} \oplus \liealg{m}$
	and define $\Hor(G) \subseteq T G$ fiber-wise by 
	\begin{equation}
		\Hor(G)_g = (T_e \ell_g) \liealg{m}, \quad g \in G.
	\end{equation}
	Then $\Hor(G)$ is a subbundle of $TG$ defining a horizontal bundle
	on $T G$, i.e. a complement of the vertical bundle
	$\Ver(G) = \ker (T \pr) \subseteq TG$
	which yields a principal connection on $\pr \colon G \to G / H$. 
	This principal connection $\mathcal{P} \in \Secinfty\big(\End(T G) \big)$ 
	corresponding to $\Hor(G)$ is given by
	\begin{equation}
		\label{equation:principal_connection_G_mod_H_defined_by_reductive_split}
		\mathcal{P}\at{g}(v_g) = T_e \ell_g \circ \pr_{\liealg{h}} \circ (T_e \ell_g)^{-1} v_g,
		\quad 
		g \in G, \quad v_g \in T_g G .
	\end{equation}
	The corresponding connection one-form $\omega \in \Secinfty(T^* G) \tensor \liealg{h}$ reads 
	\begin{equation}
		\omega\at{g}(v_g) = \pr_{\liealg{h}} \circ (T_e \ell_g)^{-1} v_g 
	\end{equation}
	for $g \in G$ and $v_g \in T_g G$.
\end{proposition}

\subsubsection{Invariant Covariant Derivatives}

In this subsection which is mostly based on~\cite[Sec. 4.1]{schlarb:2023},
partially copied word by word,
we recall some facts on invariant covariant derivatives on
reductive homogeneous spaces.
We point out that invariant covariant derivatives correspond to the well-known
invariant affine connections from~\cite{nomizu:1954}.

In the sequel, let $G / H$ be
a reductive homogeneous space with a fixed reductive
decomposition $\liealg{g} = \liealg{h} \oplus \liealg{m}$.
We start with defining $\Ad(H)$-invariant bilinear maps and
invariant covariant derivatives.

\begin{definition}
	A bilinear map
	$\alpha \colon \liealg{m} \times \liealg{m} \to \liealg{m}$
	is called $\Ad(H)$-invariant if
	\begin{equation}
		\Ad_h \big( \alpha(X, Y) \big) 
		=
		\alpha\big( \Ad_h(X), \Ad_h(Y) \big) 
	\end{equation}
	holds for all $X, Y \in \liealg{m}$ and $h \in H$.
\end{definition}

\begin{definition}
	\label{definition:invariant_covariant_derivative_on_G_H_reductive}
	A covariant derivative
	$\nabla \colon \Secinfty\big(T (G / H)\big) \times \Secinfty\big(T (G / H)\big) \to \Secinfty\big( T(G / H) \big)$ on $G / H$
	is called $G$-invariant, or invariant for short,
	if 
	\begin{equation}
		\nabla_X Y = (\tau_{g^{-1}})_* \big( \nabla_{ (\tau_g)_* X}  (\tau_g)_* Y\big)
	\end{equation}
	holds for all $g \in G$ and $X, Y \in \Secinfty\big(T (G / H) \big)$,
	where
	$(\tau_g)_* X$ denotes the push-forward
	of $X$ by $\tau_g \colon G / H \to G / H$ given by
	$(\tau_g)_* X =  T  \tau_g \circ X \circ \tau_{g^{-1}}$.
\end{definition}
In order to establish the one-to-one correspondence between $\Ad(H)$-invariant bilinear maps and invariant covariant derivatives,
we introduce the following notation.
Let $X \in \liealg{g}$.
We denote by $X_{G / H} \in \Secinfty\big(T (G  /H) \big)$ the fundamental vector
field associated with the action $\tau \colon G \times G / H \to G / H$,
i.e. $X_{G / H}$ is defined by
\begin{equation}
	X_{G / H}(p) = \tfrac{\D}{\D t} \tau_{\exp(t X)}(p) \at{t = 0}, 
	\quad
	p \in G / H .
\end{equation}

\begin{lemma}
	\label{lemma:invariant_covariant_derivative_Ad_H_invariant_bilinear_map_one_to_one_correspondence}
	Let $\alpha \colon \liealg{m} \times \liealg{m} \to \liealg{m}$ be an
	$\Ad(H)$-invariant bilinear map.
	Then $\alpha$ corresponds to an invariant covariant derivative
	$\nablaAlpha$ on $G / H$
	by requiring
	\begin{equation}
		\label{equation:invariant_covariant_derivative_alpha_definition}
		\nablaAlpha_{X_{G / H}} Y_{G / H} \at{\pr(e)}
		=
		T_e \pr \big( -[X, Y]_{\liealg{m}} + \alpha(X, Y) \big) 
	\end{equation}
	for $X, Y \in \liealg{m}$.
	Conversely, for a given invariant covariant derivative,
	Equation~\eqref{equation:invariant_covariant_derivative_alpha_definition}
	defines an $\Ad(H)$-invariant bilinear map.
\end{lemma}

\begin{definition}
	\label{definition:invariant_covariant_derivative}
	Let $\alpha \colon \liealg{m} \times \liealg{m} \to \liealg{m}$ be an 
	$\Ad(H)$-invariant bilinear map.
	Then the invariant covariant derivative
	$\nablaAlpha \colon \Secinfty\big(T (G / H) \big) \times \Secinfty \big( T(G / H) \big) \to \Secinfty \big( T (G / H) \big)$
	defined in
	Lemma~\ref{lemma:invariant_covariant_derivative_Ad_H_invariant_bilinear_map_one_to_one_correspondence}
	is called the invariant covariant derivative associated
	with $\alpha$ or corresponding to $\alpha$.
\end{definition}

\begin{remark}
	Let $\alpha \colon \liealg{m} \times \liealg{m} \to \liealg{m}$
	be an $\Ad(H)$-invariant bilinear map.
	Then the invariant
	affine connection $t^{\alpha}$
	from~\cite[Thm. 8.1]{nomizu:1954} that corresponds to $\alpha$
	is related to $\nablaAlpha$
	by $\nablaAlpha_X Y = t^{\alpha}(Y)(X)$ for all $X, Y \in \Secinfty\big( T (G / H)\big)$ according to~\cite[Prop. 4.19]{schlarb:2023}.
\end{remark}

According to~\cite[Sec. 4.1.2]{schlarb:2023}
the invariant covariant derivative $\nablaAlpha$
from Definition~\ref{definition:invariant_covariant_derivative}
can be expressed in terms of horizontally lifted vector fields as follows.
Let $X, Y \in \Secinfty\big(T (G / H) \big)$ and let $\overline{X},
\overline{Y} \in \Secinfty\big(\Hor(G) \big)$
denote their horizontal lifts with respect to
the principal connection from
Proposition~\ref{proposition:principal_connection_reductive_homogeneous_space}.
Moreover, let $\{A_1, \ldots, A_N\} \subseteq \liealg{m}$ be a basis
and let $A_1^L, \ldots, A_N^L \in \Secinfty\big(\Hor(G) \big)$
be the associated left invariant frame.
Writing $\overline{X} = x^i A_i^L$ and
$\overline{Y} = y^j A_j^L$ with uniquely determined smooth
functions $x^i, y^j \colon G \to \field{R}$
for $i, j \in \{1, \ldots, N\}$
one obtains by~\cite[Thm. 4.16]{schlarb:2023} for the horizontal lift
$\overline{\nablaAlpha_X Y}$ of $\nablaAlpha_X Y$
\begin{equation}
	\overline{\nablaAlpha_X Y} 
	=
	\big( \Lie_{\overline{X}} y^j \big) A_j^L + x^i y^j \big( \alpha (A_i, A_j) \big)^L ,
\end{equation}
where Einstein summation convention is used, as usual.
This expression paves the way for characterizing parallel
vector fields along curves on $G / H$ with respect to $\nablaAlpha$.
Indeed, we have the next proposition
which is a reformulation of~\cite[Cor. 4.27]{schlarb:2023}.

\begin{proposition}
	\label{proposition:covariant_derivative_along_curve}
	Let ${\gamma} \colon I \to G /H$ be a curve in the
	reductive homogeneous space $G / H$ with reductive decomposition
	$\liealg{g} = \liealg{h} \oplus \liealg{m}$
	and let
	$g \colon I \to G$ be a horizontal lift of $\gamma$
	with respect to the principal connection from
	Proposition~\ref{proposition:principal_connection_reductive_homogeneous_space}.
	Moreover, let $\widehat{Z} \colon I \to T (G / H)$ be a vector
	field along $\gamma$ with horizontal lift
	$\overline{Z} \colon I \ni 
	t \mapsto \big(T_{g(t)} \pr \at{\Hor(G)_{g(t)}} \big)^{-1} \widehat{Z}(t) 
	\in \Hor(G)$
	along $g \colon I \to G$.
	Define the curves $x, z \colon I \to \liealg{m}$ by
	\begin{equation}
		x(t) = (T_e \ell_{g(t)})^{-1} \dot{g}(t)
		\quad \text{ and } \quad
		z(t) = (T_e \ell_{g(t)})^{-1} \overline{Z}(t) 
	\end{equation}
	for $t \in I$.
	Then $\widehat{Z} \colon I \to T ( G / H)$ is parallel along
	$\gamma \colon I \to G / H$ with respect to $\nablaAlpha$
	defined by the $\Ad(H)$-invariant bilinear map
	$\alpha \colon \liealg{m} \times \liealg{m} \to \liealg{m}$
	iff the ODE
	\begin{equation}
		\dot{z}(t) = - \alpha\big(x(t), z(t) \big)
	\end{equation}
	is satisfied for all $t \in I$.
\end{proposition}

\subsubsection{Metric Invariant Covariant Derivatives}

If $G / H$ is equipped with an invariant pseudo-Riemmannian metric,
one has the following characterization of invariant metric covariant
derivatives which is copied from~\cite[Prop. 4.22]{schlarb:2023}.

\begin{proposition}
	\label{proposition:invariant_covariant_derivative_compatible_with_structure_skew_adjoint}
	Let $\alpha \colon \liealg{m} \times \liealg{m} \to \liealg{m}$ be
	an $\Ad(H)$-invariant bilinear map defining the invariant
	covariant derivative $\nablaAlpha$ on $G / H$.
	Then $\nablaAlpha$ is metric with respect to the
	invariant pseudo-Riemannian
	metric on $G / H$ defined by
	the $\Ad(H)$-invariant scalar product 
	$\langle \cdot, \cdot \rangle \colon \liealg{m} \times \liealg{m} \to \field{R}$ iff
	for each $X \in \liealg{m}$
	the linear map
	\begin{equation}
		\alpha(X, \cdot) \colon \liealg{m} \to \liealg{m},
		\quad
		Y \mapsto \alpha(X, Y)
	\end{equation}
	is skew-adjoint
	with respect to $\langle \cdot, \cdot \rangle$, i.e. 
	\begin{equation}
		\label{equation:lemma_invariant_covariant_derivative_compatible_with_structure_skew_adjoint_skew_adjoint}
		\big\langle \alpha(X, Y) , Z \big\rangle 
		=
		- \big\langle Y, \alpha(X, Z) \big\rangle
	\end{equation} 
	holds for all $X, Y, Z \in \liealg{m}$.
\end{proposition}

Moreover, we state the following
reformulation of~\cite[Re. 4.25]{schlarb:2023}
concerning the Levi-Civita
covariant derivative of a naturally reductive homogeneous space.

\begin{corollary}
	\label{corollary:nabla_LC_naturally_reductive_homogeneous_space}
	Let $G / H$ be a naturally reductive homogeneous space.
	Then the Levi-Civita covariant derivative $\nablaLC$
	coincides with the invariant covariant derivative $\nablaAlpha$
	corresponding to the $\Ad(H)$-invariant bilinear map
	\begin{equation}
		\alpha \colon \liealg{m} \times \liealg{m} \to \liealg{m},
		\quad
		(X, Y) \mapsto \alpha(X, Y) = \tfrac{1}{2} [X, Y]_{\liealg{m}},
	\end{equation}
	i.e. $\nablaLC = \nablaAlpha$ holds.
\end{corollary}

\subsubsection{Canonical Invariant Covariant Derivatives}

In this subsection, we consider
the so-called canonical invariant covariant derivatives
of first and second kind
which correspond to the canonical affine connections
from~\cite[Sec. 10]{nomizu:1954}.
Here we copy and adapt some
parts of~\cite[Sec. 4.6]{schlarb:2023}.

\begin{definition}
	The canonical covariant derivative of first kind
	$\nablaCan$ is defined by the $\Ad(H)$-invariant
	bilinear map given by
	\begin{equation}
		\alpha \colon \liealg{m} \times \liealg{m} \to \liealg{m},
		\quad
		(X, Y) \mapsto \alpha(X, Y) =  \tfrac{1}{2} [X,Y]_{\liealg{m}}.
	\end{equation}
	The canonical covariant derivative of second kind
	$\nablaCanSecond$ corresponds to the
	$\Ad(H)$-invariant bilinear map
	\begin{equation}
		\alpha \colon \liealg{m} \times \liealg{m} \to \liealg{m},
		\quad 
		(X, Y) \mapsto \alpha(X, Y) = 0 .
	\end{equation}
\end{definition}
These canonical
covariant derivatives correspond to
the Levi-Civita covariant derivatives
on certain pseudo-Riemannian homogeneous spaces.
More precisely, we have the next remarks
which are copied from~\cite{schlarb:2023}.

\begin{remark}
	\label{remark:horizontal_lift_covariant_derivative_naturally_redutive_and_symmetric}
	Assume that $G / H$ is a naturally reductive homogeneous space.
	Then the Levi-Civita covariant derivative coincides with the canonical
	covariant derivative of first kind, i.e $\nablaLC = \nablaCan$
	by
	Corollary~\ref{corollary:nabla_LC_naturally_reductive_homogeneous_space}.
	This has already been proven
	in~\cite[Thm. 13.1 and Eq. (13.2)]{nomizu:1954}.
\end{remark}
We briefly comment on the canonical covariant derivatives on
symmetric homogeneous spaces in the next remark
following~\cite[Thm. 15.1]{nomizu:1954}.

\begin{remark}
	\label{remark:symmetric_spaces_invariant_covariant_derivatives}
	Let $(G, H, \sigma)$ be a symmetric pair and let
	$G / H$ be the corresponding
	symmetric homogeneous space.
	Let
	$\liealg{g} = \liealg{h} \oplus \liealg{m}$
	denote the canonical reductive decomposition.
	Then $[X, Y] \in \liealg{h}$ holds for all
	$X, Y \in \liealg{m}$
	by Lemma~\ref{lemma:symmetric_pair_canonical_reductive_decomposition}.
	Therefore $\tfrac{1}{2}[X, Y]_{\liealg{m}} = 0$ 
	is fulfilled for all $X, Y \in \liealg{m}$.
	Hence $\nablaCan = \nablaCanSecond$ is fulfilled.
	Moreover, if $G / H$ is a pseudo-Riemannian symmetric space,
	then $\nablaLC = \nablaCan = \nablaCanSecond$ holds
	by Remark~\ref{remark:horizontal_lift_covariant_derivative_naturally_redutive_and_symmetric}
	combined with
	Remark~\ref{remark:pseudo-Riemannian-symmetric-naturally-reductive}.
\end{remark}

\section{Intrinsic Rolling}
\label{sec:rolling_general}

In this section, a
notion of rolling intrinsically a manifold $M$ over another manifold
$\widehat{M}$ of equal dimension $\dim(M) = n = \dim(\widehat{M})$ is
recalled from the literature and slightly generalized.
As preparation to define the configuration space, we state the following
lemma which can be regarded as
a slight generalization of the
definition of the configuration space in~\cite[Sec. 3.1]{molina.grong.markina.leite:2012}.
In particular, the definition of the map $\Psi$ in
Lemma~\ref{lemma:configuration_space_is_fiber_bundle},
Claim~\ref{item:lemma_configuration_space_is_fiber_bundle_fiber_isomorphisms_identification}, below,
is very similar to~\cite[Eq. (4)]{molina.grong.markina.leite:2012}.

\begin{lemma}
	\label{lemma:configuration_space_is_fiber_bundle}
	Let $E \to M$
	and $\widehat{E} \to \widehat{M}$ be two vector bundles
	both having typical fiber $V$
	and let $\liegroupsub(V) \subseteq \liegroup{GL}(V)$ be a closed subgroup.
	Assume that the frame bundles of $E$ and $\widehat{E}$
	admit both a $\liegroupsub(V)$-reduction
	along the canonical inclusion $\liegroupsub(V) \to \liegroup{GL}(V)$
	which
	we denote by $\liegroupsub(V, E) \to M$ and
	$\liegroupsub(V, \widehat{E}) \to \widehat{M}$, respectively.
	Let
	\begin{equation}
		Q 
		=
		\big(\liegroupsub(V, E) \times \liegroupsub(V, \widehat{E}) \big)
		/ \liegroupsub(V)
	\end{equation}
	be defined as the quotient of
	$\liegroupsub(V, E) \times \liegroupsub(V, \widehat{E})$
	by the diagonal action of $\liegroupsub(V)$,
	where the action on each component is given by the
	$\liegroupsub(V)$-principal action.
	Moreover, define
	\begin{equation}
			\pi \colon Q 
			\to M \times \widehat{M}, 
			\quad
			[f, \widehat{f}] 
			\mapsto
			\big( \pr_{\liegroupsub(V, E)}(f), \pr_{\liegroupsub(V, \widehat{E})}(\widehat{f})\big) ,
	\end{equation}
	where $[f, \widehat{f}] \in Q$ denotes the equivalence class defined by $(f, \widehat{f}) \in \liegroupsub(V, E) \times \liegroupsub(V, \widehat{E})$.
	Then the following assertions are fulfilled:
	\begin{enumerate}
		\item
		\label{item:lemma_configuration_space_is_fiber_bundle_fiber_bundle}
		$\pi \colon Q = \big(\liegroupsub(V, E) \times \liegroupsub(V, \widehat{E}) \big) / \liegroupsub(V) \to M \times \widehat{M}$
		is a $\liegroupsub(V)$-fiber bundle over $M \times \widehat{M}$.
		\item
		\label{item:lemma_configuration_space_is_fiber_bundle_fiber_isomorphisms_identification}
		Let $(x, \widehat{x}) \in M \times \widehat{M}$ and define
		\begin{equation}
			\widetilde{Q}_{(x, \widehat{x})}
			=
			\big\{
			\widetilde{q} \colon E_x \to \widehat{E}_{\widehat{x}} 
			\mid
			\widehat{f}^{-1} \circ \widetilde{q} \circ f  
			\in \liegroupsub(V)
			\text{ for all }
			(f, \widehat{f}) \in \liegroupsub(V, E)_x \times \liegroupsub(V, \widehat{E}_{\widehat{x}})
			\} .
		\end{equation}
		Then the map
		\begin{equation}
			\begin{split}
				\Psi \colon Q =
				\big(\big(\liegroupsub(V, E) \times \liegroupsub(V, \widehat{E}) \big)
				/ \liegroupsub(V) \big)_{(x, \widehat{x})}
				\ni [f, \widehat{f}] \mapsto \widehat{f} \circ f^{-1} \in
				\widetilde{Q}_{(x, \widehat{x})}			
			\end{split}
		\end{equation}
		is bijective.
	\end{enumerate}
	\begin{proof}
		The action 
		\begin{equation*}
			\begin{split}
				\big(\liegroupsub(V, E) \times \liegroupsub(V, \widehat{E}) \big) \times \liegroupsub(V) 
				&\to
				\liegroupsub(V, E) \times \liegroupsub(V, \widehat{E}), \\
				\quad
				((f, \widehat{f}), A) 
				&\mapsto (f \racts A, \widehat{f} \racts A)
			\end{split}
		\end{equation*}
		is free and proper since the action on each component is free and proper.
		Thus
		$Q = 
		\big(\liegroupsub(V, E) \times \liegroupsub(V, \widehat{E}) \big)  / \liegroupsub(V)$
		is a smooth manifold.
		Moreover, $(f, \widehat{f}) \sim (f^{\prime}, \widehat{f}^{\prime})$
		holds
		iff there is an $A \in \liegroupsub(V)$
		such that 
		$(f^{\prime}, \widehat{f}^{\prime})
		=
		(f \racts A, \widehat{f} \racts A)$
		is fulfilled.
		Let $(x, \widehat{x}) \in M \times \widehat{M}$ and let
		$U \subseteq M$ and $\widehat{U} \subseteq \widehat{M}$ be
		open neighbourhoods of $x$ and $\widehat{x}$, respectively,
		such that
		\begin{equation*}
			\varphi \colon \liegroupsub(V, E)\at{U} \to U \times \liegroupsub
			(V)
			\quad \text{ and } \quad
			\widehat{\varphi} \colon \liegroupsub(V, \widehat{E})\at{\widehat{U}} \to \widehat{U} \times \liegroupsub
			(V)
		\end{equation*}
		are local trivializations of $ \liegroupsub(V, E) \to M$ and $\liegroupsub(V, \widehat{E}) \to \widehat{M}$ as
		$\liegroupsub(V)$-principal fiber bundles, respectively.
		Locally, one obtains for the principal action
		for $A \in \liegroupsub(V)$
		\begin{equation}
			\label{equation:lemma_configuration_space_is_fiber_bundle_principal_action_local_trivialization}
			\varphi(f \racts A) = (\pr_1(\varphi(f)), \pr_2(\varphi(f)) \circ A),
			\quad
			\widehat{\varphi}(\widehat{f} \racts A) = (\pr_1(\widehat{\varphi}(\widehat{f})), \pr_2(\widehat{\varphi}(\widehat{f})) \circ A),
		\end{equation}
		see e.g.~\cite[Sec. 18, p. 211]{michor:2008}.
		We now define the local trivialization
		$\phi \colon Q\at{U \times \widehat{U}}
		\to U \times \widehat{U} \times \liegroupsub(V)$
		of $\pi \colon Q \to M \times \widehat{M}$
		by
		\begin{equation*}
			\phi\big([f, \widehat{f}]\big)
			=
			\big( 
			\pr_1(\varphi(f)),
			\pr_1(\widehat{\varphi}(\widehat{f})),
			\big(\pr_2(\varphi(f)) \big) \circ \big(\pr_2(\widehat{\varphi}(\widehat{f})) \big)^{-1} \big) ,
			\quad
			[f, \widehat{f}] \in Q\at{U \times \widehat{U}} .
		\end{equation*}
		The map $\phi$ is well-defined. Indeed, we obtain
		by~\eqref{equation:lemma_configuration_space_is_fiber_bundle_principal_action_local_trivialization}
		for
		$(f , \widehat{f} ) \in \big( \liegroupsub(V, E) \times \liegroupsub(V, \widehat{E}) \big)\at{U \times \widehat{U}}$ 
		and $A \in \liegroupsub(V)$	
		\begin{equation*}
			\begin{split}
				\phi\big([f \racts A, \widehat{f} \racts A ]\big)
				&=
				\Big(\pr_1(\varphi(f \racts A)),
				\pr_1(\widehat{\varphi}(\widehat{f} \racts A)),
				\big(\pr_2(\varphi(f \racts A)) \big) \circ \big(\pr_2(\widehat{\varphi}(\widehat{f} \racts A) ) \big)^{-1} \Big) \\
				&=
				\Big(\pr_1(\varphi(f)),
				\pr_1(\widehat{\varphi}(\widehat{f})),
				\big(\pr_2(\varphi(f)) \circ A \big) \circ \big(\pr_2(\widehat{\varphi}(\widehat{f} ) ) \circ A \big)^{-1}  \Big) \\
				&=
				\phi\big( [f, \widehat{f} ] \big) .
			\end{split}
		\end{equation*}
		Moreover, the
		canonical projection
		$\overline{\pr} \colon 
		\big(\liegroupsub(V, E) \times \liegroupsub(V, \widehat{E}) \big)
		\ni (f, \widehat{f}) \mapsto [f, \widehat{f}] \in 
		\big( \liegroupsub(V, E) \times \liegroupsub(V, \widehat{E}) \big)  / \liegroupsub(V)$ is 
		a surjective submersion and the map
		\begin{equation*}
			\overline{\phi} \colon \big(\liegroupsub(V, E) \times \liegroupsub(V, \widehat{E}) \big) \at{U \times \widehat{U}}
			\to 
			U \times \widehat{U} \times \liegroupsub(V)
		\end{equation*}
		defined for
		$(f, \widehat{f}) \in \big( \liegroupsub(V, E) \times \liegroupsub(V, \widehat{E})\big)\at{U \times \widehat{U}}$
		by 
		\begin{equation*}
			\overline{\phi}(f, \widehat{f}) 
			= 
			\Big( 
			\pr_1(\varphi(f)),
			\pr_1(\widehat{\varphi}(\widehat{f})),
			\big(\pr_2(\varphi(f)) \big) \circ \big(\pr_2(\widehat{\varphi}(\widehat{f})) \big)^{-1} \Big) 
		\end{equation*}
		is smooth as the composition of smooth maps.
		Obviously, $\phi = \overline{\pr} \circ \overline{\phi}$
		holds. Hence $\phi$ is smooth, too.
		Moreover, $\phi$ is bijective. Indeed, the map
		\begin{equation*}
			\phi^{-1} \colon U \times \widehat{U} \times \liegroupsub(V) \to Q\at{U \times \widehat{U}},
			\quad
			(x, \widehat{x}, A) 
			\mapsto
			\phi^{-1}(x, \widehat{x}, A)
			=
			\big[\varphi^{-1}(x, A), \widehat{\varphi}^{-1}(\widehat{x}, \id_{V}) \big]
		\end{equation*}
		yields the inverse
		of $\phi$ since one verifies
		by a straightforward computation
		\begin{equation*}
			\phi \circ \phi^{-1} = \id_{U \times \widehat{U} \times \liegroupsub(V)}
			\quad \text{ and } \quad
			\phi^{-1} \circ \phi = \id_{Q\at{U \times \widehat{U}}} .
		\end{equation*}
		Next we define the map
		\begin{equation*}
			\overline{\phi^{-1}} \colon U \times \widehat{U} \times \liegroupsub(V) \to \liegroupsub(V, E) \times \liegroupsub(V, \widehat{E}),
			\quad
			(x, \widehat{x}, A) \mapsto
			\big(\varphi^{-1}(x, A), \widehat{\varphi}^{-1}(\widehat{x}, A) \big)
		\end{equation*}
		which is smooth as the composition of smooth maps.
		Then $\phi^{-1} = \overline{\pr} \circ \overline{\phi^{-1}}$ is clearly fulfilled.
		Hence $\phi^{-1}$ is smooth
		by the smoothness of $\overline{\phi^{-1}}$
		since $\overline{\pr}$ is a surjective submersion.
		Thus $\phi$ is a diffeomorphism.
		This yield the desired result since for every point
		$(x, \widehat{x}) \in M \times \widehat{M}$, we can construct a
		local trivialization
		$\phi \colon Q\at{U \times \widehat{U}} \to U \times \widehat{U} \times \liegroupsub(V)$,
		i.e.
		Claim~\ref{item:lemma_configuration_space_is_fiber_bundle_fiber_bundle}
		is shown.

		It remains to prove
		Claim~\ref{item:lemma_configuration_space_is_fiber_bundle_fiber_isomorphisms_identification}.
		Let $(x, \widehat{x}) \in M \times \widehat{M}$ and
		$f \in \liegroupsub(V, E)_x$
		as well as
		$\widehat{f} \in \liegroupsub(V, \widehat{E})_{\widehat{x}}$.
		In particular, $f \colon V \to E_x$ and
		$\widehat{f} \colon V \to \widehat{E}_{\widehat{x}}$
		are invertible linear maps.
		Hence $\widehat{f} \circ f^{-1} \colon E_x \to \widehat{E}_{\widehat{x}}$
		is a linear isomorphism.
		
		Moreover,
		$\Psi\big([f, \widehat{f}]\big)$
		is independent of the representative of $[f, \widehat{f}] \in Q$ 
		due to
		\begin{equation*}
			\Psi\big([f \circ A, \widehat{f} \circ A]\big) 
			=
			(\widehat{f} \circ A) \circ (f \circ A)^{-1}
			=
			\widehat{f} \circ f^{-1}
			=
			\Psi\big([f, \widehat{f}]\big) 
		\end{equation*}	
		for all $A \in \liegroupsub(V)$.
		
		Next we show that
		$\Psi\big([f, \widehat{f}]\big) \in \widetilde{Q}_{(x, \widehat{x})}$
		holds for all $[f, \widehat{f}] \in Q_{(x, \widehat{x})}$.
		Let $[f, \widehat{f}] \in Q_{(x, \widehat{x})}$.
		By the fiber-wise transitivity of the principal
		$\liegroupsub(V)$-actions on $\liegroupsub(V, E)$
		and $\liegroupsub(V, \widehat{E})$, respectively,
		we obtain for $A, B \in \liegroupsub(V)$
		\begin{equation*}
			(\widehat{f} \circ B)^{-1} \circ \Psi\big([f , \widehat{f} ]\big) \circ (f \circ A)
			=
			B^{-1} \circ \widehat{f}^{-1} \circ \big(\widehat{f} \circ f^{-1} \big) \circ f \circ A
			=
			B^{-1} \circ A 
			\in \liegroupsub(V)
		\end{equation*}
		showing
		$\Psi\big([f, \widehat{f}]\big) \in \widetilde{Q}_{(x, \widehat{x})}$ 
		for all $[f, \widehat{f}] \in Q_{(x, \widehat{x})}$,
		i.e. $\Psi \colon Q_{(x, \widehat{x})} \to \widetilde{Q}_{(x, \widehat{x})}$
		is well-defined.
		Moreover, $\Psi$ is injective.
		Let $[f, \widehat{f}], [f^{\prime}, \widehat{f}^{\prime}] \in Q_{(x, \widehat{x})}$
		with
		$\Psi\big([f, \widehat{f}]\big) 
		=
		\Psi\big([f^{\prime}, \widehat{f}^{\prime}]\big)$.
		Since the $\liegroupsub(V)$-principal actions on
		$\liegroupsub(V, E)$ and $\liegroupsub(V, \widehat{E})$
		are free and fiber-wise transitive, 
		we can write $f^{\prime} = f \circ A$ and
		$\widehat{f}^{\prime} = \widehat{f} \circ B$
		with some uniquely determined $A, B \in \liegroupsub(V)$.
		By this notation, we obtain
		\begin{equation*}
			\widehat{f} \circ f^{-1} 
			=
			\Psi\big([f, \widehat{f}]\big) 
			=
			\Psi\big([f \circ A, \widehat{f} \circ B]\big) 
			=
			(\widehat{f} \circ B) \circ (f \circ A)^{-1}
			=
			\widehat{f} \circ (B \circ A^{-1}) \circ f^{-1},
		\end{equation*}
		implying $B \circ A^{-1} = \id_V \iff A = B$ 
		because
		$f \colon V \to E_x$ and
		$\widehat{f} \colon V \to \widehat{E}_{\widehat{x}}$ 
		are both linear isomorphisms.
		Thus $[f^{\prime}, \widehat{f}^{\prime}] = [f \circ A, \widehat{f} \circ A] = [f, \widehat{f}]$ is shown.
		It remains to show that $\Psi$ is surjective.
		To this end, let $\widetilde{q} \in \widetilde{Q}_{(x, \widehat{x})}$
		and chose some $f \in \liegroupsub(V, E)_{(x, \widehat{x})}$
		and
		$\widehat{f} \in \liegroupsub(V, \widehat{E})_{(x, \widehat{x})}$.
		Then $\widehat{f}^{-1} \circ \widetilde{q} \circ f \in \liegroupsub(V)$ holds.
		We now compute
		\begin{equation*}
			\Psi\big([f, \widehat{f} \circ (\widehat{f}^{-1} \circ \widetilde{q} \circ f)]\big)
			=
			\big(\widehat{f} \circ (\widehat{f}^{-1} \circ \widetilde{q} \circ f)\big) \circ f^{-1} 
			=
			(\widehat{f} \circ \widehat{f}^{-1}) \circ \widetilde{q} \circ (f \circ f^{-1})
			= 
			\widetilde{q} ,
		\end{equation*}
		i.e. $\Psi$ is surjective.
		This yields the desired result.
	\end{proof}
\end{lemma}

After this preparation, we consider intrinsic rollings.
Let $M$ and $\widehat{M}$ be two manifolds 
with $\dim(M) = n  = \dim(\widehat{M})$.
Moreover, let
$\liegroupsub(\field{R}^n) \subseteq \liegroup{GL}(\field{R}^n)$
be a closed subgroup
and assume that the frame bundles
$\liegroup{GL}(\field{R}^n, TM) \to M$
and $\liegroup{GL}(\field{R}^n, T \widehat{M}) \to \widehat{M}$
admit both a $\liegroupsub(\field{R}^n)$-reduction along the
canonical inclusion
$\liegroupsub(\field{R}^n) \to \liegroup{GL}(\field{R}^n)$.
These reductions are denoted by
\begin{equation}
	\label{equation:frame_bundles_reduction_M_hat_M_rolling_def}
	\liegroupsub(\field{R}^n, TM) \to \liegroup{GL}(\field{R}^n, TM)
	\quad \text{ and } \quad
	\liegroupsub(\field{R}^n, T\widehat{M}) \to \liegroup{GL}(\field{R}^n, T\widehat{M}),
\end{equation}
respectively.
In this section, we denote by
\begin{equation}
	\prQ \colon \Qspace 
	= 
	\big( \liegroupsub(\field{R}^n, T M) \times \liegroupsub(\field{R}^n, T \widehat{M})\big) / \liegroupsub(\field{R}^n) 
	\to
	M \times \widehat{M}
\end{equation}
the $\liegroupsub(\field{R}^n)$-fiber bundle over $M \times \widehat{M}$
obtained by applying Lemma~\ref{lemma:configuration_space_is_fiber_bundle}
to the frame bundles
from~\eqref{equation:frame_bundles_reduction_M_hat_M_rolling_def}.

We now define a notion of rolling of $M$ over $\widehat{M}$ intrinsically,
where $M$ and $\widehat{M}$ are both equipped with a covariant
derivative $\nabla$ and $\widehat{\nabla}$, respectively.

\begin{definition}
	\label{definition:intrinsic_rolling_reduced}
	An intrinsic ($\liegroupsub(\field{R}^n)$-reduced)
	rolling of $(M, \nabla)$ over $(\widehat{M}, \widehat{\nabla})$
	is a curve
	\begin{equation}
		q \colon I \to Q =
		\big(\liegroupsub(\field{R}^n, TM) \times \liegroupsub(\field{R}^n, TM) \big) / \liegroupsub(\field{R}^n)
	\end{equation}
	with projection
	$(x, \widehat{x}) = \prQ \circ q \colon I \to M \times \widehat{M}$
	such that the following conditions are fulfilled:
	\begin{enumerate}
		\item
		\label{item:definition_intrinsic_rolling_reduced_no_slip}
		No slip condition: $\dot{\widehat{x}}(t) = q(t) \dot{x}(t)$ for all $t \in I$.
		\item
		\label{item:definition_intrinsic_rolling_reduced_no_twist}
		No twist condition: $Z \colon I \to TM$ is a parallel vector field along
		$x$ iff $\widehat{Z} \colon I \to T \widehat{M}$ defined by
		$\widehat{Z}(t) = q(t) Z(t)$ for $t \in I$ is parallel along $\widehat{x}$.
	\end{enumerate}
	Here Lemma~\ref{lemma:configuration_space_is_fiber_bundle},
	Claim~\ref{item:lemma_configuration_space_is_fiber_bundle_fiber_isomorphisms_identification}
	is used to identify $q(t)$ with the linear isomorphism 
	$q(t) \colon T_{x(t)} M \to T_{\widehat{x}(t)} \widehat{M}$
	which is
	denoted by $q(t)$, as well.
	We call the curve $x \colon I \to M$ rolling curve.
	The curve $\widehat{x} \colon I \to \widehat{M}$
	is called development curve.
	The curve $q \colon I \to Q$ is often called rolling for short.
\end{definition}
The next remark yields an other perspective
on the intrinsic rollings from
Definition~\ref{definition:intrinsic_rolling_reduced}.

\begin{remark}
	\label{remark:intrinsic_rolling_as_triple}
	Let $q \colon I \to Q$ be a ($\liegroupsub(\field{R}^n)$-reduced)
	intrinsic rolling of $M$ over $\widehat{M}$ in the sense of
	Definition~\ref{definition:intrinsic_rolling_reduced}
	and write $(x, \widehat{x}) = \pi \circ q \colon I \to M \times \widehat{M}$.
	Then we can view this rolling as a triple
	$(x(t), \widehat{x}(t), A(t))$,
	where $A(t) = q(t) \colon T_{x(t)} M \to T_{\widehat{x}(t)} \widehat{M}$
	is the linear isomorphism defined by $q(t)$ as in
	Lemma~\ref{lemma:configuration_space_is_fiber_bundle},
	Claim~\ref{item:lemma_configuration_space_is_fiber_bundle_fiber_isomorphisms_identification}.
	This point of view allows for relating a
	rolling $q \colon I \to Q$ from
	Definition~\ref{definition:intrinsic_rolling_reduced}
	to~\cite[Def. 1]{jurdjevic.markina.leite:2023},
	where a rolling is defined as
	a triple $(x(t), \widehat{x}(t), A(t))$ satisfying
	certain properties.
\end{remark}

Definition~\ref{definition:intrinsic_rolling_reduced}
of an intrinsic rolling of $M$ over $\widehat{M}$
generalizes several notions of intrinsic rolling from the literature.

\begin{remark}
	Assume that $M$ and $\widehat{M}$ are both orientible and
	both equipped with a Riemannian metric.
	Let $\liegroup{SO}(\field{R}^n, T M)$ and 
	$\liegroup{SO}(\field{R}^n, T \widehat{M})$ be the
	corresponding reductions of their frame bundles.
	Moreover, let $M$ and $\widehat{M}$ be endowed with the Levi-Civita
	covariant derivatives $\nablaLC$ and $\widehat{\nablaLC}$
	corresponding to the Riemannian metrics on $M$ and
	$\widehat{M}$, respectively.
	Then Definition~\ref{definition:intrinsic_rolling_reduced}
	specializes to~\cite[Def. 3]{molina.grong.markina.leite:2012}.
	Here the no twist condition is rewritten as
	in~\cite[Prop. 2]{jurdjevic.markina.leite:2023}.
	More generally, if $M$ and $\widehat{M}$ are oriented and equipped with
	a pseudo-Riemannian metric,
	Definition~\ref{definition:intrinsic_rolling_reduced}
	specializes to~\cite[Def. 4]{markina.leite:2016}.	
	If $M$ and $\widehat{M}$ are both equipped with an arbitrary
	covariant derivative $\nabla$ and $\widehat{\nabla}$, respectively,
	Definition~\ref{definition:intrinsic_rolling_reduced}
	yields the definition proposed in~\cite[p. 35]{kokkonen:2012}
	and~\cite[Sec. 7]{grong:2012}
	by setting $\liegroupsub
	(\field{R}^n, TM) = \liegroup{GL}(\field{R}^n, TM)$ and
	$\liegroupsub(\field{R}^n, T\widehat{M}) = \liegroup{GL}(\field{R}^n, T\widehat{M})$.
\end{remark}

Studying properties of rollings
in the sense of 
Definition~\ref{definition:intrinsic_rolling_reduced}
for general manifolds
is out of the scope of this text.
However, in
Section~\ref{sec:intrinsic_rolling_reductive_space}
below, we
discuss intrinsic rollings in the context of reductive 
homogeneous spaces in detail.

\section{Frame Bundles of Associated Vector Bundles}
\label{sec:frame_bundles}

In this section, we identify
(certain reductions of) the frame bundle
of a reductive homogeneous space $G / H$ with certain principal fiber bundles obtained as associated bundles of the $H$-principal fiber
bundle $\pr \colon G \to G / H$.
We point out that the results of this section might be well-known
since the statement of
Corollary~\ref{corollary:reductive_homogeneous_space_frame_bundle}
can be found as an exercise in the German book~\cite[Ex. 2.7]{baum:2009}.
However, we were not able to find a reference including a proof.
Hence we provide one in this section in order to keep this text as self-contained as possible.
Here we first start with a more general situation that
is applied to reductive homogeneous spaces later.
We first determine (certain reductions of) the frame bundles
of vector bundles 
given as associated bundles of some principal fiber bundle.

\subsection{Frame Bundles of Associated Vector Bundles}

Let $P \to M$ be a $H$-principal fiber bundle.
We describe
(reductions of) the frame bundle of a vector bundle
associated to $P$
in terms of another fiber bundle associated to $P$.
To this end, we state the following lemma as preparation.

\begin{lemma}
	\label{lemma:PxHGV_principal_fiber_bundle_structure}
	Let $P \to M$ be a $H$-principal fiber bundle and let
	$\rho \colon H \to \liegroup{GL}(V)$ be a smooth representation
	of $H$ on a finite dimensional
	$\field{R}$-vector space $V$.
	Moreover, let $\liegroupsub(V) \subseteq \liegroup{GL}(V)$
	be a closed subgroup such that $\rho_h \in \liegroupsub(V)$
	is fulfilled for all $h \in H$.
	Then the following assertions are fulfilled:
	\begin{enumerate}
		\item
		\label{item:lemma_PxHGV_principal_fiber_bundle_structure_H_action}
		The Lie group $H$ acts on $\liegroupsub(V)$ via
		\begin{equation}
			H \times \liegroupsub(V) \to \liegroupsub(V),
			\quad
			(h, A) \mapsto \rho_h \circ A
		\end{equation}
		smoothly from the left.
		\item
		\label{item:lemma_PxHGV_principal_fiber_bundle_structure_GV_principal_action}
		The map
		\begin{equation}
			\racts \colon (P \times_H \liegroupsub(V)) \times \liegroupsub(V) \to P \times_H \liegroupsub(V)
			\quad
			([g, A], B ) \to [g, A \circ B] ,
		\end{equation}
		denoted by the same symbol as the principal
		action $\racts \colon P \times H \to P$,
		yields a well-defined, smooth, free and proper
		$\liegroupsub(V)$-right action
		on the associated bundle
		\begin{equation}
			\pi \colon P \times_H \liegroupsub(V) \to M
		\end{equation}
		turning
		\begin{equation}
			\widetilde{\pi} \colon P \times_H \liegroupsub(V) \to 
			\big(P \times_H \liegroupsub(V) \big) /  \liegroupsub(V)
		\end{equation}
		into a $\liegroupsub(V)$-principal fiber bundle,
		where $\widetilde{\pi}$ denotes the canonical projection.
		\item
		\label{item:lemma_PxHGV_principal_fiber_bundle_structure_isomorphism_base}
		The map
		\begin{equation}
			\phi \colon (P \times_H \liegroupsub(V) ) /  \liegroupsub(V) \ni \widetilde{\pi}([p, S])
			\mapsto \pr(p) \in M
		\end{equation}
		is a diffeomorphism such that $\phi \circ \widetilde{\pi} = \pi$ holds.
		Moreover,
		\begin{equation}
			\id_{G \times_H \liegroupsub(\liealg{m})} \colon G \times_H \liegroupsub(\liealg{m}) \to G \times_H \liegroupsub(\liealg{m})
		\end{equation}
		is an isomorphism of $\liegroupsub(V)$-principal fiber bundles
		covering $\phi$.
	\end{enumerate}
	\begin{proof}	
		Claim~\ref{item:lemma_PxHGV_principal_fiber_bundle_structure_H_action}
		is obvious.
		
		We now show
		Claim~\ref{item:lemma_PxHGV_principal_fiber_bundle_structure_GV_principal_action}.
		The  $\liegroupsub(V)$-right action
		$\racts$ on $P \times_H \liegroupsub(V)$
		is well-defined due to 
		\begin{equation*}
			[p \racts h, \rho_{h^{-1}} \circ A] \racts B 
			=
			[p \racts h, \rho_{h^{-1}} \circ A \circ B] 
			=
			[g, A \circ B] 
			=
			[g, A] \racts B
		\end{equation*}
		for all $p \in P$, $h \in H$ and $A, B \in \liegroupsub(V)$.
		Next we show that $\racts$ is smooth. To this end, we consider the diagram
		\begin{equation}
			\label{equation:lemma_GxHOm_principal_fiber_bundle_structure_diagram}
			\begin{tikzpicture}[node distance= 5.0cm, auto]
				\node (GxOmOm) {$(P \times \liegroupsub(V)) \times \liegroupsub(V)$};
				\node (GxHOmOm) [below= 1.5cm of GxOmOm] {$(P \times_H 	\liegroupsub(V) ) \times \liegroupsub(V)$};
				\node (GxHOm) [right=3.0cm of GxHOmOm] {$P \times_H \liegroupsub(V)$};
				\draw[->] (GxOmOm) to node [left] {$\overline{\pi} \times \id_{\liegroupsub(V)}$} (GxHOmOm);	
				\draw[->] (GxHOmOm) to node  [above]{$\racts$} (GxHOm);
				\draw[->] (GxOmOm) to node [pos=0.5] {$(\overline{\pi} \times 	\id_{\liegroupsub(V)}) \circ \widetilde{\racts}$} (GxHOm);			
			\end{tikzpicture}  
		\end{equation}
		where
		$\overline{\pi} \colon P \times \liegroupsub(V) \to ( P \times \liegroupsub(V)) / H
		= P \times_H \liegroupsub(V)$
		denotes the canonical projection and
		$\widetilde{\racts}$ is given by
		\begin{equation*}
			\widetilde{\racts} 
			\colon (P \times \liegroupsub(V)) \times \liegroupsub(V) \to P \times \liegroupsub(V) ,
			\quad
			((p, A), B) \mapsto (p, A \circ B)
		\end{equation*}
		which is clearly a smooth and free $\liegroupsub(V)$-right action
		on $P \times \liegroupsub(V)$.
		Moreover, the action $\widetilde{\racts}$ is proper since
		the $\liegroupsub(V)$-action on $\liegroupsub(V)$ by
		right translations is proper, 
		see e.g.~\cite[Prop. 9.29]{gallier.quaintance:2020a}.
		
		The map $\overline{\pi} \times \id_{\liegroupsub(V)}$
		is a surjective
		submersion and
		$(\overline{\pi} \times \id_{\liegroupsub(V)}) \circ \widetilde{\racts}$ is
		smooth as the composition of smooth maps. Thus the
		action $\racts$ is smooth
		since the diagram
		\eqref{equation:lemma_GxHOm_principal_fiber_bundle_structure_diagram}
		commutes.
		
		Next let $[p, A] \in P \times_H \liegroupsub(V)$ and
		$B \in \liegroupsub(V)$.
		Then
		\begin{equation*}
			[p, A] \racts B = [p, A \circ B] = [p, A] \implies B = \id_{V} 
		\end{equation*}
		holds proving that $\racts$ is free.
		
		We now show that $\racts$ is proper.
		To this end,
		we use the characterization of a proper Lie group action
		in terms of sequences,
		see e.g. \cite[Sec. 6.20]{michor:2008}.
		
		Let $([p_i, A_i])_{i \in \field{N}}$ be a convergent sequence
		in $P  \times_H \liegroupsub(V)$ with limit
		$[p, A] \in P  \times_H \liegroupsub(V)$.
		Next let $(B_i)_{i \in \field{N}}$ be a sequence in $\liegroupsub(V)$
		such that the sequence defined by
		$[p_i, A_i] \racts B_i = [p_i, A_i \circ B_i]$ converges. 
		Then the action $\racts$ is proper iff $(B_i)_{i \in \field{N}}$ has a
		convergent subsequence.
		Let $s \colon U \to P \times \liegroupsub(V)$ be a 
		local section of the $H$-principal fiber bundle
		$\overline{\pi} \colon P \times \liegroupsub(V) \to P \times_H \liegroupsub(V)$
		defined on some open $U \subseteq P \times_H \liegroupsub(V)$
		such that 
		$[p, A] \in U$ holds. Then $[p_i, A_i] \in U$ is fulfilled for all
		$i \geq N$ with $N \in  \field{N}$ sufficiently large. 
		We define the sequence
		$(\widehat{p}_i, \widehat{A}_i)_{i \in \field{N}}$ in
		$P \times \liegroupsub(V)$ by setting 
		\begin{equation*}
			(\widehat{p}_i, \widehat{A}_i) = s([p_i, A_i]), 
			\quad i \geq N
		\end{equation*}
		and choosing
		$(\widehat{p}_i, \widehat{A}_i) \in {\overline{\pi}}^{-1}([p_i, A_i])$
		for $i < N$ arbitrarily. By construction, we have
		\begin{equation}
			\label{equation:lemma_PxHGV_principal_fiber_bundle_structure_pA_phat_Ahat_relation}
			[p_i, A_i] 
			=
			(\overline{\pi} \circ s)([p_i, A_i]) 
			=
			\overline{\pi}(\widehat{p}_i, \widehat{A}_i)
			=
			[\widehat{p}_i, \widehat{A}_i]
		\end{equation}
		for all $i \in \field{N}$.
		Moreover, the sequence
		$(\widehat{p}_i, \widehat{A}_i)_{i \in \field{N}}$
		converges to
		\begin{equation*}
			(\widehat{p}, \widehat{A}) 
			= \lim_{i \to \infty}s([p_i, A_i]) 
			= s([p, A]) \in P \times \liegroupsub(V)
		\end{equation*}
		by the continuity of the local section
		$s \colon U \to P \times \liegroupsub(V)$
		and the convergence of $[p_i, A_i]$.
		Moreover, let $(B_i)_{i \in \field{N}}$ be a sequence in $\liegroupsub(V)$
		such that the sequence defined by
		\begin{equation*}
			[p_i, A_i] \racts B_i = [p_i, A_i \circ B_i],
			\quad i \in \field{N}
		\end{equation*}
		is convergent in $P  \times_H \liegroupsub(V)$.
		We denote its limit by
		$[p, C] = \lim_{i \to \infty} [p_i, A_i \circ B_i] 
		\in P \times_H \liegroupsub(V)$.
		Clearly,
		\begin{equation*}
			[p_i, A_i] \racts B_i = [\widehat{p}_i, \widehat{A}_i] \racts B_i = [\widehat{p}_i, \widehat{A} \circ B_i],
			\quad i \in \field{N}
		\end{equation*}
		holds
		by~\eqref{equation:lemma_PxHGV_principal_fiber_bundle_structure_pA_phat_Ahat_relation}.
		Next we choose a local sections
		$s_2 \colon U_2 \to P \times \liegroupsub(V)$
		of $\overline{\pi} \colon P \times \liegroupsub(V) \to P \times_H \liegroupsub(V)$
		such that
		$[p, C] \in U_2 \subseteq P \times \liegroupsub(V)$
		is fulfilled.
		Then there exists an $N_2 \in \field{N}$ 
		with $[p_i, A_i \circ B_i] \in U_2$ 
		for all $i \geq N_2$.
		We define the sequence 
		\begin{equation}
			\label{equation:lemma_GxHOm_principal_fiber_bundle_structure_tilde_p_i_C_i}
			(\widetilde{p}_i, \widetilde{C}_i) 
			=
			s_2([\widehat{p}_i, \widehat{A}_i \circ B_i]), 
			\quad
			i \geq \widetilde{N}
		\end{equation}
		and select
		$(\widetilde{p}_i, \widetilde{C}_i) \in \overline{\pi}^{-1}([p_i, \widehat{A}_i \circ B_i])$ 
		for $i < N_2$ arbitrarily.
		Recall from~\cite[Sec. 18, p. 211]{michor:2008}
		that the map 
		\begin{equation*}
			\sigma \colon P \oplus P \to H,
			\quad
			(p, p^{\prime}) \mapsto \sigma(p, p^{\prime})
		\end{equation*}
		is smooth, where $\sigma(p, p^{\prime}) \in H$ is
		defined by
		$p \racts	\sigma(p, p^{\prime}) = p^{\prime}$
		for $(p, p^{\prime}) \in P \oplus P$.
		Next we define the map
		\begin{equation*}
			\Theta \colon (P  \oplus P )  \times \liegroupsub(V)\to  P \times \liegroupsub(V),
			\quad
			\big( (p, p^{\prime} ), A \big)
			\mapsto 
			\big(p , \rho_{\sigma(p, p^{\prime})^{-1}} \circ A \big)
		\end{equation*}
		which is a smooth map as the composition of smooth maps.
		The definition of
		$(\widetilde{p}_i, \widetilde{C}_i)_{i \in \field{N}}$
		in~\eqref{equation:lemma_GxHOm_principal_fiber_bundle_structure_tilde_p_i_C_i}
		implies
		\begin{equation}
			\label{equation:lemma_GxHOm_principal_fiber_bundle_structure_definition_tilde_p_i_tilde_C_i}
			\begin{split}
				(\widetilde{p}_i, \widetilde{C}_i) 
				&=
				s_2([\widehat{p}_i, \widehat{A}_i \circ B_i]) \\
				&=
				\big(\widehat{p}_i \racts \sigma(\widehat{p}_i, \widetilde{p}_i), \rho_{(\sigma(\widehat{p}_i, \widetilde{p}_i))^{-1}} \circ \widehat{A}_i \circ B_i \big)  \\
				&=
				\big( \widetilde{p}_i, \rho_{(\sigma(\widehat{p}_i, \widetilde{p}_i))^{-1}} \circ \widehat{A}_i \circ B_i \big) 
			\end{split}
		\end{equation}
		since
		$[\widetilde{p}_i, \widetilde{C}_i] = [\widehat{p}_i, \widehat{A}_i \circ B_i]$
		holds iff there exists
		a $h_i \in H$ with
		$\widetilde{p}_i = \widehat{p}_i \racts h_i$ and
		$\widetilde{C}_i = \rho_{h_i^{-1}} \circ \widehat{A}_i \circ B_i$.
		Next consider the sequence defined for $i \in \field{N}$ by
		\begin{equation}
			\label{equation:lemma_GxHOm_principal_fiber_bundle_structure_definition_another_sequence}
			\Theta\big((\widetilde{p}_i, \widehat{p}_i), \widehat{A}_i \big)
			=
			(\widehat{p}_i \racts \sigma(\widehat{p}_i, \widetilde{p}_i), \rho_{ (\sigma(\widehat{p}_i, \widetilde{p}_i))^{-1}} \circ \widehat{A}_i)
			=
			(\widetilde{p}_i, \rho_{ (\sigma(\widehat{p}_i, \widetilde{p}_i))^{-1}} \circ \widehat{A}_i), 
		\end{equation}
		which converges by the continuity of $\Theta$
		as well as the convergence of the sequences
		$(\widetilde{p}_i, \widetilde{C}_i)_{i \in \field{N}}$ and
		$(\widehat{p}_i, \widehat{A}_i)_{i \in \field{N}}$
		in $P \times \liegroupsub(V)$.
		By~\eqref{equation:lemma_GxHOm_principal_fiber_bundle_structure_definition_another_sequence}, we obtain
		\begin{equation*}
			\begin{split}
				\Theta\big((\widetilde{p}_i, \widehat{p}_i), \widehat{A}_i \big) \widetilde{\racts} B_i
				&=
				\big(\widetilde{p}_i, \rho_{ (\sigma(\widehat{p}_i, \widetilde{p}_i))^{-1}} \circ \widehat{A}_i \big) \widetilde{\racts}  B_i \\
				&=
				\big(\widetilde{p}_i, \rho_{ (\sigma(\widehat{p}_i, \widetilde{p}_i))^{-1}} \circ \widehat{A}_i \circ B_i \big) \\
				&=
				(\widetilde{p}_i, \widetilde{C_i}) \to (\widetilde{p}, \widetilde{C}),
				\quad i \to \infty ,
			\end{split}
		\end{equation*}
		where we
		used~\eqref{equation:lemma_GxHOm_principal_fiber_bundle_structure_definition_tilde_p_i_tilde_C_i}
		to obtain last equality. 
		Since the action
		$\widetilde{\racts} \colon (P \times \liegroupsub(V) ) \times \liegroupsub(V) \to P \times \liegroupsub(V)$
		is proper and the sequence
		$\big(	\Theta((\widetilde{p}_i, \widehat{p}_i), \widehat{A}_i)\big)_{i \in \field{N}}$
		defined
		in~\eqref{equation:lemma_GxHOm_principal_fiber_bundle_structure_definition_another_sequence}
		is convergent in $P \times \liegroupsub(V)$,
		the sequence
		$(B_i)_{i \in \field{N}}$ has a convergent subsequence.
		Thus the right action
		$\racts \colon (P \times_H \liegroupsub(V)) \times \liegroupsub(V) \to P \times_H \liegroupsub(V)$
		is indeed proper.
		Therefore
		$P \times_H \liegroupsub(V) \to (P \times_H \liegroupsub(V)) / \liegroupsub(V)$ 
		is a principal fiber bundle
		by~\cite[Re. 1.1.2]{rudolph.schmidt:2017}.

		It remains to prove
		Caim~\ref{item:lemma_PxHGV_principal_fiber_bundle_structure_isomorphism_base}.
		We first show that $\phi$ is a diffeomorphism.
		The equivalence classes
		$\widetilde{\pi}([p, A]) ) 
		=
		\widetilde{\pi}([p, A \circ B]) \in 
		\big( P \times_H \liegroupsub(V)\big) / \liegroupsub(V)$
		represented
		by $[p, A], [p, A \circ B] \in P \times_H \liegroupsub(V)$
		are equal, where
		$p \in P$, $A, B \in \liegroupsub(V)$.
		Thus we have
		\begin{equation*}
			\phi\big(\widetilde{\pi}([p,  A \circ B ]) \big) 
			=
			\pr(p)
			=
			\phi\big(\widetilde{\pi}([p, A])\big)
		\end{equation*}
		showing that $\phi$ is well-defined.
		We now consider the diagrams
		\begin{equation}
			\label{equation:lemma_GxHOm_principal_fiber_bundle_structure_diagram2}
			\begin{tikzpicture}[node distance= 5.0cm, auto]
				\node (GxOmOm) {$P \times_H \liegroupsub(V)$};
				\node (GxOmOm2) [right=3.0cm of GxOmOm] {$P \times_H \liegroupsub(V)$};
				\node (GxHOmOm) [below= 1.5cm of GxOmOm] {$\big(P \times_H \liegroupsub(V) \big) / \liegroupsub(V)$};
				\node (GxHOm) [below=1.583cm of GxOmOm2] {$M$};
				\draw[->] (GxOmOm) to node[left] {$\widetilde{\pi}$} (GxHOmOm);	
				\draw[->] (GxHOmOm) to node  [above]{$\phi$} (GxHOm);
				\draw[->] (GxOmOm2) to node [right] {$\pi$} (GxHOm);
				\draw[->] 	(GxOmOm) to node [above]{$\id_{(P \times_H \liegroupsub(V))}$}  (GxOmOm2);		
			\end{tikzpicture} 
		\end{equation}
		and
		\begin{equation}
			\label{equation:lemma_GxHOm_principal_fiber_bundle_structure_diagram3}
			\begin{tikzpicture}[node distance= 5.0cm, auto]
				\node (GxOmOm) {$P \times_H \liegroupsub(V)$};
				\node (GxOmOm2) [right=3.0cm of GxOmOm] {$P \times_H \liegroupsub(V)$};
				\node (GxHOmOm) [below= 1.5cm of GxOmOm] {$\big(P \times_H \liegroupsub(V) \big)  / \liegroupsub(V)$};
				\node (GxHOm) [below=1.583cm of GxOmOm2] {$M$};
				\draw[->] (GxOmOm) to node[left] {$\widetilde{\pi}$} (GxHOmOm);	
				\draw[<-] (GxHOmOm) to node  [above]{$\phi^{-1}$} (GxHOm);
				\draw[->] (GxOmOm2) to node [below, right] {$\pi$} (GxHOm);
				\draw[<-] 	(GxOmOm) to node [above]{$\id_{(P \times_H \liegroupsub(V))}$}  (GxOmOm2);		
			\end{tikzpicture} ,
		\end{equation}
		where $\phi^{-1}$ is given by 
		\begin{equation*}
			\phi^{-1} \colon M \ni \pr(p) \mapsto \widetilde{\pi}([p, \id_{V}]) \in  \big(P \times_H \liegroupsub(V) \big) / \liegroupsub(V) .
		\end{equation*}
		Clearly $\phi^{-1} \circ \phi = \id_{M}$ and
		$\phi \circ \phi^{-1} 
		=
		\id_{(P \times_H \liegroupsub(V))/ \liegroupsub(V)}$
		holds showing that $\phi$ is bijective.
		In addition, $\phi$ and $\phi^{-1}$ are smooth
		since~\eqref{equation:lemma_GxHOm_principal_fiber_bundle_structure_diagram2}
		and~\eqref{equation:lemma_GxHOm_principal_fiber_bundle_structure_diagram3}
		commute and $\widetilde{\pi}$ as well as $\pi$ are both
		surjective submersions.
		Hence the commutativity of \eqref{equation:lemma_GxHOm_principal_fiber_bundle_structure_diagram2}
		implies that $\id_{P \times_H \liegroupsub(V)}$ is indeed an isomorphism
		of $\liegroupsub(V)$-principal 
		fiber bundles over $\phi$ as desired.	
	\end{proof}
\end{lemma}

\begin{remark}
	By Lemma \ref{lemma:PxHGV_principal_fiber_bundle_structure},
	Claim~\ref{item:lemma_PxHGV_principal_fiber_bundle_structure_isomorphism_base},
	we view 
	$P \times_H \liegroupsub(V) \to \big(P \times_H \liegroupsub(V) \big) / \liegroupsub(V) \cong M$ 
	as an $\liegroupsub(V)$-principal fiber bundle over $M$
	which is denoted
	by the same symbol as the associated bundle, i.e.
	from now on, we write 
	\begin{equation}
		\pi \colon P \times_H \liegroupsub(V) \to M ,
	\end{equation} 
	if we view $P \times_H \liegroupsub(V)$ as an
	$\liegroupsub(V)$-principal 
	fiber bundle over $M$.
\end{remark}
Under certain conditions on the representation
$\rho \colon H \to \liegroup{GL}(V)$, one can determine
a reduction of the $\liegroup{GL}(V)$-principal fiber bundle 
$P \times_H\liegroup{GL}(V) \to M$ obtained
by setting $\liegroupsub(V) = \liegroup{GL}(V)$
in
Lemma~\ref{lemma:PxHGV_principal_fiber_bundle_structure}.

\begin{corollary}
	\label{corollary:reduction_frame_bundle_embedding}
	Let $P \to M$ be a $H$-principal fiber bundle and let
	$\rho \colon H \to  \liegroup{GL}(V)$ be a representation of $H$ on $V$.
	Moreover, let $\liegroupsub(V) \subseteq \liegroup{GL}(V)$ be
	a closed subgroup such that $\rho_h \in \liegroupsub(V)$ holds
	for all $h \in H$.
	Then
	\begin{equation}
		\iota_{P \times_H \liegroupsub(V)} \colon P \times_H \liegroupsub(V)
		\to 
		P \times_H \liegroup{GL}(V),
		\quad
		[p,A] \mapsto [p, A]
	\end{equation}
	is a reduction of the $\liegroup{GL}(V)$-principal fiber bundle
	$\pi_{P \times \liegroup{GL}(V)} \colon P \times_H \liegroup{GL}(V) \to M$ along the canonical
	inclusion $\liegroupsub(V) \to \liegroup{GL}(V)$.
	\begin{proof}
		Let $\overline{\iota_{P \times_H \liegroupsub(V)}} \colon P \times \liegroupsub(V) \to P \times \liegroup{GL}(V)$
		denote the canonical inclusion.
		Consider the diagram
		\begin{equation}
			\label{equation:corollary_reduction_frame_bundle_embedding_diagram}
			\begin{tikzpicture}[node distance= 5.0cm, auto]
				\node (PxGV) {$P \times \liegroupsub(V)$};
				\node (PxGLV) [right=3.0cm of PxGV] {$P \times \liegroup{GL}(V)$};
				\node (PxHGV) [below= 1.5cm of PxGV] {$P \times_H \liegroupsub(V)$};
				\node (PxHGLV) [below=1.5cm of PxGLV] {$P \times_H \liegroup{GL}(V)$};
				\draw[->] (PxGV) to node {$\overline{\iota_{P \times_H \liegroupsub(V)}}$} (PxGLV);	
				\draw[->] (PxGV) to node  [left]{$\overline{\pi}_{P \times \liegroupsub(V)}$} (PxHGV);
				\draw[->] (PxGLV) to node [below, right] {$\overline{\pi}_{P \times \liegroup{GL}(V)}$} (PxHGLV);
				\draw[->] 	(PxHGV) to node [above]{$\iota_{P \times_H \liegroupsub(V)}$}  (PxHGLV);		
			\end{tikzpicture} 
		\end{equation}
		which clearly commutes.
		Since the map $\overline{\pi}_{P \times \liegroupsub(V)}$
		is a surjective submersion and 
		$\overline{\pi}_{\liegroup{GL}(V)} \circ \overline{\iota_{P \times \liegroupsub(V)}}$ is smooth as
		the composition of smooth maps,
		the map $\iota_{P \times_H \liegroupsub(V)}$ is smooth, as well,
		because~\eqref{equation:corollary_reduction_frame_bundle_embedding_diagram}
		commutes.
		Clearly, the map $\iota_{P \times_H \liegroup{G}(V)}$ covers the
		map $\id_M \colon M \to M$.
		We now compute for $[p, A] \in P \times_H \liegroupsub(V)$ and
		$B \in \liegroupsub(V)$
		\begin{equation*}
			\iota_{P \times_H \liegroupsub(V)}([p, A] \racts B)
			= 
			[p, A] \racts B 
			= 
			\big(\iota_{P \times_H \liegroupsub(V)}([p, A] ) \big)\racts B
		\end{equation*}
		showing that $\iota_{P \times_H \liegroupsub(V)}$ is a morphism of
		principal fiber bundles along the canonical inclusion
		$\liegroupsub(V) \to \liegroup{GL}(V)$ covering
		$\id_M$,
		i.e. $\iota_{P \times_H \liegroupsub(V)}$ a reduction of $P \times_H \liegroup{GL}(V) \to M$.
	\end{proof}
\end{corollary}
The next proposition shows that
$\pi \colon P \times_H \liegroup{GL}(V) \to M$
can be identified with the frame bundle of the
associated vector bundle 
$P \times_H V \to M$, where $H$ acts on $V$ via
the representation viewed as the left action
\begin{equation}
	\rho \colon H  \times V \to V, 
	\quad
	(h, v) \mapsto \rho_h(v) .
\end{equation}

\begin{proposition}
	\label{proposition:frame_bundle_of_associated_vector_bundle}
	Let $\pr \colon P \to M$ be a $H$-principal fiber bundle and let 
	$\rho \colon H \to \liegroup{GL}(V)$ be a representation of $H$.
	The frame bundle of the associated vector bundle
	$P \times_H V \to M$ is isomorphic to $P \times_H \liegroup{GL}(V) \to M$ 
	as $\liegroup{GL}(V)$-principal fiber bundle via the isomorphism
	\begin{equation}
		\Psi \colon P \times_H \liegroup{GL}(V) \to \liegroup{GL}(V, P \times_H V),
		\quad
		[p, A] \mapsto \Psi([p, A])
	\end{equation}
	covering $\id_M$,
	where, for fixed $[p, A] \in P \times \liegroup{GL}(V)$,
	the linear isomorphism
	$\Psi([p, A]) \colon \{\pr(p)\} \times V \cong V \to (P \times_H V)_{\pr(p)}$
	is given by
	\begin{equation}
		\big( \Psi([p, A])\big)(v) = [p, A v] 
	\end{equation}
	for all $v \in V$.
	Here we view $\liegroup{GL}(V, P \times_H V)$ as an open subset of
	the morphism bundle $\Hom(M \times V, P \times_H V) \to M$
	as in~\cite[Sec. 18.11]{michor:2008}.
	Moreover, we write
	$\big( \Psi([p, A])\big)(v) = \big(\Psi([p, A])\big)(\pr(p), v)$
	for short, i.e. we suppress the first component 
	$\pr(p) \in M$ of $(\pr(p), v) \in M \times V$
	in the notation. 
	\begin{proof}
		We start with showing that $\Psi$ is well-defined.
		Let $h \in H$.
		Indeed, $\Psi$ is independent of the chosen representative
		of $[p, A] \in P \times_H \liegroup{GL}(V)$ due to
		\begin{equation*}
			\big( \Psi([p \racts h, \rho_{h^{-1}} \circ A])\big)(v)
			=
			[p \racts h, (\rho_{h^{-1}} \circ A)(v)]
			=
			[p, A v] 
		\end{equation*}
		for all $v \in V$.
		Moreover, for fixed $[p, A] \in P \times_H \liegroup{GL}(V)$, the map
		\begin{equation*}
			V \ni v \mapsto \Psi([p, A])(v) \in (P \times_H V)_{\pr(p)} 
		\end{equation*}
		is clearly linear.
		In addition, this map is invertible and its inverse is given by
		\begin{equation*}
			\big(\Psi([p, A]) \big)^{-1} \colon (P \times_H V)_{\pr(p)}
			\ni ([p, v])
			\mapsto 
			\big(\Psi([p, A] \big)^{-1}([p, v])
			=
			A^{-1} v 
			\in 
			V.
		\end{equation*}
		Indeed, $\big(\Psi([p, A]) \big)^{-1} $ is well-defined.
		Let $h, h^{\prime} \in H$.
		Then one has
		$p \racts h^{\prime} 
		=
		\big(p \racts (h h^{-1}) \big) \racts h^{\prime} 
		=
		(p \racts h) \racts (h^{-1} h^{\prime})$.
		Thus we obtain 
		\begin{equation*}
			\begin{split}
				\big(\Psi(p \racts h, \rho_{h^{-1}} \circ A) \big)^{-1}([p \racts h^{\prime}, \rho_{{h^{\prime}}^{-1}}(v)]) 
				&=
				(\rho_{(h^{-1}h^{\prime})^{-1}} \circ A)^{-1} (\rho_{(h^{-1}h^{\prime})^{-1}}(v)) \\
				&=
				A^{-1} v \\
				&=
				\big( \Psi([p, A]) \big)^{-1}([p, v])
			\end{split}
		\end{equation*}
		for all $v \in V$
		showing that $\big(\Psi([p, A]) \big)^{-1}$ is well-defined.
		Moreover, one has
		\begin{equation*}
			\big(\Psi([p, A]) \circ \big(\Psi([p, A])\big)^{-1}\big)([p, v])
			= 
			\big( \Psi([p, A])\big)(A^{-1} v)
			=
			[p, A A^{-1} v]
			=
			[p, v]
		\end{equation*}
		as well as 
		\begin{equation*}
			\big((\Psi([p, A]))^{-1} \circ \Psi([p, A]) \big)(v)
			=
			\big(\Psi([p, A]) \big)^{-1}([p, A v])
			=
			A^{-1}( A v) 
			=
			v 
		\end{equation*}
		showing that $\Psi([p, A]) \colon V \to (P \times_H V)_{\pr(p)}$
		is a linear isomorphism
		for all $[p, A] \in P \times_H \liegroup{GL}(V)$.
		Thus
		$\Psi \colon P \times_H \liegroup{GL}(V) \to \liegroup{GL}(V, P \times_H V)$
		is well-defined.
		
		Next we show that $\Psi$ is a morphism of principal fiber bundles over $\id_M$.
		Clearly,
		$\id_M \circ \pi_{P \times_H \liegroup{GL}(V)}
		=
		\pr_{\liegroup{GL}(V, P \times_H V)}$ holds,
		i.e. $\Psi$ covers $\id_M$.
		
		We now show that $\Psi$ is smooth.
		To this end, let $P \times_H \End(V) \to M$ denote the vector bundle
		associated to the $H$-principal fiber bundle
		$\pr \colon P \to M$ with typical fiber $\End(V)$,
		where $H$ acts on
		$\End(V)$ via
		\begin{equation*}
			H \times \End(V) \ni (h, A) \mapsto \rho_h \circ A \in \End(V)
		\end{equation*}
		from the left.
		We now define the map
		\begin{equation*}
			\begin{split}
				\widetilde{\Psi} \colon P \times_H \End(V) &\to \Hom(M \times V, P \times_H V), \\
				\quad
				[p, A] &\mapsto \widetilde{\Psi}([p, A]) 
				= \big( (x, v) 
				\mapsto \big(\widetilde{\Psi}([p, A])\big)(x, v) = [p, A v]  \big) .
			\end{split}
		\end{equation*}
		An argument analogously to the one at the beginning of this proof,
		showing that $\Psi$
		is well-defined,
		proves
		that the map $\widetilde{\Psi}$
		is well-defined, i.e. $\widetilde{\Psi}$
		is independent of the representative
		$(p, A) \in P \times \End(V)$ of
		$[p, A] \in P \times_H \End(V)$
		and that $\widetilde{\Psi}$ takes values
		in $\Hom(M \times V, P \times_H V) \to M$.
		Next we show that $\widetilde{\Psi}$
		is a smooth morphism of vector bundles.
		To this end, we prove that
		\begin{equation}
			\label{equation:proposition_frame_bundle_of_associated_vector_bundle_smoothness_vb_morphism_induced_module_morphism}
			\overline{\Psi} \colon \Secinfty\big( P \times_H \End(V) \big) \to \Secinfty\big( \Hom(M \times V, P \times_H V)  \big),
			\quad
			s \mapsto \widetilde{\Psi} \circ s
		\end{equation}
		is $\Cinfty(M)$-linear.
		Then the  desired properties of $\widetilde{\Psi}$
		follow by~\cite[Lem. 10.29]{lee:2013}.

		We first show that $\overline{\Psi}$ is well-defined,
		i.e. that
		$\overline{\Psi}(s) \in \Secinfty\big( \Hom(M \times V, P \times_H V)  \big)$
		is a smooth section of $\Hom(M \times V, P \times_H V) \to M$
		for all $s \in  \Secinfty\big( P \times_H \End(V) \big)$.
		In other words, we have to show that for fixed
		$s \in \Secinfty\big( P \times_H \End(V) \big)$ the map
		$\overline{\Psi}(s)$ is a smooth vector bundle morphism
		$\overline{\Psi}(s) \colon M \times V \to P \times_H V$
		over $\id_M$.
		Obviously, $\overline{\Psi}(s)$ is fiber-wise linear and covers $\id_M$.
		It remains to prove the smoothness of $\overline{\Psi}(s)$.
		To this end, we proceed locally. 
		Let $x_0 \in M$ and
		let $U \subseteq M$ be open with $x_0 \in U$.
		Moreover, after shrinking $U$ if necessary, 
		let $\widetilde{U} \subseteq P \times_H \End(V)$
		be open with $s(x) \in \widetilde{U}$ for all $x \in U$
		such that there is a smooth local section
		$\overline{s}  \colon \widetilde{U} \to P \times \End(V)$
		of the $H$-principal fiber bundle
		$\overline{\pi}_{P \times \End(V)} 
		\colon P \times \End(V) \to P \times_H \End(V)$.
		Then
		$\overline{s} \circ s \colon U \to P \times \liegroup{GL}(V)$
		is smooth and
		$(\overline{s} \circ s)(x) = (p(x), A(x))$
		holds for all $x \in U$
		with some smooth maps $U \ni x \mapsto p(x) \in P$ and
		$U \ni x \mapsto A(x) \in \liegroup{GL}(V)$.
		Thus
		$s(x) 
		=
		(\overline{\pi}_{P \times \End(V)} \circ \overline{s} \circ s)(x)
		=
		[p(x), A(x)]$
		is fulfilled for all $x \in U$.
		By this notation, we obtain for $(x, v) \in U \times V$
		\begin{equation}
			\label{equation:proposition_frame_bundle_of_associated_vector_bundle_smoothness_vb_morphism_induced_module_morphism_well_defined_section}
			\big(\overline{\Psi}(s) \big)(x, v) 
			=
			[p(x), A(x) v]  
			=
			\overline{\pi}_{P \times V}
			\circ (\id_P \times e)
			\circ 
			((\overline{s} \circ s) \times \id_V)
			(x, v)
		\end{equation}
		with $e \colon \End(V) \times V \ni (A, v) \mapsto A v \in V$.
		Hence the map
		$\Psi(s) \at{U \times V}$
		is smooth as the composition of smooth maps
		by~\eqref{equation:proposition_frame_bundle_of_associated_vector_bundle_smoothness_vb_morphism_induced_module_morphism_well_defined_section}.
		Thus $\Psi(s)$ is smooth since $x_0 \in M$ is arbitrary.
		
		Next we prove the $\Cinfty(M)$-linearity of $\overline{\Psi}$.
		Let $s_1,  s_2 \in \Secinfty\big( P \times_H \End(V) \big)$ be two sections
		point-wise given by
		\begin{equation*}
			s_1(x) = [p(x), A_1(x)] \quad \text{ and } \quad s_2(x) = [p(x), A_2(x)], 
			\quad
			x \in M .
		\end{equation*}
		Here we assume without loss of generality that their first component
		is represented by the same element $p(x) \in P$ for all $x  \in M$.
		Moreover, let $f, g \in \Cinfty(M)$.
		By the vector bundle structure on associated vector bundles,
		see e.g.~\cite[Re. 1.2.9]{rudolph.schmidt:2017},
		we obtain for $(x, v) \in M \times V$
		\begin{equation*}
			\begin{split}
				\big( \widetilde{\Psi} \circ (f s_1 + g s_1) \big)(x, v)
				&=
				\big(\widetilde{\Psi}\big(f(x) [p(x), A_1(x)] + g(x) [p(x), A_2(x) ] \big) \big)(x, v) \\
				&=
				\big[p(x), \big(f(x) A_1(x) + g(x) A_2(x) \big) v \big] \\
				&=
				f(x) [p(x), A_1(x) v] + g(x) [p(x), A_2(x) v] \\
				&=			
				\big( f (\widetilde{\Psi} \circ s_1)\big)(x, v) + 
				\big( g (\widetilde{\Psi} \circ s_2)\big)(x, v)	
			\end{split}
		\end{equation*}
		showing the $\Cinfty(M)$-linearity
		of $\overline{\Psi}$
		by its definition
		in~\eqref{equation:proposition_frame_bundle_of_associated_vector_bundle_smoothness_vb_morphism_induced_module_morphism}.
		Hence $\widetilde{\Psi}$
		is indeed a smooth morphism of vector bundles
		by~\cite[Lem. 10.29]{lee:2013}.
		
		In order to prove the smoothness of $\Psi$, we consider the map
		\begin{equation*}
			i \colon P \times_H \liegroup{GL}(V) \to P \times_H \End(V),
			\quad
			[p, A] \mapsto [p, A]
		\end{equation*}
		whose smoothness can be proven analogously to the proof of
		Corollary~\ref{corollary:reduction_frame_bundle_embedding}
		by exploiting the smoothness of the canonical inclusion
		$P \times \liegroup{GL}(V) \to P \times \End(V)$.
		We now obtain for $[p, A] \in P \times_H \liegroup{GL}(V)$
		and $(x, v) \in M \times V$
		\begin{equation*}
			\big( (\widetilde{\Psi} \circ i)([p, A]) \big)(x, v)
			=
			\big( \widetilde{\Psi}([p, A])\big) (x, v)
			=
			[p, A v]
			=
			\big(\Psi([p, A]) \big) (x, v) .
		\end{equation*}
		Thus $\Psi = \widetilde{\Psi} \circ i$ is smooth as the composition
		of smooth maps.
		
		It remains to show that $\Psi$ is an isomorphism of
		$\liegroup{GL}(V)$-principal fiber bundles.
		To this end, we recall that the $\liegroup{GL}(V)$-action
		on $\liegroup{GL}(V, P \times_H V)$ is given by composition
		from the right, see e.g.~\cite[Sec. 18.11]{michor:2008}.
		Thus we have  for $[p, A] \in P \times_H \liegroup{GL}(V)$ and
		$B \in \liegroup{GL}(V)$ as well as $v \in V$
		\begin{equation*}
			\big(\Psi([p, A] \racts B) \big)(v)
			=
			\big(\Psi([p, A \circ B]) \big)(v)
			=
			[p, (A \circ B) v]
			=
			\Psi([p, A]) (B v)
			=
			\big(\Psi[p, A] \circ B \big)(v) .
		\end{equation*}
		proving that $\Psi$ is a morphism of $\liegroup{GL}(V)$-principal
		fiber bundles over $\id_M \colon M \to M$.
		Therefore it is an
		isomorphism of principal fiber bundles
		by~\cite [Prop. 9.23]{gallier.quaintance:2020a}.
	\end{proof}
\end{proposition}
Assuming that $P \times_H \liegroup{GL}(V)$ admits a reduction as
in Corollary~\ref{corollary:reduction_frame_bundle_embedding},
we obtain a reduction of $\liegroup{GL}(V, P \times_H V)$. 

\begin{corollary}
	\label{corollary:frame_bundle_associated_bundle_reduction}
	Let $P \to M$ be a $H$-principal fiber bundle and let
	$\rho \colon H \to \liegroup{GL}(V)$ be a representation of $H$ on $V$
	such that $\rho_h \in \liegroupsub(V)$ holds for all $h \in H$, 
	where $\liegroupsub(V) \subseteq \liegroup{GL}(V)$ is a closed subgroup.
	Moreover, let
	$\Psi \colon P \times_H \liegroup{GL}(V) 
	\to \liegroup{GL}(V, P \times_H V)$
	be the isomorphism of principal fiber bundles from
	Proposition~\ref{proposition:frame_bundle_of_associated_vector_bundle}
	and let
	$\iota_{P \times_H \liegroupsub(V)} \colon P \times_H \liegroupsub(V) 
	\to P \times_H \liegroup{GL}(V)$
	be the reduction of principal fiber bundles from
	Corollary~\ref{corollary:reduction_frame_bundle_embedding}.
	Then
	\begin{equation}
		\label{equation:corollary_frame_bundle_associated_bundle_reduction}
		P \times_H \liegroup{G}(V) \to \liegroup{GL}(V, P \times_H V),
		\quad
		[p, A] \mapsto \big(\Psi \circ \iota_{P \times_H \liegroupsub(V)} \big)([p, A])
	\end{equation}
	is a $\liegroupsub(V)$-reduction of the
	frame bundle $\liegroup{GL}(V, P \times_H V)$ along the
	canonical inclusion $\liegroupsub(V) \to  \liegroup{GL}(V)$.
	\begin{proof}
		Obviously, the map $\Psi \circ \iota_{P \times_H \liegroupsub(V)}$
		is smooth as the composition of smooth maps.
		Moreover, since $\Psi$ is an isomorphism of principal fiber bundles
		covering $\id_M$
		by Proposition~\ref{proposition:frame_bundle_of_associated_vector_bundle}
		and
		$\iota_{P \times_H \liegroupsub(V)}$ is a reduction of principal
		fiber bundles
		along the canonical inclusion 
		$\liegroupsub(V) \to \liegroup{GL}(V)$
		by
		Corollary~\ref{corollary:reduction_frame_bundle_embedding},
		one verifies by a straightforward computation
		that~\eqref{equation:corollary_frame_bundle_associated_bundle_reduction}
		is a reduction of principal fiber bundles along
		the canonical inclusion 
		$\liegroupsub(V) \to \liegroup{GL}(V)$.
	\end{proof}
\end{corollary}

\begin{corollary}
	\label{corollary:frame_bundle_associated_bundle_vb_isomorphism}
	Let $P \times_H V \to M$ be a vector bundle associated to $P \to M$,
	where $\rho \colon H \to \liegroup{GL}(V)$ is a representation.
	Moreover, let
	$E \to N$ be another vector bundle 
	and let $\Phi \colon P \times_H V \to E$ be an isomorphism
	of vector bundles covering the diffeomorphism $\phi \colon M \to N$.
	Then
	\begin{equation}
		\chi \colon P \times_H \liegroup{GL}(V) \to \liegroup{GL}(V, E),
		\quad
		[p, A] \mapsto \chi([p, A]) 
		=
		( \Phi \circ \Psi)([p, A])
	\end{equation}
	is an isomorphism of $\liegroup{GL}(V)$-principal fiber bundles
	over the diffeomorphism $\phi \colon M \to N$,
	where
	$\Psi \colon P \times_H \liegroup{GL}(V) \to \liegroup{GL}(V, P \times_H V)$
	denotes the isomorphism
	from
	Proposition~\ref{proposition:frame_bundle_of_associated_vector_bundle}.
	\begin{proof}
		Obviously, for fixed $[p, A] \in P \times_H \liegroup{GL}(V)$,
		the map
		$\chi([p, A]) = (\Phi \circ \Psi)([p, A])
		\colon V \to E_{\phi(\pr(p))}$
		is linear and invertible since $\Phi$ is an isomorphism
		of vector bundles.
		Hence $\chi$ is well-defined.
		Moreover, the map $\chi$ is smooth as the composition
		of the smooth maps $\Psi$ and $\Phi$.
		Its inverse is given by the composition of the smooth maps
		$\chi^{-1} = \Psi^{-1} \circ \Phi^{-1} \colon \liegroup{GL}(V, E) \to P \times_H \liegroup{GL}(V)$,
		i.e. $\chi^{-1}$ is clearly smooth, as well.
		Let $B \in \liegroup{GL}(V)$ and $[p, A] \in P \times_H \liegroup{GL}(V)$.
		Then 
		\begin{equation*}
			\chi([p, A] \racts B )(v)
			=
			(\Phi \circ \Psi)( [p, A \circ B])(v)
			=
			(\Phi \circ \Psi)([p, A])(B v)
			=
			(\chi([p, A]) \circ B)(v)
		\end{equation*}
		holds for all $v \in V$ by the definition of $\Psi$.
		Hence $\chi$ is an isomorphism of
		$\liegroup{GL}(V)$-principal fiber bundles
		which covers the diffeomorphism
		$\phi \colon M \to N$.
	\end{proof}
\end{corollary}

\subsection{Principal Fiber Bundles over Frame Bundles and Principal Connections}

Since the $\liegroupsub(V)$-principal fiber bundle
$\pi \colon P \times_H \liegroupsub(V) \to M$ is obtained as a
fiber bundle associated to
the $H$-principal fiber bundle $P \to M$, we have the
$H$-principal fiber bundle
$	\overline{\pi} \colon P \times \liegroupsub(V) \to P \times_H \liegroupsub(V)$ over $P \times_H \liegroupsub(V)$.
Given a principal connection on $P \to M$, we construct a principal
connection on
$\overline{\pi} \colon P \times \liegroupsub(V) \to P \times_H \liegroupsub(V)$.
This construction will be
applied to the configuration space of an intrinsic rolling
of a reductive homogeneous space in
Proposition~\ref{proposition:principal_H_fiber_bundle_over_configuration_space}
below.

\begin{proposition}
	\label{proposition:principal_H_fiber_bundle_over_frame_bundle}
	Let $\pr \colon P \to M$ be a $H$-principal fiber bundle
	and let $\rho \colon H \to \liegroup{GL}(V)$ be a representation of $H$ on the finite dimensional $\field{R}$-vector space $V$.
	Assume that there exists a closed subgroup
	$\liegroupsub(V) \subseteq \liegroup{GL}(V)$
	with Lie algebra $\liealgsub(V) \subseteq \liealg{gl}(V)$
	such that $\rho_h \in \liegroupsub(V)$
	holds for all $h \in H$.
	Moreover, let 
	\begin{equation}
		\rho^{\prime} 
		\colon \liealg{h} \ 
		\to
		\liealgsub(V) \subseteq \liealg{gl}(V),
		\quad
		\eta
		\mapsto 
		(T_e \rho) \eta = \rho^{\prime}_{\eta} \in \liealgsub(V) \subseteq \liealg{gl}(V)
	\end{equation}
	denote the induced morphism of Lie algebras.
	Consider the $H$-principal fiber bundle
	\begin{equation}
		\overline{\pi} \colon P \times \liegroupsub(V) \to P \times_H \liegroupsub(V)
	\end{equation}
	over the associated bundle $\pi \colon P \times_H \liegroupsub(V) \to M$,
	where $H$ acts on $\liegroupsub(V)$ via
	\begin{equation}
		H \times \liegroupsub(V) \ni (h, A) \mapsto \rho_h \circ A \in \liegroupsub(V) .
	\end{equation}
	Moreover,
	let $\mathcal{P} \in \Secinfty\big(\End(T P) \big)$ be a principal connection
	on $\pr \colon P \to M$
	with corresponding connection one-form $\omega \in \Secinfty(T^* M) \tensor \liealg{h}$.
	Then the following assertions are fulfilled:
	\begin{enumerate}
		\item
		\label{item:proposition_principal_H_fiber_bundle_over_frame_bundle_vertical_bundle}
		The vertical bundle $\Ver(P \times \liegroupsub(V)) \subseteq T (P \times \liegroupsub(V)) \cong T P \times T \liegroupsub(V)$
		is fiber-wise given by 
		\begin{equation}
			\label{equation:proposition_principal_H_fiber_bundle_over_configuration_space_vertical_bundle_description}
			\Ver(P \times \liegroupsub(V))_{(p, A)}
			=
			\big\{ \big(\tfrac{\D}{\D t} \big(p \racts \exp(t \eta) \big) \at{t = 0},
			- \rho^{\prime}_{\eta} \circ A \big) \mid \eta \in \liealg{h}
			\big\} 
		\end{equation}
		where $(p, A) \in P \times \liegroupsub(V)$.
		\item
		\label{item:proposition_principal_H_fiber_bundle_over_frame_bundle_principal_connection}
		Defining
		$\overline{\mathcal{P}} 
		\in \Secinfty\big( \End(T (P \times \liegroup{GL}(V)))\big)$
		for $(p, A) \in P \times \liegroupsub(V)$ and 
		$(v_p, v_A) \in T_{(p, A)} (P \times \liegroupsub(V))$ by
		\begin{equation}
			\label{equation:proposition_principal_H_fiber_bundle_over_configuration_space_def_principal_connection}
			\overline{\mathcal{P}}\at{(p, A)}(v_p, v_A)
			=
			(\mathcal{P}{\at{g}}(v_p), - \rho^{\prime}_{\omega\at{p}(v_p)} \circ A) 
		\end{equation}
		yields a principal connection on
		$\overline{\pi} \colon P \times \liegroupsub(V) \to P \times_H \liegroupsub(V)$
		with corresponding connection one-form
		$\overline{\omega} \in \Secinfty\big(T^* (P \times \liegroupsub(V)) \big) \tensor \liealg{h}$
		given by
		\begin{equation}
			\label{equation:proposition_principal_H_fiber_bundle_over_configuration_space_def_connection_one_form}
			\overline{\omega}\at{(p, A)}(v_p, v_A)
			= 
			\omega\at{p}(v_p) 
		\end{equation}
		for all $(p, A) \in P \times \liegroupsub(V)$
		and $(v_p, v_A) \in T_{(p, A)}( P \times \liegroupsub(V))$.
		\item
		\label{item:proposition_principal_H_fiber_bundle_over_configuration_space_g_component}
		Let
		$\overline{q} \colon I \ni t \mapsto \overline{q}(t) = (p(t), A(t)) 
		\in P \times \liegroupsub(V)$ be a curve which
		is horizontal with respect to the
		principal connection $\overline{\mathcal{P}}$.
		Then the curve $p \colon I \to P$ given by the first
		component of $\overline{q}$ is horizontal with respect to
		the principal connection
		$\mathcal{P}$ on $P \to M$.
	\end{enumerate}
	\begin{proof}
		First we recall that
		$\rho^\prime \colon \liealg{h} \to \liealgsub(V)$ 
		is indeed a morphism of Lie algebras,
		see e.g.~\cite[Lem. 4.13]{michor:2008}.
		Next we prove Claim~\ref{item:proposition_principal_H_fiber_bundle_over_frame_bundle_vertical_bundle}.
		To this end, we compute for $(p, A) \in P \times \liegroupsub(V)$
		\begin{equation*}
			\begin{split}
				\Ver(P \times \liegroupsub(V))_{(p, A)}
				&=
				\big\{ \tfrac{\D}{\D t}\big( (p, A)  \overline{\racts} \exp(t \eta) \big) \at{t = 0}
				\mid \eta \in \liealg{h} \big\} \\
				&=
				\big\{  \big( \tfrac{\D}{\D t} \big( p \racts (\exp(t \eta)) \big) \at{t = 0}, \tfrac{\D}{\D t} \big( \rho_{\exp(- t \eta)}  \circ A \big) \at{t = 0} \big)
				\mid \eta \in \liealg{h} \big\} \\
				&=
				\big\{ \big( \tfrac{\D}{\D t} \big( p \racts (\exp(t \eta))  \big) \at{t = 0}, - \rho^{\prime}_{\eta} \circ A \big) \mid \eta \in \liealg{h}
				\big\} 
			\end{split}
		\end{equation*}
		showing
		Claim~\ref{item:proposition_principal_H_fiber_bundle_over_frame_bundle_vertical_bundle}.
		
		We now prove
		Claim~\ref{item:proposition_principal_H_fiber_bundle_over_frame_bundle_principal_connection}.
		Obviously,
		$\overline{\mathcal{P}}  \in \Secinfty\big(\End(T (P \times \liegroupsub(V))) \big)$
		holds.	
		Next we show that $\overline{\mathcal{P}}$ is a projection,
		i.e. $\overline{\mathcal{P}}^2 = \overline{\mathcal{P}}$
		is fulfilled.
		By
		using
		the correspondence of $\mathcal{P}$ and
		$\omega$
		from~\eqref{equation:principal_connection_connection_one_form_correspondence}
		as well as $\mathcal{P}^2 = \mathcal{P}$,
		we calculate
		for $p \in P$ and $v_p \in T_p P$
		\begin{equation}
			\label{equation:proposition_principal_H_fiber_bundle_over_configuration_space_computation_principal_con_on_G_preparation}
			\begin{split}
				\omega\at{p}(\mathcal{P}\at{p}(v_p)) 
				&=
				\big(T_e (p \racts \cdot) \big)^{-1} \mathcal{P}\at{p} \big(
				\mathcal{P}\at{p}(v_p) \big) \\
				&=
				\big(T_e (p \racts \cdot) \big)^{-1} \mathcal{P}\at{p}(v_p) \\
				&= 
				\omega\at{p}(v_p) .
			\end{split}
		\end{equation}
		Using~\eqref{equation:proposition_principal_H_fiber_bundle_over_configuration_space_computation_principal_con_on_G_preparation}
		and $\mathcal{P}^2 = \mathcal{P}$,
		we have for $(p, A) \in P \times \liegroupsub(V)$
		and
		$(v_p, v_A) \in T_p P \times T_A \liegroupsub(V)$
		\begin{equation*}
			\begin{split}
				\overline{\mathcal{P}}\at{(p, A)}\big(\overline{\mathcal{P}}\at{p, A}(v_p, v_A) \big)
				&=
				\overline{\mathcal{P}}\at{(p, A)}
				\big(\mathcal{P}\at{p}(v_p), - \rho^{\prime}_{\omega\at{p}(v_p)} \circ A \big) \\
				&=
				\big(\mathcal{P}\at{p}\big(\mathcal{P}\at{p}(v_p) \big), - \rho^{\prime}_{\omega\at{p}(\mathcal{P}\at{p}(v_p))} \circ A \big) \\
				&= 
				\big(\mathcal{P}\at{p}(v_p), - \rho^{\prime}_{\omega\at{p}(v_p)} \circ A \big) \\
				&= 
				\overline{\mathcal{P}}\at{(p, A)}(v_p, A) ,
			\end{split}
		\end{equation*}
		proving that
		$\overline{\mathcal{P}}^2 = \overline{\mathcal{P}} \in \Secinfty\big(\End(T (P \times \liegroupsub(V))) \big)$
		is a projection.
		
		Moreover, $\image(\overline{\mathcal{P}}) = \Ver(P \times \liegroupsub(V))$
		holds by $\image(\mathcal{P}) = \Ver(P)$
		and the characterization of the vertical bundle
		in~\eqref{equation:proposition_principal_H_fiber_bundle_over_configuration_space_vertical_bundle_description}.
		
		We now show that $\overline{\mathcal{P}}$
		corresponds to $\overline{\omega}$.
		To this end, let $\eta \in \liealg{h}$ and denote by
		$\eta_{P \times \liegroupsub(V)} \in \Secinfty\big(T (P \times \liegroupsub(V)) \big)$
		the corresponding fundamental vector field associated to the
		$H$-principal action
		given by
		\begin{equation*}
			\eta_{P \times \liegroupsub(V)}(p, A) 
			= 
			\tfrac{\D}{\D t}
			\big( (p, A) \overline{\racts} \exp(t \eta) \big) \at{t = 0},
			\quad
			(p, A) \in P \times \liegroupsub(V) .
		\end{equation*}
		By this notation and the definition of $\overline{\omega}$
		in~\eqref{equation:proposition_principal_H_fiber_bundle_over_configuration_space_def_connection_one_form},
		we obtain
		\begin{equation}
			\label{equation:proposition_principal_H_fiber_bundle_over_configuration_space_P_omega_1}
			\begin{split}
				&\big(\overline{\omega}\at{(p, A)}(v_p, v_A)\big)_{P \times \liegroupsub(V)}(p, A) \\
				&=
				\tfrac{\D}{\D t} \big( (p, A) \overline{\racts}
				\exp\big(t \overline{\omega}\at{(p, A)}(v_p, v_A) \big) \big) \at{t = 0}  \\
				&=	
				\big(\tfrac{\D}{\D t} \big( p \racts \exp\big(t \omega\at{p}(v_p) \big) \big) \at{t = 0},
				\tfrac{\D}{\D t} \big( \rho_{\exp(- t \omega\at{p}(v_p) )} \circ A \big)\at{t = 0} \big) \\
				&=
				\big( \mathcal{P}\at{g}(v_g),  - \rho^{\prime}_{\omega\at{p}(v_p)} \circ A \big) \\
				&= 
				\overline{\mathcal{P}}\at{(u, A)}(v_p, v_A) .
			\end{split}
		\end{equation}
		Moreover, denoting by
		$\eta_P \in \Secinfty(T P)$ the fundamental vector
		field on $P$ defined by $\eta \in \liealg{h}$, as usual,
		we compute
		\begin{equation}
			\label{equation:proposition_principal_H_fiber_bundle_over_configuration_space_P_omega_2}
			\begin{split}
				\overline{\omega}\at{(p, A)}\big(\eta_{P \times \liegroupsub(V)}(p, A) \big)
				&=
				\overline{\omega}\at{(p, A)} \big(\tfrac{\D}{\D t} \big( p \racts \exp(t \eta) \big) \at{t = 0}, \tfrac{\D}{\D t} \big( \rho_{\exp(- t \eta)} \circ A \big)\at{t = 0} \big)  \\
				&=
				\omega\at{p}\big(\eta_P(p) \big) \\
				&=
				\eta 
			\end{split}
		\end{equation}
		since $\omega$, being the connection one-form associated to
		$\mathcal{P}$, fulfills $\omega(\eta_P) = \eta$ for
		all $\eta \in \liealg{h}$.
		Thus $\overline{\omega}$ is the connection one-form
		corresponding to the connection $\overline{\mathcal{P}}$
		due
		to~\eqref{equation:proposition_principal_H_fiber_bundle_over_configuration_space_P_omega_1}
		and~\eqref{equation:proposition_principal_H_fiber_bundle_over_configuration_space_P_omega_2}.
		In order to show that $\overline{\mathcal{P}}$ is a principal
		connection, we show that $\overline{\omega}$ has the desired
		equivarience-property. 
		By exploiting that
		$\omega \in \Secinfty(T^* G) \tensor \liealg{h}$ 
		is a principal connection one-form,
		we compute for $h \in H$
		\begin{equation*}
			\begin{split}
				\big((\cdot \overline{\racts} h)^* \overline{\omega} \big)
				\at{(p, A)} (v_p, v_A) 
				&=
				\overline{\omega}\at{(p, A)  \overline{\racts} h}
				\big( T_{(p, A)} ( \cdot \overline{\racts} h) (v_p, v_A) \big) \\
				&= \overline{\omega}\at{(p \racts h, \rho_{h^{-1}} \circ A)} 
				\big( T_p (\cdot \racts h) v_p, T_A (\rho_{h^{-1}} \circ (\cdot) ) v_A \big) \\
				&=
				\omega_{p \racts h} \big(T_p (\cdot \racts h) v_p \big) \\
				&=
				\Ad_{h^{-1}}\big(\omega\at{p}(v_p) \big) \\
				&=
				\Ad_{h^{-1}}\big(\overline{\omega}\at{(p, A)}(v_p, v_A) \big)
			\end{split}
		\end{equation*}
		as desired.
		
		It remains to show
		Claim~\ref{item:proposition_principal_H_fiber_bundle_over_configuration_space_g_component}.
		Let 
		$\overline{q} \colon I \ni t \mapsto \overline{q}(t) 
		= (p(t), A(t)) \in P \times \liegroupsub(V)$
		be horizontal with respect to $\overline{\mathcal{P}}$.
		Then
		\begin{equation*}
				0 
				=
				\overline{\mathcal{P}}\at{\overline{q}(t)}(\dot{\overline{q}}(t))
				=
				\big( \mathcal{P}\at{p(t)}(\dot{p}(t)), - \ad_{\omega\at{p(t)}( \dot{p}(t))} \circ A(t) \big)
		\end{equation*}
		holds.
		In particular, this implies $\mathcal{P}\at{p(t)}(\dot{p}(t)) = 0$.
		Hence $p \colon I \to P$ is horizontal with respect to
		the principal connection $\mathcal{P}$ on $P \to M$.
	\end{proof}
\end{proposition}

\subsection{Frame Bundles of Reductive Homogeneous Spaces}

We now consider (certain reductions) of the
frame bundle of a reductive homogeneous space
by applying Proposition~\ref{proposition:frame_bundle_of_associated_vector_bundle}
to the $H$-principal fiber bundle $\pr \colon G \to G / H$.
To this end, we recall
that the tangent bundle of a reductive homogeneous space $G / H$ with
reductive decomposition $\liealg{h} = \liealg{h} \oplus \liealg{m}$ is
isomorphic to the vector bundle $G \times_H \liealg{m} \to G / H$,
where $H$ acts on $\liealg{m}$ via
\begin{equation}
	H \times \liealg{m} \ni (h, X) \mapsto \Ad_h(X) \in \liealg{m} .
\end{equation}
This statement as well as the statement of
Corollary~\ref{corollary:reductive_homogeneous_space_frame_bundle}
below
seem to be well-known since they can be found
in~\cite[Ex. 2.7]{baum:2009}.
Moreover, exploiting that the isotropy representation
$H \ni h \mapsto T_{\pr(e)} \tau_h \in \liegroup{GL}(T_{\pr(e)} (G / H))$
is equivalent to the
representation $H \mapsto \Ad_h\at{\liealg{m}} \in \liegroup{GL}(\liealg{m})$,
see Lemma~\ref{lemma:isotropy_representation_equivalence},
one obtains that
\begin{equation}
	\label{equation:tantent_bundle_homogeneous_space_associated_bundle_isomorphism}
	G \times_H \liealg{m} \to T (G / H),
	\quad
	[g, X] \mapsto (T_g \pr \circ T_e \ell_g) X 
\end{equation}
is an isomorphism of vector bundles over $\id_{G / H}$
by adapting the proof in~\cite[Sec. 18.16]{michor:2008}.

\begin{corollary}
	\label{corollary:reductive_homogeneous_space_reduced_frame_bundle}
	Let $G / H$ be a reductive homogeneous space with reductive decomposition $\liealg{g} = \liealg{h} \oplus \liealg{m}$.
	Moreover, assume that $\Ad_h\at{\liealg{m}} \in \liegroupsub(\liealg{m})$
	holds for all $h \in H$, where $\liegroupsub(\liealg{m})$ is some closed
	subgroup of $\liegroup{GL}(\liealg{m})$.
	Then 
	\begin{equation}
		\label{equation:corollary_reductive_homogeneous_space_reduced_frame_bundle_reduction_of_frame_bundle}
		G \times_H \liegroupsub(\liealg{m})
		\ni [g, A] \mapsto (X \mapsto [g, A X])  \in
		\liegroup{GL}(\liealg{m}, G \times_H \liealg{m})		
	\end{equation}
	is a reduction of the frame bundle of $G \times_H \liealg{m} \to G / H$
	along the
	canonical inclusion
	$\liegroupsub(\liealg{m}) \to \liegroup{GL}(\liealg{m})$.
	Moreover, the map
	\begin{equation}
		\label{equation:corollary_reductive_homogeneous_space_reduced_frame_bundle_reduction_of_frame_bundle_linear_isomorphism}
		G \times_H \liegroupsub(\liealg{m}) \to \liegroup{GL}(\liealg{m}, T(G / H)),
		\quad
		[g, A]
		\mapsto 
		\big(X \mapsto (T_g \pr \circ T_e \ell_g \circ A) X \big)
	\end{equation}
	is a reduction of $\liegroup{GL}(\liealg{m}, T(G / H))$
	along the canonical inclusion $\liegroupsub(\liealg{m}) \to \liegroup{GL}(\liealg{m})$.
	\begin{proof}
		The map defined
		in~\eqref{equation:corollary_reductive_homogeneous_space_reduced_frame_bundle_reduction_of_frame_bundle}
		is a reduction of the frame bundle of
		$G  \times_H \liealg{m} \to G / H$ by
		Proposition~\ref{proposition:frame_bundle_of_associated_vector_bundle}.
		It remains to show
		that~\eqref{equation:corollary_reductive_homogeneous_space_reduced_frame_bundle_reduction_of_frame_bundle_linear_isomorphism}
		is a reduction of principal fiber bundles.
		In fact, 
		\begin{equation}
			\label{equation:corollary_reductive_homogeneous_space_reduced_frame_bundle_reduction_iso_of_principal_fiber_bundles}
			G  \times_H \liegroup{GL}(\liealg{m}) \to \liegroup{GL}(\liealg{m}, T (G / H)),
			\quad
			[g, A] \mapsto
			\big( X \mapsto (T_g \pr  \circ T_e \ell_g \circ A) X  \big)
		\end{equation}
		is an isomorphism of principal fiber bundles covering $\id_{G / H}$
		by Corollary~\ref{corollary:frame_bundle_associated_bundle_vb_isomorphism}
		since~\eqref{equation:tantent_bundle_homogeneous_space_associated_bundle_isomorphism}
		is an isomorphism of vector bundles
		covering $\id_{G / H}$.
		The desired result follows by exploiting 
		that~\eqref{equation:corollary_reductive_homogeneous_space_reduced_frame_bundle_reduction_of_frame_bundle_linear_isomorphism}
		is the composition of the
		isomorphism~\eqref{equation:corollary_reductive_homogeneous_space_reduced_frame_bundle_reduction_iso_of_principal_fiber_bundles}
		and the
		reduction~\eqref{equation:corollary_reductive_homogeneous_space_reduced_frame_bundle_reduction_of_frame_bundle}.
	\end{proof}
\end{corollary}

\begin{remark}
	\label{remark:reduction_frame_bundel_G_H_identification_subset}
	In the sequel, under the assumption of Corollary~\ref{corollary:reductive_homogeneous_space_reduced_frame_bundle},
	we often identify $G \times_H \liegroupsub(\liealg{m})$ with the
	image of the
	reduction~\eqref{equation:corollary_reductive_homogeneous_space_reduced_frame_bundle_reduction_of_frame_bundle_linear_isomorphism}
	from Corollary~\ref{corollary:reductive_homogeneous_space_reduced_frame_bundle}
	as in~\cite[Re. 1.1.8]{rudolph.schmidt:2017}.
	This is indicated by the notation
	$\liegroupsub(\liealg{m}, T(G / H)) 
	\subseteq \liegroup{GL}(\liealg{m}, T (G / H))$.
\end{remark}

\begin{corollary}
	\label{corollary:reductive_homogeneous_space_pseuo_Riemannian_frame_bundle}
	Let $G / H$ be a pseudo-Riemannian reductive homogeneous space whose
	invariant metric corresponds to the $\Ad(H)$-invariant scalar product
	$\langle \cdot, \cdot \rangle \colon \liealg{m} \times \liealg{m} \to \field{R}$.
	Moreover,
	denote by $\liegroup{O}(\liealg{m})$ the pseudo-orthogonal
	group of $\liealg{m}$ with respect to
	$\langle \cdot, \cdot \rangle$.
	Then
	\begin{equation}
		\label{equation:corollary_reductive_homogeneous_space_pseuo_Riemannian_frame_bundle_reduction}
		G \times_H \liegroup{O}(\liealg{m}) \to \liegroup{GL}(\liealg{m}, T( G / H)),
		\quad
		[g, S] \mapsto  \big(X \mapsto (T_g \pr \circ T_e \ell_g \circ S) X \big)
	\end{equation}
	is a reduction of the frame bundle of $T (G / H)$ along the 
	canonical inclusion
	$\liegroup{O}(\liealg{m}) \to \liegroup{GL}(\liealg{m})$.
	\begin{proof}
		This is a consequence
		of
		Corollary~\ref{corollary:reductive_homogeneous_space_reduced_frame_bundle}.
	\end{proof}
\end{corollary}
If $G \times_H \liegroup{O}(\liealg{m})$ is identified with the image
of~\eqref{equation:corollary_reductive_homogeneous_space_pseuo_Riemannian_frame_bundle_reduction}, it
is often denoted by $\liegroup{O}(\liealg{m}, T(G / H))$.

\begin{corollary}
	\label{corollary:reductive_homogeneous_space_frame_bundle}
	Let $G / H$ be a reductive homogeneous space with reductive decomposition $\liealg{g} = \liealg{h} \oplus \liealg{m}$.
	Then the frame bundle of $T (G / H)$ is isomorphic to
	$G \times_H \liegroup{GL}(\liealg{m})  \to G / H$
	as $\liegroup{GL}(\liealg{m})$-principal fiber bundle via 
	the isomorphism
	\begin{equation}
		G \times_H \liegroup{GL}(\liealg{m})
		\to
		\liegroup{GL}(\liealg{m}, T (G / H)),
		\quad
		[g, A] \mapsto
		\big( X \mapsto ( T_g \pr \circ T_e \ell_g \circ A ) X \big) 
	\end{equation}
	of $\liegroup{GL}(\liealg{m})$-principal fiber bundles.
	\begin{proof}
		This follows by setting
		$\liegroupsub(\liealg{m}) = \liegroup{GL}(\liealg{m})$
		in
		Corollary~\ref{corollary:reductive_homogeneous_space_reduced_frame_bundle}.
	\end{proof}
\end{corollary}

\begin{remark}
	Corollary~\ref{corollary:reductive_homogeneous_space_frame_bundle}
	seems to be well-known since the statement that
	$G \times_H \liegroup{GL}(\liealg{m}) \to G / H$ is isomorphic to the frame bundle of $G / H$ can be found
	in~\cite[Ex. 2.7]{baum:2009}.
\end{remark}

\section{Intrinsic Rollings of Reductive Homogeneous Spaces}
\label{sec:intrinsic_rolling_reductive_space}

Let $G / H$ be a reductive homogeneous space with fixed
reductive decomposition $\liealg{g} = \liealg{h} \oplus \liealg{m}$.
In the sequel, we always endow $\liealg{m}$ with the covariant
derivative $\nabla^{\liealg{m}}$ which is defined
in~\eqref{equation:definition_covariant_derivative_m} below.
Let
$V \colon \liealg{m} \ni v \mapsto (v, V_2(v)) \in \liealg{m} \times \liealg{m} \cong T \liealg{m}$
and 
$W \colon \liealg{m} \ni v \mapsto (v, W_2(v)) \in \liealg{m} \times \liealg{m} \cong T \liealg{m}$
be vector fields on 
$\liealg{m}$, where $V_2, W_2 \colon \liealg{m} \to \liealg{m}$
are smooth maps.
Then
$\nabla^{\liealg{m}} \colon \Secinfty(T \liealg{m}) \times \Secinfty(T \liealg{m}) \to \Secinfty(T \liealg{m})$
is defined by
\begin{equation}
	\label{equation:definition_covariant_derivative_m}
	\nabla^{\liealg{m}}_V W \at{v}
	=
	\big( v, (T_v W_2) V_2(v) \big),
	\quad
	v \in \liealg{m} .
\end{equation}
Clearly, for $\liealg{m} = \field{R}^n$ the covariant derivative
$\nabla^{\liealg{m}}$ coincides with the covariant derivative
from~\cite[Chap. 3, Def. 8]{oneill:1983}.

In this section, we consider intrinsic
($\liegroupsub(\liealg{m})$-reduced) rollings of
$(\liealg{m}, \nabla^{\liealg{m}})$ 
over $G / H$ equipped with an invariant covariant
derivative $\nablaAlpha$.
Such intrinsic rollings are called rollings of $\liealg{m}$
over $G / H$ with respect to $\nablaAlpha$, 
rollings of $G / H$ with respect to $\nablaAlpha$,
or simply rollings of $G / H$, for short.

\begin{notation}
	In the sequel, we do not explicitly refer to the
	$\liegroupsub(\liealg{m})$-reduction if this reduction is clear
	by the context, for instance by denoting the configuration space
	by $Q = \liealg{m} \times (G \times \liegroupsub(\liealg{m}))$
	as in Lemma~\ref{lemma:configuration:space_intrinsic_rolling},
	below. 
\end{notation}

\subsection{Configuration Space}
\label{subsec:configuration-space}

The goal of this subsection is to derive an
explicit description of the configuration space 
for rollings of $\liealg{m}$ over $G / H$
with respect to an invariant covariant derivative
$\nablaAlpha$.
Moreover, we consider a
$H$-principal fiber bundle over the configuration space
equipped with a suitable principal connection.
This
allows for lifting rollings,
i.e. certain curves on the configuration space,
horizontally to curves on that principal fiber bundle.
We start with investigating
the configuration space.

\begin{lemma}
	\label{lemma:configuration:space_intrinsic_rolling}
	Let $G / H$ be a reductive homogeneous space
	and let $\liegroupsub(\liealg{m}) \subseteq \liegroup{GL}(\liealg{m})$
	be a closed subgroup such that
	$\Ad_h\at{\liealg{m}} \in \liegroupsub(\liealg{m})$ holds
	for all $h \in H$.
	Then the following assertions are fulfilled:
	\begin{enumerate}
		\item
		Let
		$\Qspace = \liealg{m} \times (G \times_H \liegroupsub(\liealg{m}))$
		and define
		\begin{equation}
			\prQ \colon \Qspace \to \liealg{m} \times G / H,
			\quad
			(v, [g, S]) \mapsto
			\prQ(v, [g, S]) = (v, \pr(g)) .
		\end{equation}
		Then $\prQ \colon \Qspace \to \liealg{m} \times G / H$ is 
		isomorphic to the
		configuration space of the
		intrinsic rolling of $\liealg{m}$ over $G / H$,
		i.e. to the $\liegroupsub(\liealg{m})$-fiber bundle
		\begin{equation}
			\begin{split}
				&\big( \liegroupsub(\liealg{m}, T \liealg{m}) \times \liegroupsub(\liealg{m}, T(G / H)) \big) / \liegroupsub(\liealg{m})  \\
				&\cong
				\big( (\liealg{m} \times \liegroupsub(\liealg{m})) \times (G \times_H \liegroupsub(\liealg{m})) \big) / \liegroupsub(\liealg{m}) \to
				\liealg{m} \times G / H
			\end{split}
		\end{equation} 
		via the isomorphism of $\liegroupsub(\liealg{m})$-fiber bundles
		\begin{equation}
			\label{equation:lemma_configuration:space_intrinsic_rolling_isomorphism}
			\begin{split}
				\Psi \colon 
				\big( (\liealg{m} \times \liegroupsub(\liealg{m})) \times (G \times_H \liegroupsub(\liealg{m})) \big) / \liegroupsub(\liealg{m})
				&\to
				\liealg{m} \times (G \times_H \liegroupsub(\liealg{m})), \\
				\big[ (v, S_1), [g, S_2] \big] 
				&\mapsto
				\big( v,  [g, S_2 \circ S_1^{-1}]\big)
			\end{split}
		\end{equation}
		covering the identity
		$\id_{\liealg{m} \times G / H}
		\colon \liealg{m} \times G / H \to \liealg{m} \times G / H$
		whose inverse is given by
		\begin{equation}
			\label{equation:lemma_configuration:space_intrinsic_rolling_isomorphism_inverse}
			\begin{split}
				\Psi^{-1} \colon \liealg{m} \times (G \times_H \liegroupsub(\liealg{m}))
				&\to 
				\big((\liealg{m} \times \liegroupsub(\liealg{m})) \times (G \times_H \liegroupsub(\liealg{m})) \big) / \liegroupsub(\liealg{m}) , \\
				(v, [g, S]) 
				&\mapsto 
				\big[ (v, \id_{\liealg{m}}), [g, S] \big] .
			\end{split}
		\end{equation}
		\item
		\label{item:lemma_configuration:space_intrinsice_rolling_associated_isometry}
		Let $q = (v, [g, S]) \in \Qspace$ with
		$\prQ(q) = (v, \pr(g))$.
		Then $q$ defines the linear isomorphism
		\begin{equation}
			T_v \liealg{m} \cong \liealg{m} 
			\ni Z \mapsto
			\big(T_g \pr \circ T_e \ell_g \circ S \big) Z
			\in T_{\pr(g)} G / H
		\end{equation}
		via Lemma~\ref{lemma:configuration_space_is_fiber_bundle},
		Claim~\ref{item:lemma_configuration_space_is_fiber_bundle_fiber_isomorphisms_identification},
		where $q$ is identified with
		$\Psi^{-1}(q) 
		\in 
		\big( (\liealg{m} \times \liegroupsub(\liealg{m}))
		\times (G \times_H \liegroupsub(\liealg{m})) \big) /
		\liegroupsub(\liealg{m})$.
		In the sequel, we often denote this isomorphism by $q$, as well, i.e. we write 
		$q(Z) = q Z = \big(T_g \pr \circ T_e \ell_g \circ S \big) Z$.
	\end{enumerate}
	\begin{proof}
		By
		Corollary~\ref{corollary:reductive_homogeneous_space_reduced_frame_bundle}
		we have $\liegroupsub(\liealg{m}, T(G / H)) \cong G \times_H \liegroupsub(\liealg{m})$.
		Moreover, $\liegroupsub(\liealg{m}, T \liealg{m}) \cong \liealg{m} \times \liegroupsub(\liealg{m})$ is clearly fulfilled.
		We first show that $\Psi$ is smooth.
		Consider 
		\begin{equation}
			\label{equation:lemma_configuration:space_intrinsic_rolling_diagram}
			\begin{tikzpicture}[node distance= 4.0cm, auto]
				\node (topQ) {$(\liealg{m} \times \liegroupsub(\liealg{m})) \times (G \times_H \liegroupsub(\liealg{m}) )$};
				\node (tildeQ) [below= 1.5cm of topQ] {$\big((\liealg{m} \times \liegroupsub(\liealg{m})) \times (G \times_H \liegroupsub(\liealg{m})) \big) / \liegroupsub(\liealg{m})$};
				\node (Q) [right=3.0cm of tildeQ] {$\liealg{m} \times (G \times_H \liegroupsub(\liealg{m}))$};
				\draw[->] (topQ) to node [left] {$\overline{\pr}$} (tildeQ);
				\draw[->] (topQ) to node  [pos=0.5]{$\overline{\Psi}$} 
				(Q);
				\draw[->] (tildeQ) to node [above] {$\Psi$} 
				(Q);
			\end{tikzpicture}  ,
		\end{equation}
		where $\overline{\pr}$ is the canonical projection and
		$\overline{\Psi}$ is given by
		\begin{equation}
			\overline{\Psi}\big( (v, S_1), [g, S_2] \big) 
			=
			(v, [g, S_2 \circ S_1^{-1}]) 
		\end{equation}
		for 
		$((v, S_1), [g, S_2])
		\in 
		(\liealg{m} \times \liegroupsub(\liealg{m})) \times (G \times_H \liegroupsub(\liealg{m}))$.
		Clearly,
		since~\eqref{equation:lemma_configuration:space_intrinsic_rolling_diagram}
		commutes and the canonical projection $\overline{\pr}$
		is a surjective submersion,
		the map $\Psi$ defined
		by~\eqref{equation:lemma_configuration:space_intrinsic_rolling_isomorphism}
		is smooth by the smoothness of $\overline{\Psi}$.
		In addition, $\Psi$ maps fibers into fibers,
		i.e. it is a morphism of
		$\liegroupsub(\liealg{m})$-fiber bundles
		covering the identity of $\liealg{m} \times G / H$.
		Therefore $\Psi$ is an isomorphism of fiber bundles, see 
		e.g.~\cite[Prop. 9.3]{gallier.quaintance:2020a}.
		The
		formula~\eqref{equation:lemma_configuration:space_intrinsic_rolling_isomorphism_inverse}
		for $\Psi^{-1}$
		is verified by a straightforward calculation.
		
		It remains to show
		Claim~\ref{item:lemma_configuration:space_intrinsice_rolling_associated_isometry}.
		Let $q = (v, [g, S]) \in Q 
		=
		\liealg{m} \times (G \times_H \liegroupsub(\liealg{m}))$
		and let $Z \in T_v \liealg{m}$.
		Then $\Psi^{-1}((v, [g, S])) = \big[(v, \id_{\liealg{m}}), [g, S] \big]$ holds.
		Using the bijection from
		Lemma~\ref{lemma:configuration_space_is_fiber_bundle}, Claim~\ref{item:lemma_configuration_space_is_fiber_bundle_fiber_isomorphisms_identification},
		this element
		is identified with a linear isomorphism
		which we denote by the same symbol.
		Evaluated at $Z \in T_v \liealg{m} \cong \liealg{m}$, it is given by
		\begin{equation*}
			\begin{split}
				\liniso{(\Psi^{-1}(v, [g, S]) \big)\big)}(Z)
				&=
				\liniso{\big(\big( \big[ (v, \id_{\liealg{m}}), [g, S] \big)\big)}(Z) \\ 
				&=
				\big( \big(T_g \pr \circ T_e \ell_g \circ S \big) \circ (\id_{\liealg{m}})^{-1} \big) (Z) \\
				&=
				\big(T_g \pr \circ T_e \ell_g \circ S \big) Z ,
			\end{split}
		\end{equation*}
		where the second equality follows by
		Lemma~\ref{lemma:configuration_space_is_fiber_bundle},
		Claim~\ref{item:lemma_configuration_space_is_fiber_bundle_fiber_isomorphisms_identification}
		and
		Corollary~\ref{corollary:reductive_homogeneous_space_reduced_frame_bundle}.
	\end{proof}
\end{lemma}

\begin{remark}
	The configuration space $\prQ \colon \Qspace \to \liealg{m} \times G / H$ 
	can be viewed as a $\liegroupsub(\liealg{m})$-principal fiber bundle.
	Indeed, as a consequence of Lemma~\ref{lemma:PxHGV_principal_fiber_bundle_structure},
	the $\liegroupsub(\liealg{m})$-right action
	\begin{equation}
		\Qspace \times \liegroupsub(\liealg{m}) \to \Qspace,
		\quad
		\big((v, [g, S]), S_2 \big) \mapsto (v, [g, S \circ S_2]) 
	\end{equation}
	is a principal action.
\end{remark}
Moreover, since the configuration space
$\Qspace = \liealg{m} \times (G \times_H \liegroupsub(\liealg{m}))$ 
is the product of $\liealg{m}$ and the associated bundle
$G \times_H \liegroupsub(\liealg{m})$,
we obtain an $H$-principal fiber bundle over $Q$.

\begin{proposition}
	\label{proposition:principal_H_fiber_bundle_over_configuration_space}
	Let $G / H$ be a reductive homogeneous space and let
	$\liegroupsub(\liealg{m}) \subseteq \liegroup{GL}(\liealg{m})$
	be a closed subgroup such that
	$\Ad_h\at{\liealg{m}} \in \liegroupsub({\liealg{m}})$
	holds for all $h \in H$.
	Then the following assertions are fulfilled:
	\begin{enumerate}
		\item
		\label{item:proposition:principal_H_fiber_bundle_over_configuration_space_definition_QLift}
		Define
		$\QspaceLift = \liealg{m} \times G \times \liegroupsub(\liealg{m})$.
		Then 
		\begin{equation}
			\label{equation:proposition:principal_H_fiber_bundle_over_configuration_space_definition_QLift_action}
			\prQLift \colon \QspaceLift 
			\ni (v, g, S)
			\mapsto
			(v, [g, S]) \in \Qspace
		\end{equation}
		becomes a $H$-principal fiber bundle over
		$\Qspace = \liealg{m} \times (G \times_H \liegroupsub(\liealg{m}))$
		with $H$-principal action given by
		\begin{equation}
			\begin{split}
				\racts_{\QspaceLift} \colon 
				\QspaceLift \times H
				&\to
				\QspaceLift , \\
				(v, g, S) 
				&\mapsto
				(v, g, S) \racts_{\overline{Q}} h 
				=
				(v, g h, \Ad_{h^{-1}} \circ S) 
				= (v, g \racts h,  \Ad_{h^{-1}} \circ S) ,
			\end{split}
		\end{equation}
		where
		$\racts \colon G \times H \ni (g, h) \mapsto g \racts h = g h \in G$
		denotes the $H$-principal action
		from~\eqref{equation:H_principal-action_on_G_mod_H} on $\pr \colon G \to G / H$.
		\item
		\label{item:proposition:principal_H_fiber_bundle_over_configuration_space_definition_vertical_bundle}
		For $(v, g, S) \in \QspaceLift = \liealg{m} \times G \times \liegroupsub(\liealg{m})$
		the vertical bundle $\Ver(\QspaceLift) \subseteq T \QspaceLift$
		is given by 
		\begin{equation}
			\Ver(\QspaceLift)_{(v, g, S)}
			=
			\big\{ (0, T_e \ell_g \eta, - \ad_{\eta} \circ S) \mid \eta \in \liealg{h} \big\} 
			\subseteq T_{(v, g, S)} (\liealg{m} \times G \times \liegroupsub(\liealg{m})) .	
		\end{equation}
		\item
		\label{item:proposition:principal_H_fiber_bundle_over_configuration_space_defi_connection}
		Let $\mathcal{P} \in \Secinfty\big(\End(T G) \big)$
		and let $\omega \in  \Secinfty(T^* G ) \tensor \liealg{h}$
		denote the principal connection 
		and connection one-form from
		Proposition~\ref{proposition:principal_connection_reductive_homogeneous_space}
		on $\pr \colon G \to G / H$, respectively.
		Defining for $(v, g, S) \in \QspaceLift$ 
		and $(u, v_g, v_S) \in T_{(v, g, S)} \QspaceLift$
		\begin{equation}
			\overline{\mathcal{P}}\at{(v, g, S)}(u, v_g, v_S)
			=
			(0, \mathcal{P}{\at{g}}(v_g), - \ad_{\omega\at{g}(v_g)} \circ S) ,
		\end{equation}
		yields a principal connection on $\prQLift \colon \QspaceLift \to \Qspace$
		with corresponding connection one-form
		$\overline{\omega} \in \Secinfty(T^* \QspaceLift) \tensor \liealg{h}$
		given by
		\begin{equation}
			\overline{\omega}\at{(v, g, S)}(u, v_g, v_S)
			= 
			\omega\at{g}(v_g) ,
			\quad
			(v, g, S) \in \QspaceLift, \ (u, v_g, v_S) \in T_{(v, g, S)} \QspaceLift .
		\end{equation}
		\item
		\label{item:proposition:principal_H_fiber_bundle_over_configuration_space_definition_horizontal_lift_curve}
		Let $\overline{q} \colon I \ni t \mapsto \overline{q}(t) = (v(t), g(t), S(t)) \in \QspaceLift$ be a horizontal curve
		with respect to the
		principal connection $\overline{\mathcal{P}}$.
		Then the curve $g \colon I \to G$ defined by the second
		component of $\overline{q}$ is horizontal with respect to
		$\Hor(G)$ from
		Proposition~\ref{proposition:principal_connection_reductive_homogeneous_space}.
	\end{enumerate}
	\begin{proof}
		We consider $\pr \colon G \to G / H$ as a $H$-principal fiber bundle.
		Then 
		$\liealg{m} \times G \ni (v, g) \mapsto (v, \pr(g))
		\in \liealg{m} \times G / H$
		becomes clearly a $H$-principal fiber bundle with principal action
		\begin{equation*}
			(\liealg{m} \times G) \times H \ni 
			((v, g), h) \mapsto (v, g \racts h) = (v, g h) 
			\in \liealg{m} \times G /H 
		\end{equation*}
		and
		$\Qspace = \liealg{m} \times (G \times_H \liegroupsub(\liealg{m}))$ 
		can be viewed as a
		$\liegroupsub(\liealg{m})$-fiber bundle associated to
		the $H$-principal fiber bundle $\liealg{m} \times G / H$,
		where $H$ acts on $\liegroupsub(\liealg{m})$ via	
		\begin{equation*}
			H \times \liegroupsub(V) \to \liegroupsub(V),
			\quad
			(h, S) \mapsto \Ad_h\at{\liealg{m}} \circ S .
		\end{equation*}
		Thus, by the definition of an associated bundle,
		$\prQLift \colon \QspaceLift \to \Qspace$
		becomes a $H$-principal fiber bundle over $\Qspace$ with
		principal action given
		by~\eqref{equation:proposition:principal_H_fiber_bundle_over_configuration_space_definition_QLift_action},
		i.e.
		Claim~\ref{item:proposition:principal_H_fiber_bundle_over_configuration_space_definition_QLift}
		is shown.

		Next, let $\mathcal{P} \in \Secinfty\big(\End(T P) \big)$ be
		the principal connection on $G$ from
		Proposition~\ref{proposition:principal_connection_reductive_homogeneous_space}.
		It is straightforward to verify that
		$\widetilde{\mathcal{P}} \in \Secinfty\big(\End(T P) \big)$
		defined by
		\begin{equation*}
			\widetilde{\mathcal{P}}\at{(v, g)}(u, v_g) = \big(0, \mathcal{P}(v_g) \big),
			\quad
			(v, g) \in \liealg{m} \times G, 
			\quad
			(u, v_g) \in \liealg{m} \times T_g G  \cong T_v \liealg{m} \times T_g G 
		\end{equation*}
		yields a principal connection on
		$\liealg{m} \times G \to \liealg{m} \times G / H$
		with corresponding connection one-form given by
		$\widetilde{\omega}\at{(v,g)}(u, v_g) = \omega\at{g}(v_g)$.
		Thus
		Proposition~\ref{proposition:principal_H_fiber_bundle_over_frame_bundle}
		applied to $\liealg{m} \times G  \to \liealg{m} \times G / H$
		equipped with the principal connection $\widetilde{\mathcal{P}}$
		yields Claim~\ref{item:proposition:principal_H_fiber_bundle_over_configuration_space_definition_vertical_bundle},
		Claim~\ref{item:proposition:principal_H_fiber_bundle_over_configuration_space_defi_connection},
		and 
		Claim~\ref{item:proposition:principal_H_fiber_bundle_over_configuration_space_definition_horizontal_lift_curve}.
	\end{proof}
\end{proposition}

\subsection{The Distribution characterizing Intrinsic Rollings}

Motivated by~\cite[Sec. 4]{molina.grong.markina.leite:2012},
we determine
a distribution on $\Qspace$ 
characterizing intrinsic rollings of $\liealg{m}$ over $G /H$.
More precisely, a curve $q \colon I \to \Qspace$ is horizontal with 
respect to this distributions iff it is a rolling of $G / H$ 
with respect to $\nablaAlpha$.

Applying the description
of the tangent bundle of an associated bundle
from~\eqref{equation:identifcation_tangent_bundle_associated_bundle}
to the configuration space $\Qspace$,
we obtain for its tangent space
\begin{equation}
	\label{equation:identifcation_tangent_space_configuration_space}
	T \Qspace
	=
	T (\liealg{m} \times (G \times_H \liegroupsub(\liealg{m}) )) 
	\cong
	T \liealg{m} \times T (G \times_H \liegroupsub(\liealg{m}) ) ) 
	\cong
	T \liealg{m} \times (TG \times_{TH}  T\liegroupsub(\liealg{m})).
\end{equation}
Before we proceed,
we state a simple lemma concerning this identification.
We start with considering a
situation which is slightly more general
than~\eqref{equation:identifcation_tangent_space_configuration_space}.

\begin{lemma}
	\label{lemma:equivalence_relation_tangent_bundle_associated_bundle_group_special_case}
	Let $V$ be a finite dimensional $\field{R}$-vector space and
	let $\rho \colon H \to \liegroup{GL}(V)$ be a representation.
	Moreover, let $\liegroup{G}(V) \subseteq \liegroup{GL}(V)$ be a
	closed subgroup such that 
	$\rho_h \in \liegroupsub(V)$ holds for all $h \in H$
	and consider the associated bundle
	$G \times_H \liegroup{G}(V) \to G / H$,
	where $H$ acts on $\liegroupsub(V)$ via 
	$H \times \liegroupsub(V)\ni
	(h, S) \mapsto \rho_h \circ S \in \liegroupsub(V)$.
	Let $(v_g, v_S) \in T G \times T \liegroupsub(V)$
	and $h \in H$.
	Then
	\begin{equation}
		[v_g, v_S] = [T_g r_h v_g, \rho_{h^{-1}} \circ v_S] \in TG \times_{TH} \liegroupsub(V)
	\end{equation}
	holds.
	\begin{proof}
		We denote by
		$\racts \colon G \times H \ni (g, h) \mapsto g \racts h = g h \in G$
		the principal action on $G \to G / H$.
		Its tangent map is given by
		\begin{equation}
			\label{equation:lemma:equivalence_relation_tangent_bundle_associated_bundle_group_special_case_tangent_map_G}
			T_{(g, h)} (\cdot \racts \cdot)(v_g, v_h) 
			= T_h \ell_g v_h + T_g r_h v_g,
		\end{equation}
		due to~\eqref{equation:tangentmap_group_multiplication}.
		Moreover, the tangent map of
		\begin{equation*}
			\phi \colon H \times \liegroupsub(V) \to \liegroupsub(V),
			\quad
			(h, S) \mapsto \rho_{h^{-1}} \circ S
		\end{equation*}
		reads
		\begin{equation}
			\label{equation:lemma:equivalence_relation_tangent_bundle_associated_bundle_group_special_case_tangent_map_GV}
			T_{(h, S)} \phi (v_h, v_S)
			= T_{(h, S)} \phi (0, v_S) + T_{(h, S)} \phi (v_h, 0) ,
		\end{equation}
		where we identify $T(H \times \liegroupsub(V))
		= T H \times T \liegroupsub(V)$.
		By setting $v_h = 0$
		in~\eqref{equation:lemma:equivalence_relation_tangent_bundle_associated_bundle_group_special_case_tangent_map_G}
		and~\eqref{equation:lemma:equivalence_relation_tangent_bundle_associated_bundle_group_special_case_tangent_map_GV},
		respectively,
		we obtain
		$T_{(g, h)} (\cdot \racts \cdot)(v_g, 0) = T_g r_h v_g$
		and
		\begin{equation*}
			T_{(h, S)} \phi (v_h, v_S)
			=
			T_{(h, S)} \phi (0, v_S)
			= 
			T_S \phi(h, \cdot) v_S
			=
			T_S \big(\rho_{h^{-1}} \circ (\cdot) \big) v_S
			= 
			\rho_{h^{-1}} \circ v_S .
		\end{equation*}
		Thus the desired result follows by the definition of the equivalence
		relation in $TG \times_{TH} T \liegroupsub(V)$, i.e.
		$(v_g, v_S) \sim (v_g^{\prime}, v_S^{\prime})
		\in TG \times T \liegroupsub(V)$
		iff there exists an $h \in H$ and $v_h \in T_h H$ such that
		\begin{equation*}
			v_g^{\prime}
			= 
			T_{(g, h)} (\cdot \racts \cdot) (v_g, v_h)
			\quad \text{ and } \quad
			v_S^{\prime} = T_{(h, S)} \phi (v_h, v_S)
		\end{equation*}
		holds.
	\end{proof}
\end{lemma}

\begin{corollary}
	\label{corollary:equivalence_relation_tangent_bundle_associated_bundle_special_case}
	Let $G / H$ be a reductive homogeneous space and let
	$\liegroupsub(\liealg{m}) \subseteq \liegroup{GL}(\liealg{m})$
	be a closed subgroup such that
	$\Ad_h\at{\liealg{m}} \in \liegroupsub(\liealg{m})$ holds
	for all $h \in H$.
	Let $(v_g, v_S) \in T G \times T \liegroupsub(\liealg{m})$
	and $h \in H$. 
	Then
	\begin{equation}
		[v_g, v_S] = [T_g r_h v_g, \Ad_{h^{-1}} \circ v_S] \in TG \times_{TH} \liegroupsub(\liealg{m})
	\end{equation}
	is fulfilled.
	\begin{proof}
		Applying Lemma~\ref{lemma:equivalence_relation_tangent_bundle_associated_bundle_group_special_case}
		to the representation
		$H \ni h \mapsto \Ad_{h}\at{\liealg{m}}
		\in \liegroup{GL}(\liealg{m})$
		yields the desired result because of
		$\Ad_h\at{\liealg{m}} \in \liegroupsub(\liealg{m})$
		for all $h \in H$.
	\end{proof}
\end{corollary}

In order to determine the distribution on $\Qspace$ which characterizes
rollings of $\liealg{m}$ over $G / H$ with respect to $\nablaAlpha$,
we first define a distribution on $\QspaceLift$.
Afterwards, this distribution is used to obtain 
the desired distribution on the
configuration space $\Qspace$.

\begin{lemma}
	\label{lemma:rolling_distribution_lift}
	Let $G / H$ be a reductive homogenoues space.
	Moreover, let
	$\liegroupsub(\liealg{m}) \subseteq \liegroup{GL}(\liealg{m})$ 
	be a closed subgroup and let $\liealgsub(\liealg{m}) \subseteq \liealg{gl}(\liealg{m})$ denote its Lie algebra.
	Assume that $\Ad_h\at{\liealg{m}} \in \liegroupsub(\liealg{m})$ holds
	for all $h \in H$ and
	let $\alpha \colon \liealg{m} \times \liealg{m} \to \liealg{m}$
	be an $\Ad(H)$-invariant bilinear map such that for each
	$X \in \liealg{m}$ the linear map
	\begin{equation}
		\alpha(X, \cdot) \colon \liealg{m} \to \liealg{m},
		\quad
		Y \mapsto \alpha(X, \cdot)(Y) = \alpha(X, Y)
	\end{equation}
	is an element in $\liealgsub(\liealg{m})$,
	i.e. $\alpha(X, \cdot) \in \liealgsub(\liealg{m})$.
	Moreover, let
	$\QspaceLift = \liealg{m} \times G \times \liegroupsub(\liealg{m})$
	as in
	Proposition~\ref{proposition:principal_H_fiber_bundle_over_configuration_space}
	and define
	\begin{equation}
		\label{equation:lemma_rolling_distribution_lift_vb_morphism_defintion}
		\overline{\Psi^{\alpha}} \colon \QspaceLift \times \liealg{m} \to T \QspaceLift,
		\quad
		(\overline{q}, u) 
		= ((v, g, S), u) \mapsto 
		\big(u, (T_e \ell_g \circ S) u, - \alpha(S u, \cdot) \circ S \big) .
	\end{equation}
	Then $\overline{\Psi^{\alpha}}$ is a morphism of vector bundles
	covering $\id_{\QspaceLift} \colon \QspaceLift \to \QspaceLift$
	and
	$\overline{D^{\alpha}} = \image (\overline{\Psi^{\alpha}} ) \subseteq T \QspaceLift$
	is a regular distribution on $\QspaceLift$ given fiber-wise by
	\begin{equation}
		\label{equation:lemma_rolling_distribution_lift_definition}
		\overline{D^{\alpha}}_{(v, g, S)} 
		=
		\big\{(u, (T_e \ell_g \circ S) u, - \alpha(S u, \cdot) \circ S) \mid u \in T_v \liealg{m} \cong \liealg{m} \big\}
		\subseteq 
		T_{(v, g, S)} \QspaceLift
	\end{equation}
	for all $(v, g, S) \in \QspaceLift$.
	Moreover, $\overline{D^{\alpha}}$ is contained in the
	the horizontal bundle defined by the principal connection
	$\overline{\mathcal{P}}$ from 
	Proposition~\ref{proposition:principal_H_fiber_bundle_over_configuration_space},
	i.e. 
	\begin{equation}
		\overline{\mathcal{P}}\at{(v, g, S)}(u, v_g, v_S) = 0
		\quad \text{ for all } \quad 
		(v, g, S) \in \Qspace, 
		\quad
		(u, v_g, v_S) \in \overline{D^{\alpha}}_{(v, g, S)} 
	\end{equation}
	is fulfilled.
	\begin{proof}
		The image of $\overline{\Psi^{\alpha}}$
		defined
		by~\eqref{equation:lemma_rolling_distribution_lift_vb_morphism_defintion}
		is contained in $T \QspaceLift$.
		Indeed, by the assumption on
		$\alpha \colon \liealg{m} \times \liealg{m} \to \liealg{m}$,
		we have
		$\alpha(S u, \cdot) \in \liealgsub(\liealg{m})$ 
		for $S \in \liegroupsub(\liealg{m})$ and
		$u \in \liealg{m}$.
		Hence we obtain
		\begin{equation*}
			\alpha(S u, \cdot) \circ S 
			=
			(T_{\id_{\liealg{m}}} r_{S} ) (\alpha(Su, \cdot)) \in T_S \liegroupsub(\liealg{m}) 
		\end{equation*}
		proving
		$\overline{\Psi^{\alpha}}((v, g, S), u) 
		\in T_v\liealg{m} \times T_g G \times T_S \liegroupsub(\liealg{m})
		\cong T_{(v, g, S)} \QspaceLift$
		for all $((v, g, S), u) \in \QspaceLift \times \liealg{m}$.
		Thus $\overline{\Psi^{\alpha}}$ is clearly a smooth
		vector bundle morphism
		covering the identity.
		Furthermore, the rank of $\overline{\Psi^{\alpha}}$ is obviously
		constant.		
		Hence its image $\overline{D^{\alpha}} = \image(\overline{\Psi^{\alpha}})$
		is a vector subbundle of $T \QspaceLift$
		by~\cite[Thm. 10.34]{lee:2013}.
		The fiber-wise description of $\overline{D^{\alpha}}$
		in~\eqref{equation:lemma_rolling_distribution_lift_definition}
		holds by the definition of $\overline{\Psi^{\alpha}}$
		due to $\overline{D^{\alpha}} = \image(\overline{\Psi^{\alpha}})$.
		
		We now show that $\overline{D^{\alpha}}$ is contained
		in the horizontal bundle.
		Obviously, this is equivalent to 
		$\overline{\mathcal{P}}\at{(v, g, S)}(u, v_g, v_S) = 0$
		for all
		$(v, g, S) \in \QspaceLift$ 
		and $(u, v_g, v_S) \in \overline{D^{\alpha}}_{(v, g, S)}$.
		Using the definition of
		$\overline{\mathcal{P}} \in \Secinfty\big(\End(T \QspaceLift) \big)$
		from
		Proposition~\ref{proposition:principal_H_fiber_bundle_over_configuration_space}
		and writing $(u, v_g, v_S) \in \overline{D^{\alpha}}_{(v, g, S)}$ as
		\begin{equation*}
			(u, v_g, v_S) = (u, (T_e \ell_g \circ S) u, -\alpha(S u, \cdot) \circ S)
		\end{equation*}
		for some $u \in \liealg{m}$,
		we obtain
		\begin{equation*}
			\begin{split}
				\overline{\mathcal{P}}\at{(v, g, S)}(u, v_g, v_S)
				&=
				\overline{\mathcal{P}}\at{(v, g, S)}(u, (T_e \ell_g \circ S) u,
				- \alpha(S u, \cdot) \circ S)  \\
				&=
				(0, \mathcal{P}\at{g}(v_g), 
				- \ad_{\omega\at{g}(v_g)} \circ S) \\
				&=
				(0, 0, 0)
			\end{split}
		\end{equation*}
		due to $v_g = (T_e \ell_g \circ S) u \in \Hor(G)_g$ because of
		$S u \in \liealg{m}$,
		where we used
		$\mathcal{P}\at{g}(v_g) = 0$
		as well as $\omega\at{g}(v_g) = 0$ by the definitions
		of $\mathcal{P} \in \Secinfty\big(\End(T G)\big)$ and
		$\omega \in \Secinfty(T^* G) \tensor \liealg{h}$ in
		Proposition~\ref{proposition:principal_connection_reductive_homogeneous_space}.
	\end{proof}
\end{lemma}
Next we use the distribution $\overline{D^{\alpha}}$ on $\QspaceLift$
to construct the desired
distribution on $\Qspace$.

\begin{lemma}
	\label{lemma:intrinsic_rolling_distribution_first_properties}
	Using the notations and assumptions of
	Lemma~\ref{lemma:rolling_distribution_lift},
	we define $D^{\alpha} \subseteq T \Qspace$ by
	\begin{equation}
		\begin{split}
			D^{\alpha}
			= 
			(T \prQLift ) (\overline{D^{\alpha}})
			\subseteq
			T \Qspace .
		\end{split}
	\end{equation}
	Then the following assertions are fulfilled:
	\begin{enumerate}
		\item
		\label{item:lemma_intrinsic_rolling_distribution_point_wise}
		Let $(v, [g, S]) \in \Qspace$.
		Then $D^{\alpha}$ is fiber-wise given by 
		\begin{equation}
			D^{\alpha}_{(v, [g, S])} 
			= 
			\big\{
			(u, [ (T_e \ell_g \circ S) u, - \alpha(S u, \cdot) \circ S])
			\mid u \in T_v \liealg{m} \cong \liealg{m}
			\big\}
			\subseteq
			T_{(v, [g, S])}  \Qspace 
		\end{equation}
		using the
		identification~\eqref{equation:identifcation_tangent_space_configuration_space}
		implicitly.
		\item
		\label{item:lemma_intrinsic_rolling_curve_horizontal_lift_distribution_fiber_wise_bijective}
		Let $(v, g, S) \in \QspaceLift$. 
		Then the map
		\begin{equation}
			T_{(v, g, S)} \prQLift \at{\overline{D^{\alpha}}_{(v, g, S)}}
			\colon \overline{D^{\alpha}}_{(v, g, S)}
			\to D^{\alpha}_{(v, [g, S])}
		\end{equation}
		is a linear isomorphism.
		\item
		\label{item:lemma_intrinsic_rolling_curve_horizontal_lift_distributions}
		Let $q \colon I \to \Qspace$ be a curve
		and let $\overline{q} \colon I \to \QspaceLift$ denote a horizontal
		lift of $q$ with respect to the principal connection from
		Proposition~\ref{proposition:principal_H_fiber_bundle_over_configuration_space}.
		Then $q$ is horizontal with respect to
		$D^{\alpha}$, i.e. $\dot{q}(t) \in D^{\alpha}_{q(t)}$ iff
		$\overline{q}$ is horizontal with respect to
		$\overline{D^{\alpha}}$,
		i.e. $\dot{\overline{q}}(t) \in \overline{D^{\alpha}}_{\overline{q}(t)}$.
		\item
		\label{item:lemma_intrinsic_rolling_regular_distribution}
		$D^{\alpha}$ is the image of the morphism of vector bundles 
		\begin{equation}
			\Psi^{\alpha} \colon \Qspace \times \liealg{m} \to T \Qspace,
			\quad
			((v, [g, S]), u)\mapsto (u, [ (T_e \ell_g \circ S) u, - \alpha(Su, \cdot) \circ S]) 
		\end{equation}
		over $\id_{\Qspace} \colon \Qspace \to \Qspace$ of constant rank.
		In particular, $D^{\alpha}$ is a regular distribution on $\Qspace$.
	\end{enumerate}
	\begin{proof}
		We start with determining $D^{\alpha}$ point-wise.
		Let $(v, g, S) \in \QspaceLift
		=
		\liealg{m} \times G \times \liegroupsub(\liealg{m})$
		and
		$(u, v_g, v_S) \in T \QspaceLift
		\cong
		T \liealg{m} \times T G \times T \liegroupsub(\liealg{m})$.
		Then 
		\begin{equation}
			\label{equaiton:theorem_intrinsic_rolling_evaluation_prQLift_and_ttangent_map}
			\prQLift(v, g, S) = (v, [g, S])
			\quad \text{ and } \quad
			T \prQLift(u, v_g, v_S)
			=
			(u, [v_g, v_S])
		\end{equation}
		holds by the
		identification~\eqref{equation:identifcation_tangent_space_configuration_space}.
		Evaluating~\eqref{equaiton:theorem_intrinsic_rolling_evaluation_prQLift_and_ttangent_map}
		at $(u, v_g, v_S) \in \overline{D^{\alpha}}_{(v, g, S)}$, i.e.
		\begin{equation*}
			(u, v_g, v_S) = (u, (T_e \ell_g \circ S) u, - \alpha(S u, \cdot) \circ S)
		\end{equation*}
		for some $u \in \liealg{m}$, yields
		Claim~\ref{item:lemma_intrinsic_rolling_distribution_point_wise}
		because of
		$D^{\alpha}_{\prQLift(v, g, S)} 
		= (T \prQLift)\big(\overline{D^{\alpha}}_{(v, g, S)} \big)$.
		
		Next we show
		Claim~\ref{item:lemma_intrinsic_rolling_curve_horizontal_lift_distribution_fiber_wise_bijective},
		i.e. that the restriction
		\begin{equation}
			\label{equation:lemma:intrinsic_rolling_distribution_first_properties_tangent_map_prQLift_restriction_fiberwise_surjective}
			T_{(v, g, S)} \prQLift \at{\overline{D^{\alpha}}_{(v, g, S)}}
			\colon
			\overline{D^{\alpha}}_{(v, g, S)}
			\to
			D^{\alpha}_{(v, [g, S])}
		\end{equation}
		is bijective.
		Clearly, the linear map
		in~\eqref{equation:lemma:intrinsic_rolling_distribution_first_properties_tangent_map_prQLift_restriction_fiberwise_surjective}
		is injective since
		$\overline{D^{\alpha}}_{(v, g, S)} \subseteq \Hor(\overline{Q})_{(v, g, S)}$
		holds according to
		Lemma~\ref{lemma:rolling_distribution_lift}.
		We now show that~\eqref{equation:lemma:intrinsic_rolling_distribution_first_properties_tangent_map_prQLift_restriction_fiberwise_surjective}
		is surjective.
		Let $h \in H$. Moreover, let $(v, g, S) \in \QspaceLift$
		and $(v, g h , \Ad_{h^{-1}} \circ S) \in \QspaceLift$
		be two representatives of
		\begin{equation*}
			\prQLift\big(v, g, S \big) =
			(v, [g, S]) = (v, [g h, \Ad_{h^{-1}} \circ S] )
			= \prQLift\big(v, g h, \Ad_{h^{-1}} \circ S \big) \in \Qspace
		\end{equation*}
		and let
		$(u, v_g, v_S) \in \overline{D^{\alpha}}_{(v, g, S)}$.
		We show that there exists a
		$(u, v_g^{\prime}, v_S^{\prime}) \in \overline{D^{\alpha}}_{(v, gh, \Ad_{h^{-1}} \circ S)}$
		such that
		\begin{equation*}
			T_{(v, g, S)} \prQLift (u, v_g, v_S) 
			=
			T_{(v, g h,  \Ad_{h^{-1}} \circ S)}
			\prQLift (u, v_g^{\prime}, v_S^{\prime})
		\end{equation*}
		holds.
		To this end, we define
		\begin{equation}
			\label{equation:lemma_intrinsic_rolling_distribution_first_properties_v_g_v_S_prime_def}
			v_{g}^{\prime} = T_g (\cdot \racts h) v_g = T_g r_h v_g
			\quad \text{ and } \quad
			v_{S}^{\prime} = T_S \big( \Ad_{h^{-1}}(\cdot) \big) \circ  v_S
			= \Ad_{h^{-1}} \circ v_S .
		\end{equation}
		By using $(u, v_g, v_S) \in \overline{D^{\alpha}}_{(v, g, S)}$, i.e.
		\begin{equation}
			\label{equation:lemma_intrinsic_rolling_distribution_first_properties_v_g_v_S_def}
			v_g = (T_e \ell_g \circ S) u
			\quad \text{ and } \quad
			v_S = - \alpha(S u, \cdot) \circ S 
		\end{equation}
		for some $u \in \liealg{m}$,
		we show that
		$(u, v_g^{\prime}, v_S^{\prime}) \in \overline{D^{\alpha}}_{(v, gh , \Ad_{h^{-1}} \circ S)}$ holds.
		To this end, we calculate 
		\begin{equation}
			\label{equation:lemma:intrinsic_rolling_distribution_first_properties_v_g_prime_calc}
			\begin{split}
				v_g^{\prime}
				&=
				T_g r_h v_g \\
				&=
				T_g r_h \big( (T_e \ell_g \circ S ) u \big)\\
				&=
				T_e (r_h \circ \ell_g) S u \\
				&=
				T_e (\ell_g \circ r_h) S u \\
				&=
				T_e (\ell_{g h} \circ \ell_{h^{-1}} \circ r_h ) S u \\
				&=
				T_e \ell_{g h} \circ T_e (\ell_{h^{-1}} \circ r_h) S u \\
				&=
				(T_e \ell_{g h} \circ \Ad_{h^{-1}} \circ S) u .
			\end{split}
		\end{equation}
		Moreover, using the definition of $v_S^{\prime}$ 
		in~\eqref{equation:lemma_intrinsic_rolling_distribution_first_properties_v_g_v_S_prime_def}
		and $v_S$
		in~\eqref{equation:lemma_intrinsic_rolling_distribution_first_properties_v_g_v_S_def},
		we have by the $\Ad(H)$-invariance
		of $\alpha \colon \liealg{m} \times \liealg{m} \to \liealg{m}$
		\begin{equation}
			\label{equation:lemma:intrinsic_rolling_distribution_first_properties_v_S_prime_calc}
			\begin{split}
				v_S^{\prime}
				&=
				\Ad_{h^{-1}} \circ v_S \\
				&=
				\Ad_{h^{-1}} \circ \big(- \alpha(Su, \cdot) \big) \circ S \\
				&=
				- \alpha(\Ad_{h^{-1}} (S u), \cdot) \circ \Ad_{h^{-1}} \circ S .
			\end{split}
		\end{equation}
		By~\eqref{equation:lemma:intrinsic_rolling_distribution_first_properties_v_g_prime_calc}
		and~\eqref{equation:lemma:intrinsic_rolling_distribution_first_properties_v_S_prime_calc}, we obtain
		\begin{equation}
			\label{equation:lemma_intrinsic_rolling_distribution_first_properties_v_g_v_S_prime_result}
			(u, v_g^{\prime}, v_S^{\prime})
			=
			(u, 	(T_e \ell_{g h} \circ \Ad_{h^{-1}} \circ S) u,
			- \alpha(\Ad_{h^{-1}} (S u), \cdot) \circ \Ad_{h^{-1}} \circ S)  
		\end{equation}
		showing $(u, v_g^{\prime}, v_S^{\prime})
		\in
		\overline{D^{\alpha}}_{(v, g h, \Ad_{h^{-1}} \circ S)}$
		as desired.
		Equation~\eqref{equation:lemma_intrinsic_rolling_distribution_first_properties_v_g_v_S_prime_result}
		implies
		\begin{equation*}
			\begin{split}
				T_{(v, g, S)} \prQLift (u, v_g, v_S) 
				&=
				(u, [v_g, v_S]) \\
				&=
				(u, [T_g r_h v_g, \Ad_{h^{-1}} \circ v_S] ) \\
				&=
				T_{(v, g h,  \Ad_{h^{-1}} \circ S)}
				\prQLift (u, v_g^{\prime}, v_S^{\prime}) 
			\end{split}
		\end{equation*}
		due to
		Corollary~\ref{corollary:equivalence_relation_tangent_bundle_associated_bundle_special_case}.
		Thus
		the linear map~\eqref{equation:lemma:intrinsic_rolling_distribution_first_properties_tangent_map_prQLift_restriction_fiberwise_surjective}
		is surjective.
		Hence Claim~\ref{item:lemma_intrinsic_rolling_curve_horizontal_lift_distribution_fiber_wise_bijective}
		is proven.
		
		Next let $q \colon I \to \Qspace$ be a curve and let
		$\overline{q} \colon I \to \QspaceLift$ be a horizontal lift
		with respect to
		the principal connection $\overline{\mathcal{P}}$ from
		Proposition~\ref{proposition:principal_H_fiber_bundle_over_configuration_space}.
		In particular $\prQLift(\overline{q}(t)) = q(t)$ holds.
		Assume that $\dot{\overline{q}}(t) \in \overline{D^{\alpha}}_{\overline{q}(t)}$
		holds. This assumption yields
		\begin{equation*}
			\dot{q}(t) 
			= 
			\tfrac{\D}{\D t} (\prQLift \circ \overline{q})(t)
			=
			(T_{\overline{q}(t)} \prQLift ) \dot{\overline{q}}(t)
			\in
			(T_{\overline{q}(t)} \prQLift )
			\big( \overline{D^{\alpha}}_{\overline{q}(t)} \big)
			=
			D^{\alpha}_{q(t)}
		\end{equation*}
		by the definition of $D^{\alpha}$ since $\overline{q}(t)$ is
		horizontal with respect to $\overline{D^{\alpha}}$.
		Conversely, assume that
		$\dot{q}(t) \in D^{\alpha}_{q(t)}$ holds.
		Then $\dot{\overline{q}}(t) \in T_{q(t)} \QspaceLift$ is the unique
		horizontal tangent vector which fulfills
		$(T_{\overline{q}(t)} \prQLift) \dot{\overline{q}}(t) = \dot{q}(t)$
		or equivalently
		\begin{equation*}
			\dot{\overline{q}}(t) = \big(T_{\overline{q}(t)} \prQLift \at{\Hor(\QspaceLift)_{\overline{q}(t)}}\big)^{-1} \dot{q}(t) .
		\end{equation*}
		Since~\eqref{equation:lemma:intrinsic_rolling_distribution_first_properties_tangent_map_prQLift_restriction_fiberwise_surjective}
		is a linear isomorphism, we obtain
		$\big(T_{\overline{q}(t)} \prQLift
		\at{\overline{D^{\alpha}}_{\overline{q}(t)}} \big)^{-1}(D^{\alpha}_{q(t)})
		=
		\overline{D^{\alpha}}_{\overline{q}(t)}$.
		This
		yields
		Claim~\ref{item:lemma_intrinsic_rolling_curve_horizontal_lift_distributions}
		because of $\dot{q}(t) \in D^{\alpha}_{q(t)}$
		and $\overline{D^{\alpha}} \subseteq \Hor(\QspaceLift)$.
		
		It remains to proof Claim~\ref{item:lemma_intrinsic_rolling_regular_distribution}.
		To this end, using the
		identification
		~\eqref{equation:identifcation_tangent_space_configuration_space},
		we consider the diagram
		\begin{equation}
			\label{equation:lemma_intrinsic_rolling_distribution_first_properties_commutative_diagram}
			\begin{tikzpicture}[node distance= 5.0cm, auto]
				\node (QLiftxm) {$\QspaceLift \times \liealg{m}$};
				\node (TQLift) [right=3.0cm of QLiftxm] {$T \QspaceLift$};
				\node (Qxm) [below= 1.5cm of QLiftxm] {$\Qspace \times \liealg{m}$};
				\node (TQ) [right=3.0cm of Qxm] {$T \Qspace$};
				\draw[->] (QLiftxm) to node [left] {$\prQLift \times \id_{\liealg{m}}$} (Qxm);
				\draw[->] (Qxm) to node  [above]{$\Psi^{\alpha}$} 
				(TQ);
				\draw[->] (QLiftxm) to node [above] {$\overline{\Psi^{\alpha}}$} 
				(TQLift);
				\draw[->] (TQLift) to node [right] {$T \prQLift$} 
				(TQ);
			\end{tikzpicture}  
		\end{equation}
		which clearly commutes.
		Thus $\Psi^{\alpha}$ is smooth since
		$T \prQLift \circ \overline{\Psi^{\alpha}}$ is smooth
		and $\prQLift \times \id_{\liealg{m}}$ is a surjective submersion.
		In addition, $\Psi^{\alpha}$ is fiber-wise linear, i.e.
		$\Psi^{\alpha}$
		is a vector bundle morphism covering
		$\id_{\Qspace} \colon \Qspace \to \Qspace$.
		Moreover,
		Claim~\ref{item:lemma_intrinsic_rolling_curve_horizontal_lift_distribution_fiber_wise_bijective}
		implies that the rank of $\Psi$ is constant.
		Hence the image of $\Psi^{\alpha}$ is a subbundle of $T \Qspace$
		according to~\cite[Thm. 10.34]{lee:2013}. 
		In addition, we obtain $D^{\alpha} = \image(\Psi^{\alpha})$ due to 
		\begin{equation}
			D^{\alpha} 
			=
			(T \prQLift)(\overline{D^{\alpha}})
			=
			(T \prQLift \circ \overline{\Psi^{\alpha}})(\QspaceLift \times \liealg{m})
			=
			(\Psi^{\alpha} \circ (\prQLift \times \id_{\liealg{m}}))(\QspaceLift \times \liealg{m})
			=
			\Psi^{\alpha}(\Qspace \times \liealg{m})
		\end{equation}
		since~\eqref{equation:lemma_intrinsic_rolling_distribution_first_properties_commutative_diagram}
		commutes.
		This yields the desired result.
	\end{proof}
\end{lemma}

\begin{theorem}
	\label{theorem:intrinsic_rolling}
	Let $G / H$ be a reductive homogeneous space with reductive
	decomposition $\liealg{g} = \liealg{h} \oplus \liealg{m}$
	and let $\liegroupsub(\liealg{m}) \subseteq \liegroup{GL}(\liealg{m})$
	be a closed subgroup such that
	$\Ad_h\at{\liealg{m}} \in \liegroupsub(\liealg{m})$ holds
	for all $h \in H$.
	Moreover, let
	$\alpha \colon \liealg{m} \times \liealg{m} \to \liealg{m}$
	be an $\Ad(H)$-invariant bilinear map
	defining the invariant covariant derivative $\nablaAlpha$
	such that for each $X \in \liealg{m}$ the linear map
	\begin{equation}
		\alpha(X, \cdot) \colon \liealg{m} \to \liealg{m},
		\quad 
		Y \mapsto \alpha(X, Y)
	\end{equation}
	belongs to $\liealgsub(\liealg{m})$,
	i.e. to the Lie algebra of
	$\liegroupsub(\liealg{m})$.
	Let $\overline{D^{\alpha}}$ denote the distribution 
	on $\QspaceLift = \liealg{m} \times G \times \liegroupsub(\liealg{m})$
	from Lemma~\ref{lemma:rolling_distribution_lift}
	associated to $\alpha$
	and let $D^{\alpha} = T \prQLift (\overline{D^{\alpha}})$
	be the distribution defined in
	Lemma~\ref{lemma:intrinsic_rolling_distribution_first_properties}.
	Then the following assertions are fulfilled:
	\begin{enumerate}
		\item
		\label{item:theorem_intrinsic_rolling_horizontal_lift}
		Let $q \colon I \to \Qspace$ and let
		\begin{equation}
			(v, \gamma) \colon I \to \liealg{m} \times G /H, 
			\quad
			t \mapsto (\prQ \circ q)(t) = (v(t), \gamma(t)) .
		\end{equation}
		Let
		$\overline{q} \colon I \ni t \mapsto (v(t),g(t), S(t)) \in \QspaceLift$
		be a horizontal
		lift of $q$
		with respect to the principal connection $\overline{\mathcal{P}}$
		from Proposition~\ref{proposition:principal_H_fiber_bundle_over_configuration_space}.
		Then $q$ is horizontal with respect to $D^{\alpha}$ iff
		the ODE
		\begin{equation}
			\label{equation:theorem_intrinsic_rolling_horizontal_lift_initial_value_problem}
			\begin{split}
				\dot{S}(t) 
				&= - \alpha\big(S(t) \dot{v}(t), \cdot \big) \circ S(t) , \\
				\dot{g}(t) &= \big(T_e \ell_{g(t)} \circ S(t) \big) \dot{v}(t) 
			\end{split}
		\end{equation}
		is fulfilled.
		Moreover, the development curve
		is given by 
		$\gamma = \pr \circ g \colon I \to G / H$.
		\item
		\label{item:theorem_intrinsic_rolling_distribution_rolling}
		Let $q \colon I \to \Qspace$ be a curve and let
		$(v, \gamma) = (\prQ \circ q) \colon I \to \liealg{m} \times G / H$.
		Then $q$ is horizontal with respect to $D^{\alpha}$,
		i.e. $\dot{q}(t) \in D^{\alpha}_{q(t)}$,
		iff $q$ defines a ($\liegroupsub(\liealg{m})$-reduced)
		intrinsic rolling of $\liealg{m}$ over
		$G / H$ with respect to $\nablaAlpha$
		with rolling curve $v$ and development curve $\gamma$.
	\end{enumerate}
	\begin{proof}
		We first show
		Claim~\ref{item:theorem_intrinsic_rolling_horizontal_lift}.
		Let $q \colon I \to \Qspace$ be some curve and let
		\begin{equation}
			\label{equation:theorem_intrinsic_rolling_horizontal_lift_ode_curve_as_triple}
			\overline{q} \colon I \to \QspaceLift
			= \liealg{m} \times G \times \liegroupsub(\liealg{m}),
			\quad
			t \mapsto (v(t), g(t), S(t))
		\end{equation}
		be a horizontal lift of $q$ with respect to
		the principal connection $\overline{\mathcal{P}}$
		from
		Proposition~\ref{proposition:principal_H_fiber_bundle_over_configuration_space}.
		Clearly, the development curve $\gamma \colon I \to G / H$ defined by
		$q \colon I \to \Qspace$ is given by $\gamma = \pr \circ g$.
		Moreover, by
		Lemma~\ref{lemma:intrinsic_rolling_distribution_first_properties},
		Claim~\ref{item:lemma_intrinsic_rolling_curve_horizontal_lift_distributions},
		$q$ is horizontal with respect to $D^{\alpha}$ iff $\overline{q}$
		is horizontal with respect to $\overline{D^{\alpha}}$.
		Hence it is sufficient to show that $\overline{q}$
		from~\eqref{equation:theorem_intrinsic_rolling_horizontal_lift_ode_curve_as_triple}
		fulfills the
		ODE~\eqref{equation:theorem_intrinsic_rolling_horizontal_lift_initial_value_problem} iff $\overline{q}$ is horizontal with
		respect to $\overline{D^{\alpha}}$.
		First we assume
		$\dot{\overline{q}}(t) \in \overline{D^{\alpha}}_{\overline{q}(t)}$
		for all $t \in I$.
		Writing $\dot{\overline{q}}(t) = (\dot{v}(t), \dot{g}(t), \dot{S}(t))$
		and using the definition of $\overline{D^{\alpha}}$, one obtains
		\begin{equation}
			\label{equation:theorem_intrinsic_rolling_horizontal_lift_tangential_to_distribution}
			\dot{\overline{q}}(t)
			\in 
			\overline{D^{\alpha}}_{\overline{q}(t)}
			=
			\big\{ (u, \big(T_e \ell_{g(t)} \circ S(t) \big) u, - \alpha(S(t) u, \cdot) \circ S(t) ) \mid u \in T_{v(t)} \liealg{m} \cong \liealg{m} \big\} .
		\end{equation} 
		Thus
		$\dot{g}(t)$ and $\dot{S}(t)$ are uniquely determined by
		\begin{equation*}
			\begin{split}
				\dot{g}(t) = \big( T_e \ell_{g(t)} \circ S(t) \big) \dot{v}(t)
				\quad \text{ and } \quad
				\dot{S}(t) = - \alpha\big( S(t) \dot{v}(t), \cdot \big) \circ S(t) 
			\end{split}
		\end{equation*}
		due
		to~\eqref{equation:theorem_intrinsic_rolling_horizontal_lift_tangential_to_distribution}.
		Hence the curve $\overline{q} \colon I \to \QspaceLift$
		which is horizontal with respect to $\overline{D^{\alpha}}$
		fulfills the
		ODE~\eqref{equation:theorem_intrinsic_rolling_horizontal_lift_initial_value_problem}.
		Conversely, assume that $\overline{q} \colon I \to \QspaceLift$
		given by $\overline{q}(t) = (v(t), g(t), S(t))$
		fulfills~\eqref{equation:theorem_intrinsic_rolling_horizontal_lift_initial_value_problem}.
		Then $\overline{q}(t)$ is clearly horizontal with respect to
		$\overline{D^{\alpha}}$
		by the definition of $\overline{D^{\alpha}}$.
		Thus  Claim~\ref{item:theorem_intrinsic_rolling_horizontal_lift}
		is proven.
		
		Next we show
		Claim~\ref{item:theorem_intrinsic_rolling_distribution_rolling}.
		To this end, let $q \colon I \to \Qspace$ be horizontal with respect to $D^{\alpha}$.
		Then a horizontal lift
		$\overline{q} \colon I \ni t \mapsto \overline{q}(t) = (v(t), g(t), S(t)) \in \QspaceLift$
		of $q$
		fulfills the
		ODE~\eqref{equation:theorem_intrinsic_rolling_horizontal_lift_initial_value_problem}
		by Claim~\ref{item:theorem_intrinsic_rolling_horizontal_lift}.
		Moreover, $q \colon I \to Q$ can be represented by 
		\begin{equation*}
			q(t) = (\prQLift \circ \overline{q})(t) = (v(t), [g(t), S(t)]),
			\quad
			t \in I .
		\end{equation*}
		Hence the linear isomorphism associated with $q(t)$ is given by
		\begin{equation}
			\label{equation:theorem_intrinsic_rolling_definition_isometry}
			T_{v(t)} \liealg{m} \cong \liealg{m}
			\ni Z \mapsto
			q(t)Z
			=
			\big( T_{g(t)} \pr \circ T_e \ell_{g(t)}  \circ S(t) \big) Z \in T_{\pr(g(t))} (G / H)
		\end{equation}
		according to
		Lemma~\ref{lemma:configuration:space_intrinsic_rolling},
		Claim~\ref{item:lemma_configuration:space_intrinsice_rolling_associated_isometry}.
		Using\eqref{equation:theorem_intrinsic_rolling_horizontal_lift_initial_value_problem}
		and~\eqref{equation:theorem_intrinsic_rolling_definition_isometry}
		we obtain
		\begin{equation*}
			\dot{\gamma}(t) 
			=
			\tfrac{\D}{\D t} (\pr \circ g)(t) 
			=
			T_{g(t)} \pr \dot{g}(t) 
			=
			\big(T_{g(t)} \pr \big) \Big( (T_e \ell_{g(t)} \circ S(t) ) \dot{v}(t) \Big)
			=
			q(t) \dot{v}(t)
		\end{equation*}
		showing the no-slip condition.
		Next we prove the no-twist condition. Let
		\begin{equation*}
			Z \colon I \to T \liealg{m}, 
			\quad
			t \mapsto (v(t), Z_2(t)) \in T \liealg{m} 
		\end{equation*} 
		be a vector field along $v \colon I \to \liealg{m}$ which we identify with the map
		$I \ni t \mapsto Z_2(t) \in \liealg{m}$
		defined by its second component.
		Then $Z$ is parallel along $v$ iff $\dot{Z}_2(t) = 0$ holds,
		i.e. $Z_2(t) = Z_0$ for all $t \in I$ and some $Z_0 \in \liealg{m}$.
		We need to show that the vector field
		\begin{equation*}
			\widehat{Z} \colon I \to T(G / H), 
			\quad
			t \mapsto q(t) Z(t) = q(t) Z_0
			= \big( T_{g(t)} \pr \circ T_e \ell_{g(t)} \circ S(t) \big) Z_0
		\end{equation*}
		is parallel along $\gamma \colon I \to G / H$
		with respect to $\nablaAlpha$.
		By
		Proposition~\ref{proposition:principal_H_fiber_bundle_over_configuration_space},
		Claim~\ref{item:proposition:principal_H_fiber_bundle_over_configuration_space_definition_horizontal_lift_curve},
		the curve $g \colon I \to G$ is a horizontal
		lift of the curve 
		$\gamma \colon I \ni t \mapsto \pr(g(t)) \in G / H$.
		In addition, we have  
		\begin{equation*}
			S(t) \dot{v}(t) = \big( T_e \ell_{g(t)}\big)^{-1} \dot{g}(t)
		\end{equation*}
		by~\eqref{equation:theorem_intrinsic_rolling_horizontal_lift_initial_value_problem}.
		Moreover, the horizontal lift of $\widehat{Z}$
		along $g \colon I \to G$ is given by
		\begin{equation*}
			\overline{Z} \colon I \to \Hor(G),
			\quad
			t \mapsto
			\bar{Z}(t) 
			=
			\big(T_{g(t)} \pr\at{\Hor(G)}\big)^{-1} \widehat{Z}(t)	
			=
			\big( T_e \ell_{g(t)} \circ  S(t) \big) Z_0
		\end{equation*}
		and the curve $z \colon I \ni
		t \mapsto 
		z(t) 
		=
		\big(T_e \ell_{g(t)} \big)^{-1} \overline{Z}(t)
		\in 
		\liealg{m}$
		fulfills
		\begin{equation*}
			z(t)
			=
			\big(T_e \ell_{g(t)} \big)^{-1} \big( T_e \ell_{g(t)} \circ S(t) \big) Z_0
			=
			S(t) Z_0 .
		\end{equation*}
		Thus we obtain by
		exploiting~ \eqref{equation:theorem_intrinsic_rolling_horizontal_lift_initial_value_problem}
		\begin{equation*}
			\dot{z}(t) 
			=
			\tfrac{\D}{\D t} (S(t) Z_0) 
			=
			\dot{S}(t) Z_0 
			= 
			- \big( \alpha\big(S(t)\dot{v}(t), \cdot \big) \circ S(t) \big) Z_0
			= 
			- \alpha\big(S(t) \dot{v}(t), z(t) \big) .
		\end{equation*}
		Hence $\widehat{Z}$ is parallel along $\gamma$ by
		Proposition~\ref{proposition:covariant_derivative_along_curve}.
		
		Conversely, assume that $\widehat{Z}(t) = q(t) Z(t)$ is parallel
		along $\gamma \colon I \to G / H$,
		where $Z \colon I \ni t \mapsto (v(t), Z_2(t)) \in T \liealg{m}$
		is some vector field along
		along $v \colon I \to \liealg{m}$
		which we identify with
		the map $I \ni t \mapsto Z_2(t) \in \liealg{m}$.
		Let $\{A_1, \ldots, A_N\} \subseteq \liealg{m}$ be some basis of $\liealg{m}$.
		We define a parallel frame along $\gamma \colon I \to G / H$ by 
		$A_i(t) = q(t) A_i$ for $i \in \{1, \ldots, N\}$ and $t \in I$.
		Then $\widehat{Z} \colon I \to T (G / H)$ is parallel along $\gamma$ iff 
		its coefficient functions $z^i \colon I \to \field{R}$ defined by
		$\widehat{Z}(t) = z^i(t) A_i(t)$ are constant, i.e. 
		$z^i(t) = z^i_0$ for all $t \in I$ with some $z^i_0 \in \field{R}$, 
		see e.g.~\cite[Chap. 4, p. 109]{lee:2018}.
		By the linearity of 
		$q(t) \colon T_{v(t)} \liealg{m} \cong \liealg{m} \to T_{\gamma(t)} (G / H)$,
		we obtain
		\begin{equation*}
			\widehat{Z}(t) 
			=
			z^i_0 A_i(t) 
			=
			z^i_0 \big(q(t) A_i\big)
			= 
			q(t) \big( z^i_0 A_i \big)
			=
			q(t) Z_0
			=
			q(t) Z_2(t)
		\end{equation*}
		showing $Z_2(t) = z^i_0 A_i = Z_0$ for $t \in I$,
		where $Z_0 = z^i_0 A_i \in \liealg{m}$ is constant.
		Hence $Z \colon I \ni t \mapsto (v(t), Z_2(t)) = (v(t), Z_0) \in \liealg{m} \times \liealg{m} \cong T \liealg{m}$
		is a parallel vector field along the curve
		$v \colon I \to \liealg{m}$.
		Thus the curve $q \colon I \to \Qspace$
		which is horizontal with
		respect to $D^{\alpha}$
		is a rolling.

		It remains to proof the converse.
		Let $q \colon I \to \Qspace$ be a curve defining a rolling.
		We show that $q$ is horizontal with respect to $D^{\alpha}$.
		
		Let $\overline{q} \colon I \to \Qspace$ be a horizontal lift of $q$
		with respect to the principal connection from
		Proposition~\ref{proposition:principal_H_fiber_bundle_over_configuration_space}.
		By
		Lemma~\ref{lemma:intrinsic_rolling_distribution_first_properties}, Claim~\ref{item:lemma_intrinsic_rolling_curve_horizontal_lift_distributions}
		$q$ is horizontal with respect to $D^{\alpha}$ iff
		$\overline{q}$ is horizontal with respect to
		$\overline{D^{\alpha}}$.
		Writing $\overline{q}(t) = (v(t), g(t), S(t))$, the
		linear isomorphism
		$T_{v(t)} \liealg{m} \cong \liealg{m} \to T_{\pr(g(t))} (G / H)$
		defined by
		$q(t) = \prQLift(\overline{q}(t))$
		is given by 
		\begin{equation*}
			q(t)Z = \big( T_{g(t)} \pr \circ T_e \ell_{g(t)}  \circ S(t) \big) Z,
			\quad Z \in T_{v(t)} \liealg{m} \cong \liealg{m}
		\end{equation*}
		according to Lemma~\ref{lemma:configuration:space_intrinsic_rolling}, Claim~\ref{item:lemma_configuration:space_intrinsice_rolling_associated_isometry}.
		Hence the no slip condition yields
		\begin{equation}
			\label{equation:theorem_intrinsic_rolling_curve_config_space_horizontal_no_slip}
			\dot{\gamma}(t) 
			=
			\big( T_{g(t)} \pr \circ T_e \ell_{g(t)}  \circ S(t) \big)
			\dot{v}(t) .
		\end{equation}
		By 	Proposition~\ref{proposition:principal_H_fiber_bundle_over_configuration_space},
		Claim~\ref{item:proposition_principal_H_fiber_bundle_over_configuration_space_g_component}
		the curve  $g \colon I \to G$ is horizontal with
		respect to $\Hor(G)$ from
		Proposition~\ref{proposition:principal_connection_reductive_homogeneous_space}.
		In addition, $\gamma = \pr \circ g$
		holds, i.e.
		$g \colon I \to G$ is a horizontal lift of $\gamma \colon I \to G / H$.
		Thus $g \colon I \to G$ fulfills the ODE
		$\dot{g}(t) = \big( T_e \ell_{g(t)} \circ S(t)\big) \dot{v}(t)$
		by~\eqref{equation:theorem_intrinsic_rolling_curve_config_space_horizontal_no_slip}.
		Moreover, since $g \colon I \to G$ is
		a horizontal lift of $\gamma \colon I \to G / H$,
		the no twist condition yields
		\begin{equation}
			\label{equation:theorem_intrinsic_rolling_curve_config_space_horizontal_no_twist}
			S(t)Z_0 = - \alpha\big(S(t) \dot{v}(t), S(t) Z_0 \big)
		\end{equation}
		for all $Z_0 \in \liealg{m}$
		by Proposition~\ref{proposition:covariant_derivative_along_curve}.
		Clearly,~\eqref{equation:theorem_intrinsic_rolling_curve_config_space_horizontal_no_twist}
		is equivalent to the ODE
		\begin{equation*}
			\dot{S}(t) = - \alpha\big(S(t) \dot{v}(t), \cdot \big) \circ S(t) 
		\end{equation*}
		for $S \colon I \to \liegroupsub(\liealg{m})$.
		Thus $\overline{q} \colon I \to \QspaceLift$ is horizontal with respect to $\overline{D^{\alpha}}$ by
		Claim~\ref{item:theorem_intrinsic_rolling_horizontal_lift}.
		Therefore $q \colon I\to \Qspace$ is horizontal with respect to
		$D^{\alpha}$ due to
		Lemma~\ref{lemma:intrinsic_rolling_distribution_first_properties},
		Claim~\ref{item:lemma_intrinsic_rolling_curve_horizontal_lift_distributions}.
	\end{proof}
\end{theorem}
In particular, Theorem~\ref{theorem:intrinsic_rolling}
applies to (pseudo-)Riemannian reductive homogeneous spaces.
We comment on this particular situation in the next remark.

\begin{remark}
	\label{remark:intrinsic_rolling_pseudo_riemannian_metric_Qspace_reduced}
	Let $G / H$ be a reductive homogeneous space
	equipped with an invariant pseudo-Riemannian metric
	and let $\langle \cdot, \cdot \rangle \colon \liealg{m} \times \liealg{m}
	\to \field{R}$ denote the corresponding $\Ad(H)$-invariant scalar product.
	Moreover, let $\nablaAlpha$ be a \emph{metric} invariant covariant
	derivative on $G / H$.
	Then
	Proposition~\ref{proposition:invariant_covariant_derivative_compatible_with_structure_skew_adjoint}
	yields $\alpha(X, \cdot) \in \liealg{so}(\liealg{m})$
	for all $X \in \liealg{m}$.
	Thus Theorem~\ref{theorem:intrinsic_rolling}
	can be applied to $\liegroupsub(\liealg{m}) = \liegroup{O}(\liealg{m})$, 
	i.e. the configuration space can be reduced to
	$\Qspace = \liealg{m} \times (G \times_H \liegroup{O}(\liealg{m}))$
	since $\Ad_h\at{\liealg{m}} \in \liegroup{O}(\liealg{m})$
	holds for all $h \in H$.
\end{remark} Remark~\ref{remark:intrinsic_rolling_pseudo_riemannian_metric_Qspace_reduced}
can be specialized further to naturally reductive homogeneous space
equipped with the Levi-Civita covariant derivative.

\begin{remark}
	Let $G / H$ be a naturally reductive homogeneous space.	
	Then Theorem~\ref{theorem:intrinsic_rolling}
	can be applied to $G / H$ equipped with $\nablaLC$, where
	the configuration space can be reduced to
	$\Qspace = \liealg{m} \times (G \times_H \liegroup{O}(\liealg{m}))$
	and $\alpha \colon \liealg{m} \times \liealg{m} \to \liealg{m}$
	is given by
	$\alpha(X, Y) = \tfrac{1}{2} \pr_{\liealg{m}} \circ \ad_X (Y)$
	for $X, Y \in \liealg{m}$
	since  $\nablaLC = \nablaCan$ holds by
	Remark~\ref{remark:horizontal_lift_covariant_derivative_naturally_redutive_and_symmetric}.
\end{remark}

\subsection{Kinematic Equations and  Control Theoretic Perspective}

Throughout this section we denote by $G / H$ a reductive homogeneous space
and we assume that $\liegroupsub(\liealg{m}) \subseteq \liegroup{GL}(\liealg{m})$ is a closed subgroup such that $\Ad_{h} \at{\liealg{m}} \in \liegroupsub(\liealg{m})$ holds.
Moreover, we assume that the $\Ad(H)$-invariant bilinear
map $\alpha \colon \liealg{m} \times \liealg{m} \to \liealg{m}$
fulfills $\alpha(X, \cdot) \in \liealgsub(\liealg{m})$
for all $X \in \liealg{m}$.
If not indicated otherwise, we consider the ``reduced''
configuration space
$\Qspace = \liealg{m} \times (G \times_H \liegroupsub(\liealg{m}))$.

We start with
relating rollings of $\liealg{m}$ over $G / H$
to a control system.

\begin{remark}
	\label{remark:intrinsic_rolling_control_system}
	Let $G / H$ be a reductive homogeneous space equipped with an
	invariant covariant derivative $\nablaAlpha$
	defined by the $\Ad(H)$-invariant bilinear map
	$\alpha \colon \liealg{m} \times \liealg{m} \to \liealg{m}$.
	Recall the definition of the morphism of vector bundles $\Psi^{\alpha}$
	in Lemma~\ref{lemma:intrinsic_rolling_distribution_first_properties},
	Claim~\ref{item:lemma_intrinsic_rolling_regular_distribution}, i.e.
	\begin{equation}
		\begin{split}
			\Psi^{\alpha}  \colon \Qspace \times \liealg{m}
			&\to T \Qspace, \\
			((v, [g, S]), u) 
			&\mapsto 
			\Psi^{\alpha}((v, [g, S]), u)
			= (u, [(T_e \ell_g \circ  S) u, - \alpha(S u, \cdot) \circ S]) .
		\end{split}
	\end{equation}
	Then $\Psi^{\alpha}$ defines a control system in the sense
	of~\cite[p. 21]{jurdjevic:1997}
	with state space $\Qspace$ and control set $\liealg{m}$.
	Obviously, for each $u \in \liealg{m}$,
	the map $\Psi^{\alpha}(\cdot, u) \colon \Qspace \to T \Qspace$
	is a section of $D^{\alpha}$, where $D^{\alpha} \subseteq T \Qspace$
	is the distribution characterizing the rolling of
	$\liealg{m}$ over $G / H$ with respect to $\nablaAlpha$.
\end{remark}
Moreover, if $G / H$ is equipped with an invariant
pseudo-Riemannian metric and an invariant metric covariant derivative $\nablaAlpha$, we can endow $\Qspace$ with an additional
structure which is similar to a sub-Riemannian structure.
We refer to~\cite[Def. 3.2]{agrachev.barilari.boscain:2020}
for a definition of sub-Riemannian structures.

\begin{remark}
	Let $G / H$ be a reductive homogeneous space
	equipped with a $G$-invariant pseudo-Riemannian metric
	corresponding to the scalar product
	$\langle \cdot, \cdot \rangle \colon \liealg{m} \times \liealg{m} \to \field{R}$.
	Moreover, let
	$\alpha \colon \liealg{m} \times \liealg{m} \to \liealg{m}$
	be an $\Ad(H)$-invariant bilinear map defining the \emph{metric}
	invariant covariant derivative $\nabla^{\alpha}$.
	As in Remark~\ref{remark:intrinsic_rolling_pseudo_riemannian_metric_Qspace_reduced},
	we
	set $\Qspace = \liealg{m} \times (G \times_H \liegroup{O}(\liealg{m}))$.
	Moreover, motivated by~\cite[Eq. (3)]{grong:2012},
	we equip the trivial vector bundle
	$\Qspace \times \liealg{m} \to \Qspace$ with the fiber metric
	$h \in \Secinfty\big(\Sym^2 (\Qspace \times \liealg{m})^* \big)$
	defined by $h_q(X, Y) = \langle X, Y \rangle$ for $q \in \Qspace$
	and $X, Y \in \liealg{m}$.
	Then the pair $(\Psi^{\alpha}, \Qspace \times \liealg{m})$
	is formally similar to a sub-Riemannian structure on $\Qspace$ 
	except for the following facts:
	\begin{enumerate}
		\item
		The fiber metric on $\Qspace \times \liealg{m} \to \liealg{m}$ is allowed to be indefinite.
		\item
		In general, the manifold $\Qspace$ is \emph{not} connected.
		\item
		The distribution $D^{\alpha} = \image(\Psi)$ might be
		\emph{not} bracket generating.
	\end{enumerate}
	However, by imposing further restrictions on $G / H$, $\Qspace$ and
	$\langle \cdot, \cdot \rangle$ it might be possible to obtain a
	sub-Riemannian structure on $\Qspace$.
	In particular, if we assume that $G / H$ is a Riemannian reductive
	homogeneous space, the fiber metric $h$ on $\Qspace$ is positive definite.
	Moreover, if we assume that $G$ is connected and
	$\Ad_h\at{\liealg{m}} \colon \liealg{m} \to \liealg{m}$
	is an orientation preserving isometry, i.e.
	$\Ad_h\at{\liealg{m}} \in \liegroup{SO}(\liealg{m})$
	for all $h \in H$, 
	the configuration space can be reduced to
	$\Qspace = \liealg{m} \times (G \times_H \liegroup{SO}(\liealg{m}))$,
	which is obviously connected.
	Under these assumptions, the pair
	$(\Psi^{\alpha}, \Qspace \times \liealg{m})$ defines a
	structure on $\Qspace$ which fulfills the requirements
	of a sub-Riemannian structure on $\Qspace$ in the sense
	of~\cite[Def. 3.2]{agrachev.barilari.boscain:2020}
	except for the 
	fact that $D^{\alpha}$ might be not bracket generating.
	Investigating conditions on $G / H$ and $\alpha$ such
	that $D^{\alpha}$ is bracket generating is out of the scope of this text.
	Nevertheless, in this context, we refer to~\cite{grong:2012},
	where the controllability of rollings of
	oriented Riemannian manifolds are considered.
	Moreover, we mention~\cite{jurdjevic:2022},
	where optimal control
	problems associated to rollings of certain manifolds
	are considered.
\end{remark}

Using terminologies of control theory, we call a curve
$u \colon I \to \liealg{m}$
a control curve. 
Such a curve can be used to determine a rolling 
of $\liealg{m}$ over $G / H$, 
where the rolling curve $v \colon I \to \liealg{m}$
satisfies the ODE $\dot{v}(t) = u(t)$.
Inspired by the terminology used in~\cite{hueper.leite:2007},
we introduce a notion of a kinematic equation
for rollings of $\liealg{m}$ over $G / H$
with respect to $\nablaAlpha$.
To this end, we first state the following proposition.

\begin{proposition}
	\label{proposition:kinematic_equation_intrinsic_rolling}
	Let $u \colon I \to \liealg{m}$ be a control curve
	and let $\alpha \colon \liealg{m} \times \liealg{m} \to \liealg{m}$
	be an $\Ad(H)$-invariant bilinear map satisfying
	$\alpha(X, \cdot) \in \liealgsub(\liealg{m})$.
	Moreover, let
	$\overline{q} \colon I \ni t \mapsto \overline{q}(t)
	=
	(v(t), g(t), S(t)) 
	\in \liealg{m} \times G \times \liegroupsub(\liealg{m}) = \QspaceLift$
	be a curve satisfying the ODE
	\begin{equation}
		\label{equation:proposition_kinematic_equation_intrinsic_rolling}
		\begin{split}
			\begin{split}
				\dot{v}(t)
				&= u(t) , \\
				\dot{S}(t) 
				&= - \alpha\big(S(t) u(t), \cdot \big) \circ S(t) , \\
				\dot{g}(t) 
				&= \big( T_e \ell_{g(t)} \circ S(t) \big) u(t) .
			\end{split}
		\end{split}
	\end{equation}
	Then
	$q \colon I \ni t \mapsto (\prQLift \circ \overline{q})(t) \in \Qspace$
	is a rolling of $\liealg{m}$ over $G / H$ with respect to $\nablaAlpha$.
	\begin{proof}
		The curve $q \colon I \ni t \mapsto (\prQLift \circ \overline{q})(t) 
		= (v(t), [g(t), S(t)]) \in \Qspace$
		is horizontal with respect to $D^{\alpha} \subseteq T \Qspace$
		by Theorem \ref{theorem:intrinsic_rolling},
		Claim~\ref{item:theorem_intrinsic_rolling_horizontal_lift}
		because of $u(t) = \dot{v}(t)$ for all $t \in I$.
		Hence $q \colon I \to Q$ is a rolling of $\liealg{m}$ over $G / H$ by
		Theorem~\ref{theorem:intrinsic_rolling},
		Claim~\ref{item:theorem_intrinsic_rolling_distribution_rolling}.
	\end{proof}
\end{proposition}

\begin{definition}
	\label{definition:kinematic_equations_intrinsic_rolling}
	The 
	ODE~\eqref{equation:proposition_kinematic_equation_intrinsic_rolling} 
	in
	Proposition~\ref{proposition:kinematic_equation_intrinsic_rolling}
	is called the kinematic equation for
	($\liegroupsub(\liealg{m})$-reduced) rollings of $\liealg{m}$
	over $G / H$ with respect to $\nablaAlpha$.
	An initial value problem associated
	with the ODE~\eqref{equation:proposition_kinematic_equation_intrinsic_rolling}
	with some initial condition
	$(v(t_0), g(t_0), S(t_0)) \in \liealg{m} \times G \times \liegroupsub(\liealg{m})$ for some $t_0 \in I$
	is called kinematic equation, as well.
\end{definition}
In the sequel, we are mainly interested in the initial value problem
associated with~\eqref{equation:proposition_kinematic_equation_intrinsic_rolling}
defined by the initial condition
$(v(0), g(0), S(0)) = (0, e, \id_{\liealg{m}}) 
\in \QspaceLift = \liealg{m} \times G \times \liegroupsub(\liealg{m})$.

\begin{remark}
	\label{remark:kinematic_equation_special_cases_intrinsic}
	By specializing $\nablaAlpha$ in Definition~\ref{definition:kinematic_equations_intrinsic_rolling},
	one obtains:
	\begin{enumerate}
		\item
		\label{item:remark_kinematic_equation_special_cases_intrinsic_nabla_Can}
		The kinematic equation with respect to $\nablaCan$ reads
		\begin{equation}
			\label{equation:remark:kinematic_equation_special_cases_intrinsic_nabla_Can_1}
			\begin{split}
				\dot{v}(t)
				&= u(t), \\
				\dot{S}(t) 
				&=
				- \tfrac{1}{2} \pr_{\liealg{m}} \circ \ad_{S(t) u(t)} \circ S(t), \\
				\dot{g}(t) 
				&=
				\big( T_e \ell_{g(t)} \circ S(t) \big) u(t) .
			\end{split}
		\end{equation}
		\item
		\label{item:remark_kinematic_equation_special_cases_intrinsic_nabla_Can_Second}
		For $\nablaCanSecond$ one obtains the kinematic equation
		\begin{equation}
			\begin{split}
				\dot{v}(t)
				&= u(t) , \\
				\dot{S}(t) 
				&=  0 , \\
				\dot{g}(t) 
				&= \big( T_e \ell_{g(t)} \circ S(0) \big) u(t)
			\end{split}
		\end{equation}
		since $\dot{S}(t) = 0$ implies 
		$S(t) = S(t_0)$ for all $t \in I$ and some $t_0 \in I$.
		We point out that setting $S(t) = \id_{\liealg{m}}$ for all $t \in I$
		yields an expression which is similar to the ODE describing
		rollings of a symmetric space over a flat space
		obtained~\cite[Sec. 4.2]{jurdjevic.markina.leite:2023}.
	\end{enumerate}
\end{remark}
Next, we state the kinematic equation for a naturally
reductive homogeneous space.

\begin{corollary}
	\label{corollary:intrinsic_rolling_kinematic_equation_naturally_redudctive_homogeneous_spaces}
	Let $G / H$ be a naturally reductive homogeneous space
	and let
	$\langle \cdot, \cdot \rangle \colon \liealg{m} \times \liealg{m} \to \field{R}$
	be the $\Ad(H)$-invariant scalar product
	corresponding to the pseudo-Riemannian metric.
	Then the kinematic equation for ($\liegroup{O}(\liealg{m})$-reduced) rollings of $\liealg{m}$
	over $G/H$ with respect to $\nablaLC$ is given
	by~\eqref{equation:remark:kinematic_equation_special_cases_intrinsic_nabla_Can_1}
	from Remark~\ref{remark:kinematic_equation_special_cases_intrinsic},
	Claim~\ref{item:remark_kinematic_equation_special_cases_intrinsic_nabla_Can}.
	In particular, the curve $S \colon I \to \liegroup{O}(\liealg{m})$ 
	takes values in the
	pseudo-Euclidean group of $\liealg{m}$ with respect to the
	$\langle \cdot, \cdot \rangle$ provided that the initial condition 
	lies in $\liegroup{O}(\liealg{m})$.
	\begin{proof}
		Since $G / H$ has an invariant pseudo-Riemannian metric corresponding
		to the $\Ad(H)$-invariant scalar product
		$\langle \cdot, \cdot \rangle$ on $\liealg{m}$, one obtains
		$\Ad_h \at{\liealg{m}} \in \liegroup{O}(\liealg{m})$ for all $h \in H$.
		Moreover,
		Remark~\ref{remark:kinematic_equation_special_cases_intrinsic},
		Claim~\ref{item:remark_kinematic_equation_special_cases_intrinsic_nabla_Can}
		yields the desired result since
		$\nablaLC$ is metric and $\nablaLC = \nablaCan$ holds
		for naturally reductive homogeneous spaces by
		Remark~\ref{remark:horizontal_lift_covariant_derivative_naturally_redutive_and_symmetric}.
	\end{proof}
\end{corollary}

\begin{remark}
	\label{remark:completeness_of_kinematic_equations_riemannian_case}
	In general, it is not clear to us whether 
	the maximal solution of the initial value problem
	\begin{equation}
		\label{remark:remark_completeness_of_kinematic_equations_riemannian_case_time_dependent_ODE}
		\dot{S}(t) 
		= - \alpha\big(S(t) u(t),  \cdot \big) \circ S(t),
		\quad
		S(0) =\id_{\liealg{m}}
	\end{equation}
	associated with the kinematic equation
	in Proposition~\ref{proposition:kinematic_equation_intrinsic_rolling}
	is defined on the whole interval $I$.
	In principal, it could only be defined on a proper
	subinterval $I_1 \subsetneq I$.

	However, if we assume that $G / H$ is equipped with an invariant
	Riemannian metric,
	i.e. an invariant positive definite metric, $\nablaAlpha$
	is metric and the control
	curve $u \colon \field{R}  \to \liealg{m}$
	is bounded,
	we can prove that the time-independent vector field
	on $\field{R} \times \liegroup{O}(\liealg{m})$
	associated
	to~\eqref{remark:remark_completeness_of_kinematic_equations_riemannian_case_time_dependent_ODE},
	see
	e.g.~\cite[Sec. 3.30]{michor:2008},
	given by
	\begin{equation}
		X(t, S) = \big(1,  - \alpha\big(S u(t), \cdot \big)  \circ S \big),
		\quad
		(t, S) \in \field{R} \times \liegroup{O}(\liealg{m})
	\end{equation}
	is complete.
	To this end, we show that this vector field is bounded
	in a complete Riemannian metric on $\field{R} \times \liegroup{O}(\liealg{m})$.
	Then the completeness
	of
	$X \in \Secinfty\big(T (\field{R} \times \liegroup{O}(\liealg{m})) \big)$
	follows by~\cite[Prop. 23.9]{michor:2008}.
	We view $\liegroup{O}(\liealg{m})$ as subset
	of $\End(\liealg{m})$
	and denote by
	$\langle \cdot, \cdot \rangle \colon \liealg{m} \times \liealg{m} \to \liealg{m}$
	the $\Ad(H)$-invariant
	inner product corresponding to the Riemannian metric on $G / H$.
	The norm on $\liealg{m}$ induced by
	$\langle \cdot, \cdot \rangle$ is denoted by $\Vert \cdot \Vert$.
	We denote an extension of these maps to $\liealg{g}$ by the same symbols.
	We now endow $\End(\liealg{m})$ with the Frobenius scalar product
	given by $\langle S, T \rangle_F = \tr(S^{\top} T)$,
	where $S^{\top}$ is the adjoint of $S$ with respect to $\langle \cdot, \cdot \rangle$.
	Then $\langle \cdot, \cdot \rangle_F$ induces
	a bi-invariant and hence a complete Riemannian
	metric on $\liegroup{O}(\liealg{m})$.
	Moreover, the norm $\Vert \cdot \Vert_F$ defined by the Frobenius
	scalar product is equivalent to the operator norm
	$\Vert \cdot \Vert_2$.
	In particular, there is a $C > 0$ such that $\Vert S \Vert_F \leq C \Vert S \Vert_2$
	holds for all $S \in \End(\liealg{m})$.
	We now endow $\field{R} \times \liegroup{O}(\liealg{m})$
	with the Riemannian metric defined
	for $(s, S) \in \field{R} \times \liegroup{O}(\liealg{m})$
	and
	$(v, V), (w, W) \in T_{(s, S)} ( \field{R} \times \liegroup{O}(\liealg{m}))$
	 by
	\begin{equation}
		\big\langle (v, V), (w, W) \big\rangle_{(s, S)}
		=
		v w + \tr\big(V^{\top} W \big),
	\end{equation}
	which is clearly complete.
	Moreover, $\alpha \colon \liealg{m} \times \liealg{m} \to \liealg{m}$ is
	bounded since $\liealg{m}$ is finite dimensional.
	Hence there exists a $C^{\prime} \geq 0$ with
	$\Vert \alpha(X, Y) \Vert \leq C^{\prime} \Vert X \Vert \Vert Y \Vert$
	for $X, Y \in \liealg{m}$.
	Thus, for fixed $X \in \liealg{m}$,
	the operator norm of the linear map
	$\alpha(X, \cdot) \colon \liealg{m} \to \liealg{m}$
	is bounded by
	\begin{equation}
		\Vert \alpha(X, \cdot) \Vert_2 \leq C^{\prime} \Vert X \Vert .
	\end{equation}
	By this notation, we obtain 
	\begin{equation}
		\label{equation:remark:completeness_of_kinematic_equations_riemannian_case_estimation_vf}
		\begin{split}
			\Vert X(t, S) \Vert^2
			&=
			1 + \big\Vert \alpha\big(S u(t), \cdot \big) \circ S \big\Vert_F^2 \\
			&\leq
			1 + C^2 \big\Vert \alpha\big(S u(t), \cdot\big) \circ S \big\Vert_2^2 \\
			&\leq
			1 + C^2 \big\Vert \alpha \big(S u(t), \cdot \big) \big\Vert_2^2 \Vert S \Vert_2^2 \\
			&\leq
			1 + (C C^{\prime})^2  \Vert S \Vert_2^2 \Vert u(t) \Vert^2 \\
			&\leq
			1+ (C C^{\prime})^2 \Vert u \Vert_{\infty}^2 < \infty,
		\end{split}
	\end{equation}
	where we exploited
	$\Vert S \Vert_2 = 1$ for all $S \in \liegroup{O}(\liealg{m})$
	and $\Vert u \Vert_{\infty}$ denotes the supremum norm of $u \colon I \to \liealg{m}$.
	Equation~\eqref{equation:remark:completeness_of_kinematic_equations_riemannian_case_estimation_vf}
	shows that
	$X \in \Secinfty\big(T ( \field{R} \times \liegroup{O}(\liealg{m})) \big)$
	is bounded in a complete Riemannian metric
	on $\field{R} \times \liegroup{O}(\liealg{m})$ as desired.
	Thus the maximal solution of the initial value problem
	\begin{equation}
		\dot{S}(t) = - \alpha(S(t) u(t), \cdot) \circ S(t),
		\quad
		S(0) = S_0 \in \liegroup{O}(\liealg{m})
	\end{equation}
	is defined on $\field{R}$.
\end{remark}

\subsection{Rolling along Special Curves}
\label{subsec:rolling_along_special_curves}

Next we consider a rolling of $\liealg{m}$ over $G / H$ 
along a curve such that the development curve $\gamma \colon I \to G / H$
is the projection of a not necessarily horizontal one-parameter subgroup
of $G$, i.e.
\begin{equation}
	\gamma(t) = \pr( \exp(t \xi) ),
	\quad
	t \in I
\end{equation}
for some $\xi  \in \liealg{g}$.
In this subsection, we focus on
the invariant covariant derivatives
$\nablaCan$ and $\nablaCanSecond$ on $G / H$.
This discussion is motivated by the rolling
and unwrapping technique
for solving interpolation problems, see
e.g.~\cite{hueper.leite:2007,hueper.krakowski.leite:2020},
for which an explicit expression for a 
rolling along a curve connecting two points
is desirable.
A natural choice for such a curve would be
a projection of a horizontal one-parameter subgroup in $G$,
i.e. a geodesic with respect to $\nablaCan$ or $\nablaCanSecond$.
However, even if such a curve connecting two given points exists,
as far as we know, in general,
no closed-formula for such curves are known.
In this context, we refer to~\cite{kakowski.machado.leite:2017},
where the problem of
connecting two points $X, X_1 \in \Stiefel{n}{k} \subseteq \matR{n}{k}$
on the Stiefel manifold $\Stiefel{n}{k}$ by a curve of the form
$t \mapsto \e^{t \xi_1} X \e^{t \xi_2}$ 
with some suitable $(\xi_1, \xi_2) \in \liealg{so}(n) \times \liealg{so}(k)$
is addressed.

As preparation for deriving the desired rollings,
we state the following lemma as preparation.

\begin{lemma}
	\label{lemma:horizontal_lift_projection_one_parameter_group}
	Let $G / H$ be a reductive homogeneous space with
	reductive decomposition
	$\liealg{g} = \liealg{h} \oplus \liealg{m}$,
	let $\xi \in \liealg{g}$
	and let $\gamma \colon I \ni t \mapsto \pr(\exp(t \xi)) \in G / H$.
	Then the curve 
	\begin{equation}
		g \colon I \to G,
		\quad
		t \mapsto g(t) = \exp(t \xi) \exp(-t \xi_{  \liealg{h}})
	\end{equation}
	is the horizontal lift of $\gamma$	
	through $g(0) = e$
	with respect to the principal connection from
	Proposition~\ref{proposition:principal_connection_reductive_homogeneous_space}.
	Moreover, $g$ is the solution of the initial value problem
	\begin{equation}
		\label{equation:lemma_horizontal_lift_projection_one_parameter_group_IVP}
		\dot{g}(t) =  T_e \ell_{g(t)}  \big(\Ad_{\exp(t \xi_{\liealg{h}})}(\xi_{\liealg{m}}) \big), 
		\quad
		g(0) = e. 
	\end{equation}
	\begin{proof}
		Obviously, $\gamma(t) = \pr(\exp(t \xi)) = \pr(g(t))$ holds for all $t \in I$ due to $\exp(- t \xi_{  \liealg{h}}) \in H$
		for all $t \in I$.
		Moreover, $g(0) = e$ is fulfilled.
		It remains to prove that $g \colon I \to G$ is horizontal.
		To this end,
		we compute by
		exploiting~\eqref{equation:tangentmap_group_multiplication}
		and~\eqref{equation:exponential_map_lie_group}
		\begin{equation*}
			\begin{split}
				\dot{g}(t) 
				&=
				\tfrac{\D}{\D t} \big( \exp(t \xi) \exp(- t \xi_{  \liealg{h}}) \big) \\
				&=
				T_{\exp(t \xi)} r_{\exp(- t \xi_{\liealg{h}})}  \tfrac{\D}{\D t} \exp(t \xi)
				+
				T_{\exp(- t \xi_{  \liealg{h}})} \ell_{\exp(t \xi)}  \tfrac{\D }{\D t} \exp(- t \xi_{  \liealg{h}}) \\
				&= 
				\big( T_{\exp(t \xi)} r_{\exp(- t \xi_{\liealg{h})}}  \circ  T_e \ell_{\exp(t \xi)} \big) \xi 
				+
				\big(T_{\exp(- t \xi_{  \liealg{h}})} \ell_{\exp(t \xi)} 
				\circ 
				T_e r_{\exp(- t \xi_{ \liealg{h}})} \big) (- \xi_{ \liealg{h}}) \\
				&=
				T_e \big( r_{\exp(- t \xi_{\liealg{h}})} \circ \ell_{\exp(t \xi)} \big) \xi
				- 
				T_e \big( \ell_{\exp(t \xi)} \circ r_{\exp(- t \xi_{  \liealg{h}})} \big) \xi_{  \liealg{h}} \\
				&=
				T_e \big(\ell_{\exp(t \xi)} \circ r_{\exp(- t \xi_{  \liealg{h}})} \big) (\xi - \xi_{\liealg{h}}) \\
				&=
				T_e \big(\ell_{\exp(t \xi)} \circ r_{\exp(- t \xi_{  \liealg{h}})} \big) \xi_{\liealg{m}} .
			\end{split}
		\end{equation*}
		Consequently, we obtain by the chain rule
		\begin{equation*}
			\begin{split}
				(T_e \ell_{g(t)})^{-1} \dot{g}(t)
				&=
				\big(T_e \ell_{\exp(t \xi) \exp(- t \xi_{  \liealg{h}}) }\big)^{-1}
				\dot{g}(t) \\
				&=
				\big( T_{\exp(- t \xi_{  \liealg{h}})} \ell_{\exp(t \xi)} \circ T_e \ell_{\exp(- t \xi_{  \liealg{h}})}\big)^{-1}
				\dot{g}(t) \\
				&=
				\big( T_e \ell_{\exp(- t \xi_{  \liealg{h}})} \big)^{-1} 
				\circ
				\big(  T_{\exp(- t \xi_{  \liealg{h}})} \ell_{\exp(t \xi)} \big)^{-1} 
				\dot{g}(t) \\
				&=
				\big(  T_e \ell_{\exp(- t \xi_{  \liealg{h}})} \big)^{-1} 
				\big( T_{\exp(t \xi) \exp(- t \xi_{\liealg{h})} } \ell_{\exp( - t \xi)} \big)
				\big(  T_e \big(\ell_{\exp(t \xi)} \circ r_{\exp(- t \xi_{  \liealg{h}})} \big) \xi_{\liealg{m}} \big) \\
				&=
				\big( T_e \ell_{\exp(- t \xi_{  \liealg{h}})} \big)^{-1} 
				T_e \big( \ell_{\exp(- t \xi) } \circ \ell_{\exp(t \xi)}  \circ r_{\exp(- t \xi_{  \liealg{h}})}\big) \xi_{\liealg{m}} \\
				&=
				\big( T_e \ell_{\exp(- t \xi_{  \liealg{h}})} \big)^{-1} 
				T_e \big( r_{\exp(- t \xi_{  \liealg{h}})}\big) \xi_{\liealg{m}} \\
				&= \big(T_{\exp(- t \xi_{\liealg{h}})} \ell_{\exp( t \xi_{  \liealg{h}})} 
				\circ 
				T_e  r_{\exp(- t \xi_{  \liealg{h}})} \big) \xi_{\liealg{m}} \\
				&=
				T_e \big( \ell_{\exp( t \xi_{  \liealg{h}}) } \circ  r_{\exp(- t \xi_{  \liealg{h}})} \big)  \xi_{\liealg{m}} \\
				&=
				\Ad_{\exp(t \xi_{  \liealg{h}})}(\xi_{\liealg{m}}) .
			\end{split}
		\end{equation*}
		Since $G / H$ is reductive,
		this implies
		\begin{equation}
			\label{equation_lemma:horizontal_lift_projection_one_parameter_group_comp_2_conclusion}
			(T_e \ell_{g(t)})^{-1} \dot{g}(t) 
			=
			\Ad_{\exp(t \xi_{\liealg{h}})} (\xi_{\liealg{m}})
			\in \liealg{m} 
		\end{equation}
		due to $\exp(t \xi_{  \liealg{h}}) \in H$.
		Thus $g$ is horizontal.
		Moreover, the curve $g \colon I \to G$ is a solution of
		the initial value
		problem~\eqref{equation:lemma_horizontal_lift_projection_one_parameter_group_IVP}
		by~\eqref{equation_lemma:horizontal_lift_projection_one_parameter_group_comp_2_conclusion}.
	\end{proof}
\end{lemma}

\begin{remark}
	Let $G$ be equipped with a bi-invariant metric which induces
	a positive definite fiber metric on $\Hor(G)$,
	i.e. a sub-Riemannian structure on $G$,
	and let $H \subseteq G$ be a closed subgroup such that its
	Lie algebra $\liealg{h} \subseteq \liealg{g}$ is non-degenerated
	with respect to the scalar product corresponding to the metric.
	Then the curve given by
	$g(t) = \exp(t \xi)\exp(-t \xi_{\liealg{h}})$
	from Lemma~\ref{lemma:horizontal_lift_projection_one_parameter_group}
	is a sub-Riemannian geodesic on $G$
	according to~\cite[Sec. 11.3.7]{montgomery:2002}.
\end{remark}

\subsubsection{Rolling along Special Curves with respect to $\nablaCan$}
\label{subsubsec:rolling_along_special_curves}

We now derive an expression for a rolling of $\liealg{m}$ over $G / H$
with respect to $\nablaCan$
such that the development curve is given by
$\gamma \colon I \ni t \mapsto \pr(\exp(t \xi)) \in G / H$ 
for some $\xi \in \liealg{g}$.
To this end, we determine a curve
\begin{equation}
	\overline{q} \colon I \to \QspaceLift
	= 
	\liealg{m} \times G \times \liegroup{GL}(\liealg{m}), 
	\quad t \mapsto (v(t), g(t), S(t)) 
\end{equation}
which is horizontal with respect to the principal
connection $\overline{\mathcal{P}}$
from
Proposition~\ref{proposition:principal_H_fiber_bundle_over_configuration_space},
Claim~\ref{item:proposition:principal_H_fiber_bundle_over_configuration_space_defi_connection}
such that $q = \prQLift \circ \overline{q} \colon I \to \Qspace$
is the desired rolling.
In particular,
\begin{equation}
	\pr(g(t) ) = \gamma(t) = \pr(\exp(t \xi))
\end{equation}
has to be fulfilled
and $q$ has to be tangent to the distribution
$D^{\alpha}$
by Theorem~\ref{theorem:intrinsic_rolling},
Claim~\ref{item:theorem_intrinsic_rolling_distribution_rolling}.
Thus $\overline{q} \colon I \to \QspaceLift$ has to be
tangent to $\overline{D^{\alpha}}$ by Lemma~\ref{lemma:intrinsic_rolling_distribution_first_properties},
Claim~\ref{item:lemma_intrinsic_rolling_curve_horizontal_lift_distributions}.
Furthermore, by
Proposition~\ref{proposition:principal_H_fiber_bundle_over_configuration_space},
Claim~\ref{item:proposition:principal_H_fiber_bundle_over_configuration_space_definition_horizontal_lift_curve},
the curve $g \colon I \to G$
fulfilling $\gamma = \pr \circ g$ 
is tangent to 
$\Hor(G) \subseteq T G$ from
Proposition~\ref{proposition:principal_connection_reductive_homogeneous_space},
i.e. $g$ is a horizontal lift of $\gamma$.
Hence Lemma~\ref{lemma:horizontal_lift_projection_one_parameter_group}
yields
\begin{equation}
	g(t) = \exp(t  \xi) \exp(- t \xi_{  \liealg{h}})
\end{equation}
for all $t \in I$
which fulfills the initial value problem
\begin{equation}
	\label{equation:rolling_allong_special_curves_g_ode}
	\dot{g}(t) 
	=
	T_e \ell_{g(t)} \Ad_{\exp(t \xi_{  \liealg{h}})} (\xi_{\liealg{m}}) ,
	\quad
	g(0) = e .
\end{equation}
Next we recall the kinematic equation
from Remark~\ref{remark:kinematic_equation_special_cases_intrinsic},
Claim~\ref{item:remark_kinematic_equation_special_cases_intrinsic_nabla_Can}
for convenience
\begin{equation}
	\label{equation:rolling_along_special_curves_kinematic_eq_recall}
	\begin{split}
		\dot{v}(t)
		&= u(t), \\
		\dot{S}(t) 
		&=
		- \tfrac{1}{2} \pr_{\liealg{m}} \circ \ad_{S(t) u(t)} \circ S(t), \\
		\dot{g}(t) 
		&=
		\big( T_e \ell_{g(t)} \circ S(t) \big) u(t) .
	\end{split}
\end{equation}
By
comparing~\eqref{equation:rolling_along_special_curves_kinematic_eq_recall}
with~\eqref{equation:rolling_allong_special_curves_g_ode},
we obtain
\begin{equation}
	\label{equation:rolling_along_special_curves_Su_equal_Ad_exp_xi_m}
	S(t) u(t) = \Ad_{\exp(t \xi_{  \liealg{h}})} (\xi_{\liealg{m}}) .
\end{equation}
Thus the ODE for $S \colon I \to \liegroup{GL}(\liealg{m})$
in~\eqref{equation:rolling_along_special_curves_kinematic_eq_recall}
becomes
\begin{equation}
	\label{equation:rolling_along_special_curves_ode_S}
	\dot{S}(t) 
	=
	- \tfrac{1}{2} \pr_{\liealg{m}} \circ \ad_{S(t) u(t)} \circ S(t) 
	=
	- \tfrac{1}{2} \pr_{\liealg{m}} \circ \ad_{\Ad_{\exp(t \xi_{\liealg{h}})}(\xi_{\liealg{m}})} \circ S(t) .
\end{equation}
In order to obtain the desired rolling, we need to solve the initial value problem associated with~\eqref{equation:rolling_along_special_curves_ode_S}
explicitly. This is the next lemma.

\begin{lemma}
	\label{lemma:solution_ODE_S_t}
	The solution of the initial value problem 
	\begin{equation}
		\dot{S}(t) 
		=
		- \tfrac{1}{2} \pr_{\liealg{m}} \circ \ad_{\Ad_{\exp(t \xi_{\liealg{h}})}(\xi_{\liealg{m}})} \circ S(t), 
		\quad S(0) = S_0 \in \liegroup{GL}(\liealg{m})
	\end{equation}
	is given by 
	\begin{equation}
		S \colon I \to \liegroup{GL}(\liealg{m}),
		\quad 
		t \mapsto
		\Ad_{\exp(t \xi_{\liealg{h}})} \circ
		\exp\Big( -  t \pr_{\liealg{m}} \circ \ad_{\xi_{\liealg{h}} + \tfrac{1}{2} \xi_{\liealg{m}} } \Big) \circ S_0
	\end{equation}
	\begin{proof}
		We make the following Ansatz.
		We set
		$S(t) = \Ad_{\exp(t \xi_{\liealg{h}})} \circ \widetilde{S}(t)$
		for all $t \in I$,
		where $\widetilde{S} \colon I \to O(\liealg{m})$ is given by
		\begin{equation*}
				\widetilde{S}(t) 
				=
				\exp\big( - t (\ad_{\xi_{\liealg{h}}} + \tfrac{1}{2} (\pr_{\liealg{m}} \circ \ad_{\xi_{\liealg{m}}} ) )\big) \circ S_0 
		\end{equation*}
		for $t \in I$ and some $S_0 \in O(\liealg{m})$. Obviously, 
		\begin{equation*}
			\dot{\widetilde{S}}(t) 
			=
			- \big(\ad_{\xi_{\liealg{h}}} + \tfrac{1}{2} (\pr_{\liealg{m}} \circ \ad_{\xi_{\liealg{m}}} ) \big) \circ \widetilde{S}(t)
		\end{equation*}
		holds for all $t \in I$.
		Using the well-known identity
		$\Ad_{\exp(t \xi_{\liealg{h}})}
		=
		\e^{t \ad_{\xi_{\liealg{h}}}}$,
		we compute
		\begin{equation*}
			\tfrac{\D}{\D t} \Ad_{\exp(t \xi_{\liealg{h}})}
			=
			\tfrac{\D}{\D t} \e^{t \ad_{\xi_{\liealg{h}}}}
			=
			\e^{t \ad_{\xi_{\liealg{h}}}} \circ \ad_{\xi_{\liealg{h}}}
			= 
			\Ad_{\exp(t \xi_{  \liealg{h}})} \circ \ad_{\xi_{\liealg{h}}}  .
		\end{equation*} 
		Thus we obtain
		\begin{equation*}
			\begin{split}
				\dot{S}(t) 
				&= 
				\tfrac{\D}{\D t} 
				\big(  \Ad_{\exp(t \xi_{\liealg{h}})} \circ \widetilde{S}(t) \big) \\
				&=
				\big( \tfrac{\D}{\D t} \Ad_{\exp(t \xi_{\liealg{h}})} \big) \circ \widetilde{S}(t) 
				+
				\Ad_{\exp(t \xi_{\liealg{h}})} \circ \dot{\widetilde{S}}(t) \\
				&=
				\Ad_{\exp(t \xi_{\liealg{h}})} 
				\circ \ad_{\xi_{\liealg{h}}} \circ \widetilde{S}(t)
				- \Ad_{\exp(t \xi_{\liealg{h}})} \circ 
				\big(  \ad_{\xi_{\liealg{h}}} 
				+ \tfrac{1}{2} (\pr_{\liealg{m}} \circ \ad_{\xi_{\liealg{m}}} )  \big) \circ \widetilde{S}(t) \\
				&= 
				- \tfrac{1}{2} \Ad_{\exp(t \xi_{\liealg{h}})} 
				\circ \pr_{\liealg{m}} \circ \ad_{\xi_{\liealg{m}}} \circ \widetilde{S}(t) \\
				&=
				- \tfrac{1}{2} \pr_{\liealg{m}} \circ
				\ad_{\Ad_{\exp(t \xi_{\liealg{h}})}(\xi_{\liealg{m}})} \circ \Ad_{\exp(t \xi_{\liealg{h}})} \circ \widetilde{S}(t) \\
				&= 
				- \tfrac{1}{2} \pr_{\liealg{m}} \circ 
				\ad_{\Ad_{\exp(t \xi_{\liealg{h}})}(\xi_{\liealg{m}})} \circ S(t),
			\end{split}
		\end{equation*}
		where we exploited that $\pr_{\liealg{m}}$ and $\Ad_h$ commutes
		and that $\Ad_h$ is a morphism of Lie algebras
		for all $h \in H$.
		Moreover, $S(0) = \widetilde{S}(0) = S_0$ is clearly fulfilled. 
		Hence 
		\begin{equation*}
			S \colon I \to \liegroup{GL}(\liealg{m}),
			\quad 
			t \mapsto \Ad_{\exp(t \xi_{\liealg{h}})} \circ 
			\exp\big( - t (\ad_{\xi_{\liealg{h}}} + \tfrac{1}{2} (\pr_{\liealg{m}} \circ \ad_{\xi_{\liealg{m}}} ) )\big) \circ S_0
		\end{equation*}
		is the unique solution of the initial value problem 
		fulfilling $S(0) = S_0$.
		The desired result follows by
		\begin{equation*}
			\ad_{\xi_{\liealg{h}}} + \tfrac{1}{2} \pr_{\liealg{m}} \circ \ad_{\xi_{\liealg{m}}} \at[\Big]{\liealg{m}}
			=
			\pr_{\liealg{m}} \circ \ad_{\xi_{\liealg{h}} + \tfrac{1}{2} \ad_{\xi_{\liealg{m}}}} \at[\Big]{\liealg{m}},
		\end{equation*}
		where $[\liealg{h}, \liealg{m}] \subseteq \liealg{m}$
		is exploited.
	\end{proof}
\end{lemma}
We now choose $S(0) = \id_{\liealg{m}}$ in the expression for $S(t)$ 
from Lemma \ref{lemma:solution_ODE_S_t} and obtain 
\begin{equation}
	S(t) 
	=
	\Ad_{\exp(t \xi_{\liealg{h}})} \circ
	\exp\Big( -  t \pr_{\liealg{m}} \circ \ad_{\xi_{\liealg{h}} + \tfrac{1}{2} \xi_{\liealg{m}} } \Big), 
	\quad 
	t \in I .
\end{equation}
Clearly, the inverse of $S(t)$ is given by 
\begin{equation}
	S(t)^{-1} 
	=
	\exp\big( t \pr_{\liealg{m}} \circ \ad_{\xi_{\liealg{h}} + \tfrac{1}{2} \xi_{\liealg{m}} } \big)
	\circ \Ad_{\exp(- t \xi_{\liealg{h}})} .
\end{equation}
By~\eqref{equation:rolling_along_special_curves_Su_equal_Ad_exp_xi_m}, i.e.
$S(t) u(t) = \Ad_{\exp(t \xi_{  \liealg{h}})}(\xi_{\liealg{m}})$,
we have
\begin{equation}
		u(t) 
		=
		S(t)^{-1} \big( \Ad_{\exp(t \xi_{  \liealg{h}})}(\xi_{\liealg{m}}) \big) 
		=
		\exp\Big( t  \pr_{\liealg{m}} \circ \ad_{\xi_{\liealg{h}} + \tfrac{1}{2} \xi_{\liealg{m}}}
		\Big) (\xi_{\liealg{m}}
		) .
\end{equation}
According to~\eqref{equation:rolling_along_special_curves_kinematic_eq_recall},
the rolling curve
$v \colon I \to \liealg{m}$ is defined by $\dot{v}(t) = u(t)$.
Choosing $v(0) = 0$ as initial condition yields
\begin{equation}
	v(t) 
	=
	\int_{0}^{t}
	\exp\Big( s  \pr_{\liealg{m}} \circ \ad_{\xi_{\liealg{h}} + \tfrac{1}{2} \xi_{\liealg{m}}}
	\Big) (\xi_{\liealg{m}}) \D s ,
	\quad t \in I .
\end{equation}
Thus the desired rolling is determined.
The discussion above is summarized in the next proposition.

\begin{proposition}
	\label{proposition:rolling_along_special_curve}
	Let $G / H$ be a reductive homogeneous space
	with reductive decomposition
	$\liealg{g} = \liealg{h} \oplus \liealg{m}$
	and let $\xi \in \liealg{g}$ be arbitrary.
	Moreover, let $\overline{q} \colon I \ni t \mapsto \overline{q}(t)
	= (v(t), g(t), S(t)) \in \liealg{m} \times G \times \liegroup{GL}(\liealg{m}) = \QspaceLift
	$
	be defined by
	\begin{equation}
		\begin{split}
			v(t)
			&=
			\int_{0}^{t}
			\exp\Big( s  \pr_{\liealg{m}} \circ \ad_{\xi_{\liealg{h}} + \tfrac{1}{2} \xi_{\liealg{m}}}
			\Big) (\xi_{\liealg{m}}) \D s \\
			g(t) 
			&= 
			\exp(t \xi) \exp(- t \xi_{\liealg{h}}) \\
			S(t)
			&=
			\Ad_{\exp(t \xi_{\liealg{h}})} \circ
			\exp\Big( -  t \pr_{\liealg{m}} \circ \ad_{\xi_{\liealg{h}} + \tfrac{1}{2} \xi_{\liealg{m}} } \Big)
		\end{split}
	\end{equation}
	for $t \in I$.
	Then $q \colon I \ni t \mapsto (\prQLift \circ \overline{q})(t) 
	= (v(t), [g(t), S(t)]) \in \Qspace$
	defines a rolling of $\liealg{m}$ over $G / H$ 
	with rolling curve $v \colon I \to \liealg{m}$
	and development curve
	\begin{equation}
		\gamma \colon I \to G / H,
		\quad
		t \mapsto \pr(g(t)) = \pr(\exp(t \xi)) .
	\end{equation}
	Furthermore, this intrinsic rolling viewed as a triple as in
	Remark~\ref{remark:intrinsic_rolling_as_triple}
	is given by $(v(t), \gamma(t), A(t))$, where
	$A(t)$ reads
	\begin{equation}
		A(t) \colon T_{v(t)} \liealg{m} \cong \liealg{m} \to T_{\gamma(t)} (G  / H), 
		\quad
		Z \mapsto \big( T_{g(t)} \pr \circ T_e \ell_{g(t)} \circ S(t) \big) Z .
	\end{equation}
	\begin{proof}
		This is a consequent of the above discussion.
	\end{proof}
\end{proposition}

\begin{remark}
	Assume that
	$\liegroupsub(\liealg{m}) \subseteq \liegroup{GL}(\liealg{m})$
	is a closed subgroup such that
	$\Ad_h\at{\liealg{m}} \in \liegroupsub(\liealg{m})$ holds for
	all $h \in H$ and that
	$\pr_{\liealg{m}} \circ \ad_X \at{\liealg{m}} 
	\in \liegroupsub(\liealg{m})$ 
	holds for all $X \in \liealg{m}$.
	Then the curve $S \colon I \to \liegroup{GL}(\liealg{m})$
	from Proposition~\ref{proposition:rolling_along_special_curve}
	is actually contained in $\liegroupsub(\liealg{m})$, i.e.
	$S(t) \in \liegroupsub(\liealg{m})$ for all $t \in I$.
	In particular, if $G / H$ is a naturally
	reductive homogeneous space,
	ones has $S(t) \in \liegroup{O}(\liealg{m})$
	for all $t \in I$.
\end{remark}

\begin{corollary}
	Let $\xi_{\liealg{m}} \in \liealg{m}$
	and define the curves
	\begin{equation}
		\begin{split}
			&v \colon I \to \liealg{m}, 
			\quad t \mapsto v(t) = t  \xi_{\liealg{m}} \\
			&S \colon I \to \liegroup{GL}(\liealg{m}), \quad
			t \mapsto S(t) = \exp(- \tfrac{1}{2} t \pr_{\liealg{m}} \circ \ad_{\xi_{\liealg{m}}} ) \\
			&g \colon I \to G, 
			\quad t \mapsto g(t) = \exp(t  \xi_{\liealg{m}}) 
		\end{split}
	\end{equation}
	Then $q \colon I \ni t \mapsto q(t) = (v(t), [g(t), S(t)]) \in Q$
	is an intrinsic rolling
	with respect to $\nablaCan$
	whose development curve is a geodesic with respect to $\nablaCan$.
	\begin{proof}
		This is an immediate consequence of
		Proposition~\ref{proposition:rolling_along_special_curve}
		due to $\xi_{\liealg{h}} = 0$.
	\end{proof}
\end{corollary}

\subsubsection{Rolling along Special Curves with respect to $\nablaCanSecond$}

We now consider a rolling
of a reductive homogeneous space with respect to the
covariant derivative $\nablaCanSecond$
such that the development curve is given by
$\gamma \colon I \ni t \mapsto \pr(\exp(t \xi))$ for some $\xi \in \liealg{g}$.
This is the next proposition.

\begin{proposition}
	\label{proposition:rolling_along_special_curves_nabla_second}
	Let $\xi \in \liealg{g}$ be arbitrary and define
	$\overline{q} \colon I \ni t \mapsto (v(t), g(t), S(t)) \in  \liealg{m} \times G \times \liegroup{GL}(\liealg{m}) = \QspaceLift$
	by
	\begin{equation}
		\label{equation:proposition_rolling_along_special_curves_nabla_second_representatives}
		\begin{split}
			v(t) 
			&=
			\int_{0}^t \Ad_{\exp(s \xi_{\liealg{h}})}(\xi_{\liealg{m}}) \D s 		\\	
			S(t) 
			&= 
			\id_{\liealg{m}} \\
			g(t)
			&=
			\exp(t \xi) \exp(- t \xi_{ \liealg{h}})	.
		\end{split}
	\end{equation}
	Then $q = \prQLift \circ \overline{q} \colon I  \to \Qspace$
	is an intrinsic rolling of $\liealg{m}$ over $G / H$ with
	respect to $\nablaCanSecond$ whose development curve is given by 
	$\gamma  \colon I \ni t \mapsto  \pr(\exp(t \xi)) \in G / H$.
	This intrinsic rolling viewed as a triple as in
	Remark~\ref{remark:intrinsic_rolling_as_triple}
	is given by $(v(t), \gamma(t), A(t))$, where
	$A(t)$ reads
	\begin{equation}
		A(t) \colon T_{v(t)} \liealg{m} \cong \liealg{m} \to T_{\gamma(t)} (G  / H), 
		\quad
		Z \mapsto \big( T_{g(t)} \pr \circ T_e \ell_{g(t)} \big) Z .
	\end{equation}
	\begin{proof}
		The curve
		$g \colon I \ni t \mapsto \exp(t \xi) \exp(- t \xi_{\liealg{h}}) \in G$
		is the horizontal lift of $\gamma$ through $g(0) = e$ by
		Lemma~\ref{lemma:horizontal_lift_projection_one_parameter_group}.
		We now show that $\overline{q} \colon I \to \QspaceLift$ defined
		by~\eqref{equation:proposition_rolling_along_special_curves_nabla_second_representatives}
		fulfills the kinematic equation from
		Remark~\ref{remark:kinematic_equation_special_cases_intrinsic},
		Claim~\ref{item:remark_kinematic_equation_special_cases_intrinsic_nabla_Can_Second}.
		Indeed, we have
		\begin{equation*}
			u(t) 
			= 
			\dot{v}(t) 
			=
			\Ad_{\exp(t \xi_{ \liealg{h}})}(\xi_{\liealg{m}}) 
		\end{equation*}
		and $g(t) = \exp(t \xi) \exp(- t \xi_{\liealg{h}})$ is the solution of the initial value problem
		\begin{equation*}
			\dot{g}(t) 
			=
			T_e \ell_{g(t)} \Ad_{\exp(t \xi_{ \liealg{h}})}(\xi_{\liealg{m}})
			=  
			T_e \ell_{g(t)} u(t)
			=
			\big( T_e \ell_{g(t)} \circ \id_{\liealg{m}} \big) u(t) ,
			\quad
			g(0) = e
		\end{equation*}
		by Lemma~\ref{lemma:horizontal_lift_projection_one_parameter_group}
		as desired.
		Therefore $q = \prQLift \circ \overline{q} \colon I \to \Qspace$
		is indeed a
		rolling of $\liealg{m}$ over $G / H$ with respect to
		$\nablaCanSecond$ whose  developement curve is given
		by $\gamma(t) = \pr(\exp(t \xi))$.
	\end{proof}
\end{proposition}

\section{Applications and Examples}
\label{sec:applications_and_examples}

In this section, we consider some examples.
Before we study Lie groups and Stiefel manifolds
in detail, we briefly comment on symmetric homogeneous spaces.
By recalling $\nablaCan = \nablaCanSecond$ 
for symmetric homogeneous space and $\nablaLC = \nablaCan = \nablaCanSecond$
for pseudo-Riemannian symmetric homogeneous spaces from
Remark~\ref{remark:symmetric_spaces_invariant_covariant_derivatives},
the kinematic equation from
Remark~\ref{remark:kinematic_equation_special_cases_intrinsic}, 
Claim~\ref{item:remark_kinematic_equation_special_cases_intrinsic_nabla_Can_Second}
yields the next lemma.

\begin{lemma}
	\label{lemma:rolling_symmetric_spaces}
	Let $(G, H, \sigma)$ be a symmetric pair and let $G / H$
	be the corresponding symmetric homogeneous space
	with canonical reductive decomposition
	$\liealg{g} = \liealg{h} \oplus \liealg{m}$.
	Moreover, let $u \colon I \to \liealg{m}$ be a curve.
	Define the curve
	$\overline{q} \colon I \ni t \mapsto (v(t), g(t), \id_{\liealg{m}}) \in \QspaceLift
	= \liealg{m} \times G \times \liegroup{GL}(\liealg{m})$
	by the initial value problem
	\begin{equation}
		\label{equation:lemma_rolling_symmetric_spaces_kinematic_equation}
		\begin{split}
			\dot{v}(t) &= u(t), \qquad\quad \ \ v(0) = 0 \\
			\dot{g}(t)
			&= T_e \ell_{g(t)} u(t),  \quad g(0) = e.
		\end{split}
	\end{equation}
	Then $q = \prQLift \circ \overline{q} \colon 
	I \ni t \mapsto (v(t), [g(t), \id_{\liealg{m}}]) 
	\in \Qspace = 
	\liealg{m} \times (G \times_H \liegroup{GL}(\liealg{m}))$
	is a rolling of $\liealg{m}$ over $G / H$ with respect
	$\nablaCan = \nablaCanSecond$.
\end{lemma}
If $G / H$ is a pseudo-Riemannian symmetric space, we
can consider an $\liegroup{O}(\liealg{m})$-reduced rolling, i.e. we
can take $\Qspace = \liealg{m} \times (G \times \liegroup{O}(\liealg{m}))$
in Lemma~\ref{lemma:rolling_symmetric_spaces}
and $q \colon I \to \Qspace$ is a rolling of
$\liealg{m}$ over $G / H$ 
with respect to the covariant derivative
$\nablaLC = \nablaCan = \nablaCanSecond$.
In this case, Lemma~\ref{lemma:rolling_symmetric_spaces}
is very similar the result obtained
in~\cite[Sec. 4.2]{jurdjevic.markina.leite:2023}.

\subsection{Rolling Lie Groups}

In this subsection, we discuss intrinsic rollings of Lie groups.
First we discuss rollings of $\liealg{g}$ over $G$, 
where we view $G$ as the reductive homogeneous space $G / \{e\}$
equipped with the covariant derivative $\nablaCan$.
Afterwards, we discuss rollings of
a connected Lie group $G$ viewed as the symmetric
homogeneous space $(G \times G) / \Delta G$
equipped
with $\nablaCan = \nablaCanSecond$.
It turns out that both points of view are closely related.

\subsubsection{Rollings of Lie Groups as Reductive Homogeneous Spaces}

We first consider the rolling of a Lie-group $G$
viewed as a reductive homogeneous space $G / \{e\}$ equipped 
with the covariant derivative $\nablaCan$.
Obviously, the reductive decomposition is given by $\liealg{h} = \{ 0\}$
and $\liealg{m} = \liealg{g}$. 
Clearly, this implies
$\pr_{\liealg{m}} = \id_{\liealg{m}} = \id_{\liealg{g}}$.
Moreover, the configuration space becomes
\begin{equation}
	\Qspace
	= \liealg{g} \times (G \times_{\{e\}} \liegroup{GL}(\liealg{g})) = \liealg{g}  \times G \times \liegroup{GL}(\liealg{g}) .
\end{equation}
We now determine a rolling 
$q \colon I 
\ni t \mapsto (v(t), g(t), S(t))
\in Q = \liealg{g} \times G \times  \liegroup{GL}(\liealg{g})$
of $\liealg{g}$ over $G$
with respect to $\nablaCan = \nablaAlpha$,
where
\begin{equation}
	\alpha
	\colon 
	\liealg{g} \times \liealg{g} \to \field{R},
	\quad
	(X, Y) \mapsto
	\alpha(X, Y) = \tfrac{1}{2} [X, Y] .
\end{equation}
To this end, we solve the following initial value problem associated with the
kinematic equation 
from Remark~\ref{remark:kinematic_equation_special_cases_intrinsic},
Claim~\ref{item:remark_kinematic_equation_special_cases_intrinsic_nabla_Can}
\begin{equation}
	\label{equation:kinmetic_equation_Lie_group_as_nat_red_space}
	\begin{split}
		\dot{v}(t) 
		&=
		u(t) , \qquad\qquad\qquad\quad \ \, v(0) = 0 , \\
		\dot{S}(t) 
		&=
		- \tfrac{1}{2}
		\ad_{S(t) u(t)} \circ S(t) , \quad \, S(0) = \id_{\liealg{g}} , \\
		\dot{g}(t)
		&=
		\big(T_e \ell_{g(t)} \circ S(t) \big) u(t) , \, \quad g(0) = g_0 ,
	\end{split}
\end{equation}
where $u \colon I \to \liealg{g}$ denotes a
prescribed control curve.
Motivated by \cite[Sec. 3.2]{hueper.leite:2007},
where rollings of $\liegroup{SO}(n)$ over one if its affine tangent
spaces are determined
by using an extrinsic point of view, we make the following Ansatz.

We define the curves
$k \colon I \to G$ and $W \colon I \to G$
by the initial value problems
\begin{equation}
	\label{equation:rolling_lie_groups_ODE_k_W}
	\dot{k}(t) = \tfrac{1}{2}  T_e \ell_{k(t)}  u(t), 
	\quad k(0) = g_0
	\quad \text{ and } \quad	
	\dot{W}(t) = - \tfrac{1}{2} T_e \ell_{W(t)}  u(t) , 
	\quad W(0) = e .
\end{equation}
Moreover, we set
\begin{equation}
	\label{equation:definition_S_Lie_group_as_nat_red_space}
	S \colon I \to \liegroup{GL}(\liealg {g}),
	\quad
	t \mapsto S(t) = \Ad_{W(t)}
\end{equation}
as well as 
\begin{equation}
	\label{equation:definition_g_Lie_group_as_nat_red_space}
	g \colon I \to G,
	\quad 
	t \mapsto
	g(t) = k(t) W(t)^{-1} .
\end{equation}
Clearly, $S(0) = \Ad_{e} = \id_{\liealg{g}}$
and $g(0) = g_0 e^{-1} = g_0$ holds.
Next we show that $S \colon I \to \liegroup{GL}(\liealg{g})$ defined
by~\eqref{equation:definition_S_Lie_group_as_nat_red_space}
is a solution of~\eqref{equation:kinmetic_equation_Lie_group_as_nat_red_space}.
To this end, we calculate
\begin{equation}
	\begin{split}
		\label{equation:rolling_liegroups_derivative_W_curve}
		\dot{W}(t) 
		=
		-\tfrac{1}{2} T_e \ell_{W(t)}  u(t) 
		=
		\tfrac{\D}{\D s}\big( \ell_{W(t)} \big( \exp(- \tfrac{1}{2} s u(t)) \big) \big) \at{s = 0} ,
	\end{split}
\end{equation}
where we used the chain-rule and
exploited the definition of $W$
in~\eqref{equation:rolling_lie_groups_ODE_k_W}.
In other words, the smooth curve
\begin{equation}
	\label{equation:definition_gamma_Lie_group_as_nat_red_space}
	\gamma \colon \field{R} \to G, 
	\quad
	s \mapsto
	\ell_{W(t)}\big(\exp(- \tfrac{1}{2} s u(t)) \big)
	=
	W(t) \exp(- \tfrac{1}{2} s u(t))
\end{equation}
fulfills $\gamma(0) = W(t)$
and $\tfrac{\D}{\D s} \gamma(s) \at{s = 0} = \dot{W}(t)$
according to~\ref{equation:rolling_liegroups_derivative_W_curve}.
Thus we calculate for $Z \in \liealg{g}$
by using the definition of $S \colon I \to \liegroup{GL}(\liealg{g})$
from~\eqref{equation:definition_S_Lie_group_as_nat_red_space}
and the chain rule
\begin{equation}
	\begin{split}
		\dot{S}(t)Z
		&=
		\tfrac{\D}{\D t} \big( \Ad_{W(t)} (Z) \big) \\
		&=
		T_{W(t)} \Ad_{(\cdot)}(Z)  \dot{W}(t) \\
		&=
		\tfrac{\D}{\D s} \Ad_{\gamma(s)}(Z) \at{s = 0} \\
		&=
		\tfrac{\D}{\D s} \Ad_{W(t) \exp(- \tfrac{1}{2} s u(t))}(Z)\at{s =0} \\
		&=
		\tfrac{\D}{\D s} \Ad_{W(t)} \Big( \Ad_{\exp(- \tfrac{1}{2} s u(t))} (Z) \Big) \at{s = 0} \\
		&=
		\Ad_{W(t)} \Big( \tfrac{\D}{\D s} \Ad_{\exp(- \tfrac{1}{2} s u(t))}(Z)  \at{s = 0} \Big) \\
		&=
		- \tfrac{1}{2} \big( \Ad_{W(t)} \circ \ad_{u(t)}  \big) (Z) \\
		&=
		- \tfrac{1}{2} \ad_{\Ad_{W(t)}(u(t))} \circ \Ad_{W(t)} (Z) \\
		&=
		- \tfrac{1}{2} \ad_{S(t) u(t)} \circ S(t) Z
	\end{split}
\end{equation}
as desired, where we exploited that
$\Ad_g \colon \liealg{g} \to \liealg{g}$ is a morphism of Lie algebras
for all $g \in G$
and $\gamma \colon \field{R} \to G$ is defined
by~\eqref{equation:definition_gamma_Lie_group_as_nat_red_space}.

It remains to show that $g \colon I \to G$ defined
in~\eqref{equation:definition_g_Lie_group_as_nat_red_space}
fulfills~\eqref{equation:kinmetic_equation_Lie_group_as_nat_red_space}.
To this end, using the chain-rule several times,
we obtain by
$g(t) = k(t) W(t)$
and~\eqref{equation:tangentmap_group_multiplication}
as well as~\eqref{equation:tangentmap_group_inversion}
and
the definition of $k \colon I \to G$
and $W \colon I \to G$
in~\eqref{equation:rolling_lie_groups_ODE_k_W}
\begin{equation}
	\begin{split}
		\dot{g}(t)
		&= 
		\tfrac{\D}{\D t} \big(k(t) W(t)^{-1} \big) \\
		&=
		T_{k(t)} r_{W(t)^{-1}}  \dot{k}(t) 
		+  T_{W(t)^{-1}} \ell_{k(t)}  \tfrac{\D}{\D t} \inv(W(t)) \\
		&=
		T_{k(t)} r_{W(t)^{-1}}  \dot{k}(t) +
		T_{W(t)^{-1}} \ell_{k(t)}  
		\big( - T_e \ell_{W(t)^{-1}} \circ T_{W(t)} r_{W(t)^{-1}} \big) \dot{W}(t) \\
		&=
		\tfrac{1}{2} \big( T_{k(t)} r_{W(t)^{-1}}   \circ  T_e \ell_{k(t)} \big) u(t) \\
		&\quad
		+ \tfrac{1}{2} 
		\big(
		T_{W(t)^{-1}} \ell_{k(t)}   \circ 
		T_e \ell_{W(t)^{-1}} \circ T_{W(t)} r_{W(t)^{-1}}  \circ  T_e \ell_{W(t)} \big) u(t) \\
		&=
		\tfrac{1}{2}  \big( T_e (r_{W(t)^{-1}} \circ \ell_{k(t)}) \big) u(t) \\
		&\quad 
		+ \tfrac{1}{2} \big( T_{W(t)^{-1}} \ell_{k(t)}  \circ
		T_e (\ell_{W(t)^{-1}} \circ r_{W(t)^{-1}} \circ \ell_{W(t)}) \big) u(t) \\
		&=
		\tfrac{1}{2}  T_e (r_{W(t)^{-1}} \circ \ell_{k(t)})  u(t)
		+ 
		\tfrac{1}{2} \big( T_{W(t)^{-1}} \ell_{k(t)} 
		\circ  T_e r_{W(t)^{-1}} \big)  u(t) \\
		&= 
		\tfrac{1}{2}  T_e \big(r_{W(t)^{-1}} \circ \ell_{{g}(t) W(t)} \big)  u(t)
		+ \tfrac{1}{2} \big( T_{W(t)^{-1}} \ell_{g(t) W(t)}  \circ  T_e r_{W(t)^{-1}} \big)  u(t) \\
		&=
		\tfrac{1}{2}  T_e \big(r_{W(t)^{-1}} \circ \ell_{{g}(t)} \circ \ell_{W(t)} \big)  u(t) 
		+
		\tfrac{1}{2} 
		 T_{e} \big( \ell_{g(t)} \circ \ell_{W(t)} \circ r_{W(t)^{-1}}  \big)  u(t) \\
		&=
		\tfrac{1}{2}
		T_e \big(\ell_{{g}(t)} \circ \ell_{W(t)} \circ r_{W(t)^{-1}}   \big) u(t)
		+ 
		\tfrac{1}{2} 
		T_{e} \big( \ell_{g(t)} \circ \ell_{W(t)} \circ r_{W(t)^{-1}} \big)   u(t) \\
		&=  \big( T_e \ell_{g(t)} \circ  T_e \Conj_{W(t)} \big)  u(t) \\
		&=  T_e \ell_{g(t)} \circ \Ad_{W(t)} ( u(t) ) \\
		&= \big( T_e \ell_{{g}(t)} \circ  S(t) \big) u(t)
	\end{split}
\end{equation}
as desired.
Hence
\begin{equation}
	q \colon I \ni t \mapsto (v(t), g(t), S(t)) \in \Qspace
\end{equation}
is an intrinsic rolling of $\liealg{g}$ over $G / H$ 
with respect to $\nablaCan$, where $S$ and $g$ are defined
in~\eqref{equation:definition_S_Lie_group_as_nat_red_space}
and~\eqref{equation:definition_g_Lie_group_as_nat_red_space},
respectively.
Moreover,
$v \colon I \to \liealg{g}$ is 
determined by $\dot{v}(t) = u(t)$ and the initial value $v(0) = 0$.
We summarize the above discussion in the next proposition.

\begin{proposition}
	\label{proposition:Lie-groups-intrinsic-rolling}
	Let $G$ be a Lie group
	viewed as reductive homogeneous space $G / \{e\}$
	equipped with $\nablaCan$.
	Let $u \colon I \to \liealg{g}$ be some control curve and
	define $k \colon I \to G$ as well as $W \colon I \to G$ by
	the initial value problems
	\begin{equation}
		\label{equation:proposition_Lie-groups-intrinsic-rolling_ivp_k_W}
		\dot{k}(t) = \tfrac{1}{2}  T_e \ell_{k(t)}  u(t), 
		\quad k(0) = g_0
		\quad \text{ and } \quad
		\dot{W}(t) = - \tfrac{1}{2}  T_e \ell_{W(t)}  u(t), 
		\quad W(0) = e .
	\end{equation}
	Then
	\begin{equation}
		q \colon I \ni t \mapsto (v(t), g(t), S(t)) 
		=
		\big(v(t), k(t) W(t)^{-1}, \Ad_{W(t)} \big) 
		\in
		\liealg{g} \times G \times \liegroup{GL}(\liealg{g})
		=
		\Qspace
	\end{equation}
	is an intrinsic rolling of $\liealg{g}$ over $G$, 
	where the development curve $v \colon I \to \liealg{g}$ is defined by
	\begin{equation}
		v(t) = \int_{0}^t u(s) \D s 
	\end{equation}
	and the rolling curve is given by $g \colon I \ni t \mapsto k(t) W(t)^{-1} \in G$.
	This rolling can be viewed as a triple $(v(t), g(t), A(t))$ as in Remark~\ref{remark:intrinsic_rolling_as_triple},
	where the linear isomorphism
	$A(t) \colon T_{v(t)} \liealg{g} \cong \liealg{g} \to T_{g(t)} G$
	is given by
	\begin{equation}
		\label{equation:proposition:proposition_Lie-groups-intrinsic-rolling_linear_isomorphism}
		A(t) Z 
		=
		\big( T_e \ell_{g(t)} \circ \Ad_{W(t)} \big) Z
		=
		\big( T_e \ell_{k(t) W(t)^{-1}} \circ \Ad_{W(t)} \big) Z
	\end{equation}
	for all $Z \in \liealg{g}$.
\end{proposition}

\begin{remark}
	\label{remark:rolling_lie_groups_reductive_existence_of_solution_on_whole_intervall}
	Let $u \colon I \to \liealg{g}$ be a control curve.
	Then the intrinsic rolling
	$q \colon I \to \liealg{m} \times G \times \liegroup{GL}(\liealg{g})$
	of $\liealg{g}$ over $G$ with respect
	to $\nablaCan$
	is defined on the whole interval $I$
	by the form of the initial value problem
	in~\eqref{equation:proposition_Lie-groups-intrinsic-rolling_ivp_k_W}.
\end{remark}

\begin{corollary}
	\label{corollary:lie_groups_rolling_bi-invariant_metric}
	Let $G$ be a Lie group equipped with a bi-invariant
	pseudo-Riemannian metric.
	Then rollings of $\liealg{g}$ over $G$ with respect to $\nablaLC$ 
	with a prescribed control curve $u \colon I \to \liealg{g}$
	are given by
	Proposition~\ref{proposition:Lie-groups-intrinsic-rolling}.
	\begin{proof}
		Let $G$ be equipped with a pseudo-Riemannian bi-invariant metric.
		Then the corresponding scalar product
		$\langle \cdot, \cdot \rangle \colon \liealg{g} \times \liealg{g} \to \field{R}$ is
		$\Ad(G)$-invariant,
		see e.g.~\cite[Chap. 11, Prop. 9]{oneill:1983},
		i.e. $G / \{e\}$ is a naturally reductive homogeneous space.
		Thus we have $\nablaLC = \nablaCan$ by
		Remark~\ref{remark:horizontal_lift_covariant_derivative_naturally_redutive_and_symmetric}.
		This yields the desired result.
	\end{proof}
\end{corollary}

\begin{remark}
	For the special case $G = \liegroup{SO}(n) \subseteq \matR{n}{n}$,
	equipped with
	the bi-invariant metric induced by the Frobenius scalar product
	on $\matR{n}{n}$,
	expressions for extrinsic rollings considered as curves
	in the Euclidean group 
	are derived
	in~\cite[Thm. 3.2]{hueper.leite:2007}.
	Rollings of pseudo-orthogonal groups have been studied
	in~\cite{crouch.leite:2012}.
	The tangential part of these
	rollings is very similar to the result of
	Proposition~\ref{proposition:Lie-groups-intrinsic-rolling}.
	Indeed, the linear isomorphism defined by the rolling
	from Proposition~\ref{proposition:Lie-groups-intrinsic-rolling}
	in~\eqref{equation:proposition:proposition_Lie-groups-intrinsic-rolling_linear_isomorphism}
	simplifies for a matrix Lie
	group to
	\begin{equation}
		A(t) Z
		=
		\big( T_e \ell_{k(t) W(t)^{-1}} \circ \Ad_{W(t)} \big) Z
		=
		k(t) Z W(t)^{-1}
	\end{equation}
	for all $Z \in \liealg{g}$.
\end{remark}

\subsubsection{Rollings of Lie Groups as Symmetric Homogeneous Spaces}

We now identify $G$ with the symmetric homogeneous
spaces $(G \times G) / \Delta G$ and study the rolling of $\liealg{m}$
over $(G \times G) / \Delta G$
with respect to $\nablaCan = \nablaCanSecond$.
To this end, we state the next lemma as preparation
which is an adaption
of~\cite[Sec. 23.9.5]{gallier.quaintance:2020},
see also~\cite[Chap. IV, 6]{helgason:1978},
where it is stated for the Riemannian case.

\begin{lemma}
	\label{lemma:Lie_group_symmetric_homogeneous_space}
	Let $G$ be a connected Lie group and define
	\begin{equation}
		\sigma \colon G \times G \to G \times G, 
		\quad
		(g_1, g_2) \mapsto (g_2, g_1) .
	\end{equation}
	Then $\sigma$ is an involutive automorphism of $G \times G$ and $\Delta G = \{(g, g) \mid g \in G\} \subseteq G \times G$
	is the set of fixed points of $\sigma$.
	Moreover, $(G \times G) / \Delta G$ is a symmetric homogeneous space
	and the corresponding canonical reductive decomposition
	$\liealg{g} \times  \liealg{g} = \liealg{h} \oplus \liealg{m}$
	is given by
	\begin{equation}
		\liealg{h} = \{(X, X) \mid X \in \liealg{g}\}
		\quad \text{ and } \quad
		\liealg{m} = \{(X, -X) \mid X \in \liealg{g}\}.
	\end{equation}
	In addition, the map
	\begin{equation}
		\phi \colon (G \times G) / \Delta G \to G,
		\quad (g_1, g_2) \cdot \Delta G \mapsto g_1 g_2^{-1}
	\end{equation}
	is a diffeomorphism
	and the map 
	\begin{equation}
		\overline{\phi} \colon G \times G \to G,
		\quad
		(g_1, g_2) \mapsto g_1 g_2^{-1}
	\end{equation}
	is a surjective submersion which fulfills
	$\phi \circ \pr = \overline{\phi}$.
\end{lemma}

Next we determine the tangent map of $\overline{\phi}$
evaluated at elements in $\Hor(G \times G ) \subseteq T G \times TG$.
We point out that the identity
$\big( T_{(e, e)} \overline{\phi} \big)(X, X) = 2 X$
for all $X \in \liealg{g}$ is well-known,
see e.g.~\cite[Sec. 23.9.5]{gallier.quaintance:2020}
or~\cite[Chap. IV, 6]{helgason:1978}.

\begin{lemma}
	\label{lemma:lie_group_symmetric_space_diffeomorphism}
	Let $G$ be a connected Lie group and let $X \in \liealg{g}$.
	Then $(T_e \ell_{g_1} X, - T_e \ell_{g_2} X) \in \Hor(G \times G)_{(g_1,g_2)}$ holds.
	Moreover, the tangent map of
	$\overline{\phi} \colon G \times G \ni (g_1, g_2) \mapsto g_1 g_2^{-1} \in G$
	fulfills
	\begin{equation}
		\label{equation:lemma_lie_group_symmetric_space_diffeomorphism_tangent_map_overline_phi}
		\big( T_{(g_1, g_2)} \overline{\phi} \big)
		(T_e \ell_{g_1} X, - T_e \ell_{g_2} X)
		=
		\big(T_e \ell_{g_1 g_2^{-1}} \circ \Ad_{g_2} \big)(2 X) .
	\end{equation}
	In particular $T_{(e, e)} \overline{\phi}(X, -  X) = 2 X$ holds
	and 
	\begin{equation}
		\label{equation:lemma_lie_group_symmetric_space_diffeomorphism_inverse_at_identity}
		\big( T_{(e, e)} \overline{\phi} \at{\liealg{m}} \big)^{-1}(X)
		=
		\big( \tfrac{1}{2} X, - \tfrac{1}{2} X\big)
	\end{equation}
	is satisfied for all $X \in \liealg{g}$.
	\begin{proof}
		Obviously,
		$(T_e \ell_{g_1} X, - T_e \ell_{g_2} X)
		\in \Hor(G \times G)_{(g_1,g_2)}$ is satisfied by
		$\Hor(G \times G)_{(g_1, g_2)} 
		=
		T_{(e, e)} \ell_{(g_1, g_2)} \liealg{m}$
		and the definition of
		$\liealg{m} \subseteq \liealg{g} \times \liealg{g}$ in Lemma~\ref{lemma:Lie_group_symmetric_homogeneous_space}.
		Next we
		prove~\eqref{equation:lemma_lie_group_symmetric_space_diffeomorphism_tangent_map_overline_phi}.
		To this end, we consider the curve
		\begin{equation*}
			\gamma \colon \field{R} \to G \times G,
			\quad
			t \mapsto \big( g_1 \exp(t X), g_2 \exp(- t X))
		\end{equation*}
		which fulfills $\gamma(0) = (g_1, g_2)$ and
		$\dot{\gamma}(0) = (T_e \ell_{g_1} X, - T_e \ell_{g_2} X)$.
		Next we calculate
		\begin{equation*}
			\begin{split}
				\big( T_{(g_1, g_2)} \overline{\phi} \big)
				(T \ell_{g_1} X , -T_e \ell_{g_2} X)
				&=
				\tfrac{\D}{\D t} \overline{\phi} (\gamma(t)) \at{t = 0} \\
				&=
				\tfrac{\D}{\D t} \overline{\phi}\big( g_1 \exp(t X), g_2(\exp(- t X))\big) \at{t = 0} \\
				&= 
				\tfrac{\D}{\D t} g_1 \exp(t X) \big( g_2(\exp(- t X))\big)^{-1} \at{t = 0} \\
				&=
				\tfrac{\D}{\D t} g_1 \exp(t X) \exp(t X) g_2^{-1} \at{t = 0} \\
				&=
				\tfrac{\D}{\D t} g_1 \exp( 2 t X) g_2^{-1} \at{t = 0} \\
				&=
				\tfrac{\D}{\D t} g_1 g_2^{-1} g_{2} \exp( 2 t X) g_2^{-1} \at{t = 0} \\
				&=
				\tfrac{\D}{\D t} g_1 g_2^{-1} \Conj_{g_2} \big( \exp( 2 t X) \big) \at{t = 0} \\
				&=
				\tfrac{\D}{\D t} g_1 g_2^{-1} \exp\big(\Ad_{g_2}(2 t X) \big) \at{t = 0} \\
				&=
				T_e \ell_{g_1 g_2^{-1}} \tfrac{\D}{\D t} \exp\big(t \Ad_{g_2}(2 X)\big) \at{t = 0} \\
				&=
				T_e \ell_{g_1 g_2^{-1}}  \Ad_{g_2}(2 X) 
			\end{split}
		\end{equation*}
		proving~\eqref{equation:lemma_lie_group_symmetric_space_diffeomorphism_tangent_map_overline_phi}
		as desired.
		Evaluating~\eqref{equation:lemma_lie_group_symmetric_space_diffeomorphism_tangent_map_overline_phi}
		at $(g_1, g_2) = (e, e) \in G \times G$ yields
		\begin{equation*}
			T_{(e, e)} \overline{\phi} (X, -X) = (T_e \ell_e \circ \Ad_e) (2 X) = 2 X
		\end{equation*}
		for all $X \in \liealg{g}$.
		Now~\eqref{equation:lemma_lie_group_symmetric_space_diffeomorphism_inverse_at_identity}
		is verified by a straightforward calculation.
	\end{proof}
\end{lemma}

Next we consider intrinsic rollings of $\liealg{m}$
over $(G \times G) / \Delta G$ with respect to $\nablaCanSecond$ and
relate them to the intrinsic rollings
of $\liealg{g}$ over $G$ with respect to $\nablaCan$.
This is the next proposition.

\begin{proposition}
	\label{proposition:rolling_lie_groups_symmetric_space_reductive_space_relation}
	Let $G$ be a connected Lie group and let
	$u \colon I \to \liealg{g}$ be a control curve.
	Consider the initial value problem
	\begin{equation}
		\begin{split}
			\dot{v}(t)
			&= u(t), 
			\qquad\qquad\quad \ v(0) = 0, \\
			\dot{g}_1(t) 
			&= \tfrac{1}{2} T_e \ell_{g_1(t)} u(t), 
			\quad \ \ g_1(0) = g_0, \\
			\dot{g}_2(t)
			&= - \tfrac{1}{2} T_e \ell_{g_1(t)} u(t),
			\quad g_2(0) = e.
		\end{split}
	\end{equation}
	Then the following assertions are fulfilled:
	\begin{enumerate}
		\item
		\label{item:lemma:rolling_lie_groups_symmetric_space_reductive_space_relation_rolling_symmetric_space}
		The curve $\widetilde{q} \colon I \to \liealg{m} \times ( (G \times G) \times_{\Delta G} \liegroup{GL}(\liealg{m}))$ defined for $t \in I$ by
		\begin{equation}
			\widetilde{q}(t) 
			=
			\big( (\tfrac{1}{2} v(t), - \tfrac{1}{2} v(t)), 
			[(g_1(t), g_2(t)), \id_{\liealg{m}}] \big)
		\end{equation}
		is a rolling of $\liealg{m}$ over $(G \times G) / \Delta G$
		with respect to $\nablaCanSecond$
		with
		rolling curve $\widetilde{v} \colon I \ni t \mapsto \widetilde{v} = (\tfrac{1}{2} v(t), - \tfrac{1}{2} v(t)) \in \liealg{m}$
		and development curve $\widetilde{\gamma} \colon I \ni t  \mapsto \widetilde{\gamma}(t) = \pr(g_1(t), g_2(t)) \in (G \times G) / \Delta G$.
		\item
		\label{item:lemma:rolling_lie_groups_symmetric_space_reductive_space_relation_rolling_reductive_space}
		The curve
		\begin{equation}
			q \colon I \ni t \mapsto 
			\big(v(t), g_1(t) g_2(t)^{-1}, \Ad_{g_2(t)} \big)
			\in \liealg{g} \times G \times \liegroup{GL}(\liealg{g})
		\end{equation}
		is a rolling of $\liealg{g}$ over $G$ with respect to $\nablaCan$
		with rolling curve $v \colon I \to \liealg{g}$ and
		development curve
		$g \colon I \ni t \mapsto g(t) = g_1(t) g_2(t)^{-1} \in G$.
		\item
		\label{item:lemma:rolling_lie_groups_symmetric_space_reductive_space_relation_relation}
		Let $\phi \colon (G \times G) / \Delta G \to G$ be the diffeomorphism
		from Lemma~\ref{lemma:Lie_group_symmetric_homogeneous_space}.
		Then one has for all $Z \in \liealg{g} \cong T_{v(t)} \liealg{g}$
		and $t \in I$
		\begin{equation}
			q(t) Z 
			=
			T \phi \circ 
			\widetilde{q}(t) \circ 
			\big(T_{(e, e)} (\phi \circ \pr) \at{\liealg{m}} \big)^{-1} Z ,
		\end{equation}
		where $q(t)$ as well as $\widetilde{q}(t)$ are
		identified with the linear isomorphisms given by
		\begin{equation}
			q(t) \colon T_{v(t)} \liealg{g} \cong \liealg{g} \to T_{g(t)} G,
			\quad 
			Z \mapsto \big( T_{g(t)} \pr \circ T_e \ell_{g(t)} \circ \Ad_{g_2(t)} \big) Z,
		\end{equation}
		where $g(t) = g_1(t) g_2(t)^{-1}$, and 
		\begin{equation}
			\begin{split}
				\widetilde{q}(t) \colon T_{\widetilde{v}(t)} \liealg{m} \cong \liealg{m} 
				&\to T_{\pr(g_1(t), g_2(t))} (G \times G) / \Delta G,
				\\
				(Z, - Z) &\mapsto 
				\big( T_{(g_1(t), g_2(t))} \pr \circ (T_e \ell_{g_1(t)}, T_e \ell_{g_2(t)} ) \big) (Z, - Z),
			\end{split}
		\end{equation}
		respectively.
	\end{enumerate}
	\begin{proof}
		Claim~\ref{item:lemma:rolling_lie_groups_symmetric_space_reductive_space_relation_rolling_symmetric_space}
		follows by
		Remark~\ref{remark:kinematic_equation_special_cases_intrinsic},
		Claim~\ref{item:remark_kinematic_equation_special_cases_intrinsic_nabla_Can_Second}.
		Moreover,
		Claim~\ref{item:lemma:rolling_lie_groups_symmetric_space_reductive_space_relation_rolling_reductive_space}
		is a consequence of
		Proposition~\ref{proposition:Lie-groups-intrinsic-rolling}.
		
		It remains to show Claim~\ref{item:lemma:rolling_lie_groups_symmetric_space_reductive_space_relation_relation}.
		Let $Z \in \liealg{g}$.
		Then one has
		\begin{equation*}
			(T_{(e, e)} (\phi \circ \pr)\at{\liealg{m}})^{-1} Z
			=
			(T_{(e, e)} \overline{\phi}\at{\liealg{m}})^{-1} Z
			=
			\big(\tfrac{1}{2} Z, - \tfrac{1}{2} Z \big)
		\end{equation*}
		according to Lemma~\ref{lemma:lie_group_symmetric_space_diffeomorphism}.
		Moreover
		$\overline{\phi} = \phi \circ \pr$
		holds by
		Lemma~\ref{lemma:Lie_group_symmetric_homogeneous_space}
		implying $T \overline{\phi} = T \phi \circ T \pr$.
		Therefore we obtain by
		Lemma~\ref{lemma:lie_group_symmetric_space_diffeomorphism}
		\begin{equation*}
			\begin{split}
				T \phi \circ \widetilde{q}(t) \circ (T_{(e, e)} (\phi \circ \pr)\at{\liealg{m}})^{-1} Z
				&=
				T \phi \circ \widetilde{q}(t) 	\big(\tfrac{1}{2} Z, - \tfrac{1}{2} Z \big) \\
				&=
				T \phi \circ 	\big( T_{(g_1(t), g_2(t))} \pr \circ (T_e \ell_{g_1(t)}, T_e \ell_{g_2(t)} ) \big) \big(\tfrac{1}{2} Z, - \tfrac{1}{2} Z \big) \\
				&=
				T_{(g_1(t), g_2(t))} (\phi \circ \pr)  
				\big(\tfrac{1}{2} T_e \ell_{g_1(t)} Z, - \tfrac{1}{2} T_e \ell_{g_2(t)} Z \big) \\
				&=
				\big( T_{ (g_1(t), g_2(t) )} \overline{\phi} \big)
				\big(\tfrac{1}{2} T_e \ell_{g_1(t)} Z, - \tfrac{1}{2} T_e \ell_{g_2(t)} Z \big)  \\
				&=
				\big(T_e \ell_{ (g_1(t) g_2(t)^{-1})} \circ \Ad_{g_2(t)} \big) Z \\
				&=
				q(t) Z 
			\end{split}
		\end{equation*}
		for all $Z \in \liealg{g}$ as desired.
	\end{proof}
\end{proposition}

\subsection{Rolling Stiefel Manifolds}
\label{sec:rolling_stiefel_manifolds}

Rollings of Stiefel manifolds have been already
considered in the literature in~\cite{hueper.kleinsteuber.leite:2008}
and~\cite{jurdjevic.markina.leite:2023},
however not from an intrinsic point of view.
In this section we apply the general theory developed in
Section~\ref{sec:intrinsic_rolling_reductive_space}
to the Stiefel manifold 
$\Stiefel{n}{k}$ endowed with
the Levi-Civita covariant derivative defined by
a so-called $\alpha$-metric.
These metrics
have been recently introduce
in~\cite{hueper.markina.leite:2021}.

\begin{remark}
	We point out that in contrast to the previous sections,
	where $\alpha$ denotes a bilinear map
	$\liealg{m} \times \liealg{m} \to \liealg{m}$
	defining an invariant covariant derivative, in this section
	$\alpha$ denotes an element in $\field{R} \setminus \{0\}$.
	There is no danger of confusion because in this section, we consider
	rollings of Stiefel manifolds exclusively with respect to the
	Levi-Civita covariant derivative $\nablaLC$ defined by
	an $\alpha$-metric.
	Since the Stiefel manifold $\Stiefel{n}{k}$
	equipped with an $\alpha$-metric is a
	naturally reductive homogeneous space, see
	Lemma~\ref{lemma:stiefel_alpha_metric_naturally_reductive} below,
	the Levi-Civita covariant derivative
	$\nablaLC$ corresponds to
	the invariant
	covariant derivative defined by the bilinear map
	$\liealg{m} \times \liealg{m} 
	\ni (X, Y) \mapsto \tfrac{1}{2} [X, Y]_{\liealg{m}} 
	\in \liealg{m}$
	according to
	Remark~\ref{remark:horizontal_lift_covariant_derivative_naturally_redutive_and_symmetric}.
\end{remark}

\subsubsection{Stiefel Manifolds Equipped with $\alpha$-Metrics}

We start with recalling some results from~\cite{hueper.markina.leite:2021},
in particular \cite[Sec. 2-3]{hueper.markina.leite:2021}. 
The Stiefel manifold $\Stiefel{n}{k}$
can be considered as the embedded submanifold
of $\matR{n}{k}$ given by
\begin{equation}
	\Stiefel{n}{k}
	=
	\{ X \in \matR{n}{k} \mid X^{\top} X = I_k
	\}, 
	\quad 1 \leq k \leq n.
\end{equation}
In the sequel, we write
$\liegroup{O}(n) = \{R \in \matR{n}{n} \mid R^{\top} R = I_n \}$
for the orthogonal group.
We now identify $\Stiefel{n}{k}$ with a normal naturally
reductive space
$G / H$, where $G = \liegroup{O}(n) \times \liegroup{O}(k)$
is equipped with a suitable
biinvariant pseudo-Riemannian metric.
To this end, we consider the action of the Lie group
$G = \liegroup{O}(n) \times \liegroup{O}(k)$ on $\matR{n}{k}$
from the left via
\begin{equation}
	\begin{split}
		\Phi \colon (\liegroup{O}(n) \times \liegroup{O}(k)) \times \matR{n}{k} 
		\to \matR{n}{k},
		\quad
		((R, \theta), X) 
		\mapsto  \Phi\big((R, \theta), X \big) = R X \theta^{\top} .
	\end{split}
\end{equation}
For fixed $(R, \theta) \in \liegroup{O}(n) \times \liegroup{O}(k)$,
the induced diffeomorphism
\begin{equation}
	\Phi_{(R, \theta)} \colon \matR{n}{k} \to \matR{n}{k},
	\quad X \mapsto R X \theta^{\top}
\end{equation}
is clearly linear.
Restricting the second argument of $\Phi$ to $\Stiefel{n}{k}$
yields the action
\begin{equation}
	(\liegroup{O}(n) \times \liegroup{O}(k)) \times \Stiefel{n}{k} \to \Stiefel{n}{k},
	\quad 
	((R, \theta), X) \mapsto  \Phi\big((R, \theta), X \big) = R X \theta^{\top} ,
\end{equation}
which is known to be transitive.
This action is denoted by $\Phi$, as well.

Let $X \in \Stiefel{n}{k}$ be fixed and denote by $H = \Stab(X)$
the stabilizer subgroup of $X$ under the action $\Phi$.
We identify
$\Stiefel{n}{k} \cong (\liegroup{O}(n) \times \liegroup{O}(k) ) /H$
via the $(\liegroup{O}(n) \times \liegroup{O}(k))$-equivariant
diffeomorphism
\begin{equation}
	\label{equation:Stiefel_manifold_embedding}
	\iota_X \colon G / H \to \Stiefel{n}{k},
	\quad 
	(R, \theta) \cdot H \mapsto
	\Phi\big( (R, \theta), X \big)
	=
	R X \theta^{\top} ,
\end{equation}
where $ (R, \theta) \cdot H \in (\liegroup{O}(n) \times \liegroup{O}(k)) / H$ denotes the coset defined
by $(R, \theta) \in\liegroup{O}(n) \times \liegroup{O}(k)$.
Moreover, the map 
\begin{equation}
	\label{equation:definition_stiefel_pr_X}
	\pr_X \colon \liegroup{O}(n) \times \liegroup{O}(k) \to \Stiefel{n}{k},
	\quad
	(R, \theta) \mapsto R X \theta^{\top} 
\end{equation}
is a surjective submersion.
Note that $\iota_X \colon G / H  \to \Stiefel{n}{k} \subseteq \matR{n}{k}$
becomes a $(\liegroup{O}(n) \times \liegroup{O}(k))$-equivariant embedding
and
\begin{equation}
	\pr_X = \iota_X \circ \pr
\end{equation}
holds, where
$\pr \colon \liegroup{O}(n) \times \liegroup{O}(k) 
\to \big( \liegroup{O}(n) \times \liegroup{O}(k)\big) / H$
denotes the canonical projection.
The Lie algebra of $H$ is given by 
\begin{equation}
	\label{equation:stiefel_alpha_lie_algebra_h_at_X_defintion}
	\liealg{h} = \ker \big( T_{(I_n, I_k)} \pr_X \big) \subseteq \liealg{g} = \liealg{so}(n) \times \liealg{so}(k).
\end{equation}
By~\cite[Eq. (14)]{hueper.markina.leite:2021}, the stabilizer
subgroup $H \subseteq \liegroup{O}(n) \times \liegroup{O}(k)$
is isomorphic to the
Lie group $\liegroup{O}(n - k) \times \liegroup{O}(k)$.

Next we recall the definition of the so-called $\alpha$-metrics
from \cite{hueper.markina.leite:2021}.
To this end, a bi-invariant metric on
$\liealg{so}(n) \times \liealg{so}(k)$ is introduced
following \cite[Def. 3.1]{hueper.markina.leite:2021}.
Define for $0 \neq \alpha \in \field{R}$
\begin{equation}
	\label{equation:alpha_metric_scalar_product_on_group_definition}
	\begin{split}
		\langle \cdot, \cdot \rangle^{\alpha} \colon ( \liealg{so}(n) \times \liealg{so}(k) ) \times ( \liealg{so}(n) \times \liealg{so}(k)) 
		&\to
		\field{R}, \\
		\big( (\Omega_1, \Omega_2), (\eta_1, \eta_2) \big)
		&\mapsto
		- \tr(\Omega_1 \Omega_2) - \tfrac{1}{\alpha} \tr(\eta_1 \eta_2) .
	\end{split}
\end{equation}
By\cite[Prop. 2]{hueper.markina.leite:2021}, the subspace
$\liealg{h} \subseteq \liealg{so}(n) \times \liealg{so}(k)$
defined in \eqref{equation:stiefel_alpha_lie_algebra_h_at_X_defintion}
is non-degenerated iff $\alpha \neq -1$ holds.
In this case, we write $\liealg{m} = \liealg{h}^{\perp}$ and 
\begin{equation}
	\liealg{m} \oplus \liealg{h} = \liealg{g} = \liealg{so}(n) \times \liealg{so}(k)
\end{equation}
is fulfilled.
Next we reformulate \cite[Def. 3.3]{hueper.markina.leite:2021}.

\begin{definition}
	\label{definition:stiefel_alpha_metrics}
	Let $\alpha \in \field{R} \setminus \{0, -1\}$
	and let $\liegroup{O}(n) \times \liegroup{O}(k)$
	be equipped with the bi-invariant metric defined
	by the scalar product $\langle \cdot, \cdot \rangle^{\alpha}$
	from~\eqref{equation:alpha_metric_scalar_product_on_group_definition}.
	The metric on $\Stiefel{n}{k}$ defined
	by requiring that the map
	$\pr_X \colon \liegroup{O}(n) \times \liegroup{O}(k) \to \Stiefel{n}{k}$
	from~\eqref{equation:definition_stiefel_pr_X}
	is a pseudo-Riemannian submersion
	is called $\alpha$-metric.
\end{definition}
The Stiefel manifold equipped with an $\alpha$-metric is
a naturally reductive homogeneous space.

\begin{lemma}
	\label{lemma:stiefel_alpha_metric_naturally_reductive}
	Let $\alpha \in \field{R} \setminus \{-1, 0\}$.
	Then $(\liegroup{O}(n) \times \liegroup{O}(k)) / H \cong \Stiefel{n}{k}$
	equipped with an $\alpha$-metric from Definition~\ref{definition:stiefel_alpha_metrics}
	is a naturally reductive homogeneous space.
	\begin{proof}
		Obviously,
		the scalar product $\langle \cdot, \cdot \rangle^{\alpha}$
		on $\liealg{so}(n) \times \liealg{so}(k)$ 
		from Definition~\ref{definition:stiefel_alpha_metrics}
		is $\Ad(\liegroup{O}(n) \times \liegroup{O}(k))$-invariant
		for $\alpha \in \field{R} \setminus \{0\}$.
		In addition, the subspace
		$\liealg{h} = \ker(T_{(I_n, I_k)}\pr_X) \subseteq \liealg{so}(n) \times \liealg{so}(k)$
		is non-degenerated
		for $\alpha \in \field{R} \setminus \{0, - 1\}$
		by~\cite[Prop. 2]{hueper.markina.leite:2021}.
		Thus
		Lemma~\ref{lemma:normal_naturally_reductive_is_naturally_reductive}
		yields the desired result.
	\end{proof}
\end{lemma}

In the sequel, an explicit expression for the orthogonal projection
$\pr_{\liealg{m}} \colon \liealg{so}(n) \times \liealg{so}(k) \to \liealg{m}$
with respect to the scalar product $\langle \cdot, \cdot \rangle^{\alpha}$
is needed.
Therefore we state the next lemma which is taken
from~\cite[Lem. 3.2]{hueper.markina.leite:2021}.

\begin{lemma}
	\label{lemma:stiefel_orthogonal_projection_on_m}
	Let $\alpha \in \field{R} \setminus \{-1, 0\}$ and let $X \in \Stiefel{n}{k}$.
	The orthogonal projection
	\begin{equation}
		\pr_{\liealg{m}} 
		\colon \liealg{so}(n) \times \liealg{so}(k)
		\ni (\Omega, \eta) \mapsto (\Omega^{\perp_X}, \eta^{\perp_X})
		\in
		\liealg{m} \subseteq \liealg{so}(n) \times \liealg{so}(k)
	\end{equation}
	is given by
	\begin{equation}
		\begin{split}
			\Omega^{\perp_X} 
			&=
			X X^{\top} \Omega 
			+ \Omega X X^{\top} 
			- \tfrac{2 \alpha + 1}{\alpha +1} X X^{\top} \Omega X X^{\top} 
			- \tfrac{1}{\alpha + 1} X \eta X^{\top} ,
			\\
			\eta^{\perp_X} 
			&=
			\tfrac{\alpha}{\alpha + 1} \big( \eta - X^{\top} \Omega X  \big) .
		\end{split}
	\end{equation}
	\begin{proof}
		This is just a reformulation of \cite[Lem. 3.2]{hueper.markina.leite:2021}.
	\end{proof}
\end{lemma}
Furthermore, the following lemma is a trivial reformulation
of~\cite[Prop. 3]{hueper.markina.leite:2021}.

\begin{lemma}
	\label{lemma:stiefel_identification_tangent_space_m}
	Let $\alpha \in \field{R} \setminus \{-1, 0\}$
	and let $X \in \Stiefel{n}{k}$.
	The map 
	\begin{equation}
		\big(T_{(I_n, I_k)} (\iota_X \circ \pr) \at{\liealg{m}} \big)^{-1} 
		\colon T_X \Stiefel{n}{k}
		\ni V \mapsto (\Omega(V)^{\perp_X}, \eta(V)^{\perp_X}) 
		\in \liealg{m} \subseteq  \liealg{so}(n) \times \liealg{so}(k)
	\end{equation}
	is given by 
	\begin{equation}
		\begin{split}
			\Omega(V)^{\perp_X} 
			&=
			V X^{\top} 
			- X V^{\top} 
			+ \tfrac{2 \alpha + 1}{\alpha + 1} X V^{\top} X X^{\top} ,
			\\
			\eta(V)^{\perp_X}
			&= 
			-\tfrac{\alpha}{\alpha + 1} X^{\top} V 
		\end{split}
	\end{equation}
	for all $V \in T_X \Stiefel{n}{k}$.
	\begin{proof}
		This is a consequent of \cite[Prop. 3]{hueper.markina.leite:2021}.
	\end{proof}
\end{lemma}

\subsubsection{Intrinsic Rolling}
\label{subsec:stiefel_intrinsic_rolling}

We now determine intrinsic 
rollings of the Stiefel manifold equipped with an $\alpha$-metric
over one of its tangent spaces.

By Lemma~\ref{lemma:configuration:space_intrinsic_rolling},
the configuration space for rolling $T_X \Stiefel{n}{k} \cong \liealg{m}$
over $\Stiefel{n}{k}$ intrinsically is given by the fiber bundle 
\begin{equation}
	\prQ \colon \Qspace 
	=
	\liealg{m} \times 
	\big((\liegroup{O}(n) \times \liegroup{O}(k) ) 
	\times_H \liegroup{O}(\liealg{m}) \big) 
	\to \liealg{m} \times (\liegroup{O}(n) \times \liegroup{O}(k) ) / H ,
\end{equation}
where $H = \Stab(X) \subseteq \liegroup{O}(n) \times \liegroup{O}(k) = G$.
By identifying $T_X \Stiefel{n}{k} \cong \liealg{m}$
via the linear isometry $T_{(I_n, I_k)} \liealg{m} \to T_X \Stiefel{n}{k}$
from Lemma~\ref{lemma:stiefel_identification_tangent_space_m}
and $\Stiefel{n}{k} \cong (\liegroup{O}(n) \times \liegroup{O}(k) )/ H$
via the $(\liegroup{O}(n) \times \liegroup{O}(k))$-equivariant isometry
$\iota_X \colon (\liegroup{O}(n) \times \liegroup{O}(k)) / H 
\to \Stiefel{n}{k}$,
we obtain the following proposition describing intrinsic rollings
of $T_X \Stiefel{n}{k}$ over $\Stiefel{n}{k}$.

\begin{proposition}
	\label{proposition:stiefel_manifold_intrinsic_rolling}
	Let $\Stiefel{n}{k}$ be equipped with an $\alpha$-metric
	for $\alpha \in \field{R} \setminus \{-1, 0\}$ and let
	\begin{equation}
		V \colon I \to T_X \Stiefel{n}{k}, 
		\quad t \mapsto V(t)
	\end{equation}
	be a given rolling curve.
	Denote by $v \colon I \to \liealg{m}$
	the corresponding curve in $\liealg{m}$ given by 
	\begin{equation}
		\begin{split}
			v(t)
			&=
			\big(T_{I_{n, k}} (\iota_X \circ \pr\at{\liealg{m}}) \big)^{-1} V(t) \\
			&= \big(	V(t) X^{\top} 
			- X V(t)^{\top} 
			+ \tfrac{2 \alpha + 1}{\alpha + 1} X V(t)^{\top} X X^{\top}, 	
			- \tfrac{\alpha}{\alpha + 1} X^{\top} V(t)\big) 
		\end{split}
	\end{equation}
	for $t \in I$.
	Then the kinematic equation for the intrinsic rolling
	of $\Stiefel{n}{k}$ over $\liealg{m} \cong T_X \Stiefel{n}{k}$
	with respect to $\nablaLC$ defined by the $\alpha$-metric
	along $v \colon I \to \liealg{m}$ is given by
	\begin{equation}
		\label{equation:proposition_stiefel_manifold_intrinsic_rolling_kinematic_equations}
		\begin{split}
			\dot{v}(t) 
			& =
			u(t), \\
			\dot{S}(t) 
			&=
			- \tfrac{1}{2} \pr_{\liealg{m}} \circ \ad_{S(t) u(t)} \circ S(t),  \\
			\dot{g}(t) 
			&=
			\big( T_e \ell_{g(t)} \circ  S(t) \big) u(t) ,
		\end{split}
	\end{equation}
	where
	$\pr_{\liealg{m}} \colon  \liealg{so}(n) \times \liealg{so}(k) \to \liealg{m}$ 
	is explicitly given by
	Lemma~\ref{lemma:stiefel_orthogonal_projection_on_m}.
	Let
	$\overline{q} \colon I \ni t \mapsto (v(t), g(t), S(t)) \in 
	\QspaceLift = \liealg{m} \times (\liegroup{O}(n) \times \liegroup{O}(k)) \times \liegroup{O}(\liealg{m})$
	be a curve
	satisfying~\eqref{equation:proposition_stiefel_manifold_intrinsic_rolling_kinematic_equations}.
	Then
	\begin{equation}
		q \colon I \to Q, 
		\quad
		t \mapsto q(t) = (\prQLift \circ \overline{q})(t)
		= (v(t), [g(t), S(t)])
	\end{equation}
	is an intrinsic rolling of $T_X \Stiefel{n}{k} \cong \liealg{m}$
	over $\Stiefel{n}{k}$ with respect to the given $\alpha$-metric
	along  the rolling curve $v$.
	The development curve
	$I  \ni t \mapsto \pr(g(t)) = (R(t), \theta(t)) \in \liegroup{O}(n) \times \liegroup{O}(k)$
	is mapped by the embedding
	$\iota_X \colon (\liegroup{O}(n) \times \liegroup{O}(k)) / H \to \matR{n}{k}$
	to the curve
	\begin{equation}
		\gamma \colon I \to \Stiefel{n}{k},
		\quad 
		t \mapsto \gamma(t) = (\iota_X \circ \pr)(g(t))
		=
		\pr_X(g(t)) 
		= R(t) X \theta(t)^{\top}.
	\end{equation}
	\begin{proof}
		Since $\Stiefel{n}{k}$ equipped with an $\alpha$-metric is a
		naturally reductive homogeneous space by
		Lemma~\ref{lemma:stiefel_alpha_metric_naturally_reductive},
		this is a direct consequence of
		Corollary~\ref{corollary:intrinsic_rolling_kinematic_equation_naturally_redudctive_homogeneous_spaces}
		combined with
		Lemma~\ref{lemma:stiefel_orthogonal_projection_on_m} and
		Lemma~\ref{lemma:stiefel_identification_tangent_space_m}.
	\end{proof}
\end{proposition}

Next we consider the intrinsic rolling of the Stiefel manifolds along curves of a special form
by using Section~\ref{subsubsec:rolling_along_special_curves}.
This yields the next remark.
\begin{remark}
	\label{remark:stiefel_rolling_along_special_curve}
	Let $\xi = (\xi_1, \xi_2) \in \liealg{so}(n) \times \liealg{so}(k)$. 
	Then 
	\begin{equation}
		q \colon I \to \liealg{m} \times (G \times_H \liegroup{O}(\liealg{m})), 
		\quad 
		t \mapsto (v(t),  [g(t), S(t)]) ,
	\end{equation}	
	where 
	\begin{equation}
		\begin{split}
			g(t) 
			&=
			\exp(t \xi) \exp(-t \xi_{\liealg{h}}) , \\
			S(t)
			&= \Ad_{\exp(t \xi_{\liealg{h}})} \circ 
			\exp\Big(- t \pr_{\liealg{m}} \circ \ad_{\xi_{\liealg{h}} + \tfrac{1}{2} \xi_{\liealg{m}}}\Big) , \\
			v(t)
			&=
			\int_{0}^{t}
			\exp\Big(s  \Big(\pr_{\liealg{m}} \circ 
			\ad_{\xi_{\liealg{h}} + \tfrac{1}{2} \xi_{\liealg{m}}} \Big) \Big)(\xi_{\liealg{m}}) \D s
		\end{split}
	\end{equation}
	is an intrinsic rolling of $\liealg{m}$ along the rolling curve
	$v \colon I \to \liealg{m}$ with development curve
	$\gamma(t) = \pr(g(t)) = \pr(\exp(t \xi))$.
	Identifying
	$\Stiefel{n}{k} \cong (\liegroup{O}(n) \times \liegroup{O}(k)) / H$
	with the embedded submanifold $\Stiefel{n}{k} \subseteq \matR{n}{k}$
	via
	$\iota_X \colon (\liegroup{O}(n) \times \liegroup{O}(k)) / H 
	\to \Stiefel{n}{k} \subseteq \matR{n}{k}$,
	the development curve is given by 
	\begin{equation}
		\gamma(t) = \e^{t \xi_1} X e^{- t \xi_2}
	\end{equation}
	and the rolling curve reads 
	\begin{equation}
		V(t) = T_X (\iota_X \circ \pr) v(t) = v_1(t) X  - X v_2(t),
	\end{equation}
	where we write $v(t) = (v_1(t), v_2(t)) \in \liealg{m}$ with
	$v_1(t) \in \liealg{so}(n)$ and $v_2(t) \in \liealg{so}(k)$
	for all $t \in I$.
\end{remark}

\section{Conclusion}

In this text, we investigated intrinsic rollings of
reductive homogeneous spaces equipped with invariant
covariant derivatives.
As preparation,
we considered
frame bundles of vector bundles associated to principal fiber bundles
in detail.
Afterwards, using an abstract definition of intrinsic rolling
as starting point,
we investigated rollings
of $\liealg{m}$ over the reductive homogeneous spaces $G / H$
with respect to an invariant covariant derivative $\nablaAlpha$.
For a given control curve, we obtained the so-called
kinematic equation which is a time-variant explicit ODE on a Lie group,
whose solutions describe rollings of
$\liealg{m}$ over $G / H$.
Moreover, for the case, where the development curve is
the projection of a one-parameter subgroup,
we provided explicit solutions
of the kinematic equation
describing intrinsic rollings of 
$\liealg{m}$ over $G / H$
with respect to the canonical covariant derivative of first kind
and second kind, respectively.
As examples, we discussed intrinsic rollings of Lie groups
and Stiefel manifolds.

\section*{Acknowledgments}
The author would like to thank Knut H{\"u}per, Irina Markina
and F{\'a}tima Silva Leite
for some fruitful discussions about rolling manifolds.
This work has been supported by the 
German Federal Ministry of Education and Research
(BMBF-Projekt 05M20WWA: Verbundprojekt 05M2020 - DyCA).

\end{document}